\documentclass[a4paper,11pt]{article}
\usepackage[french]{babel}
\usepackage{latexsym}
\usepackage{frbib}
\usepackage{amssymb}
\usepackage{amsfonts}
\usepackage{amssymb}
\usepackage{amsmath}
\usepackage{amsthm}
\usepackage{mathrsfs}
\usepackage{graphicx}
\usepackage{psfrag}
\usepackage{url}
\usepackage{caption}
\usepackage{enumerate}
\usepackage[ps2pdf]{hyperref}

\input cyracc.def

\newtheorem{teor}{Th\'eor\`eme}
\newtheorem{lema}[teor]{Lemme}
\newtheorem{prop}[teor]{Proposition}
\newtheorem{coro}[teor]{Corollaire}
\newcommand{\summ}{\displaystyle\sum}
\newcommand{\prodd}{\displaystyle\prod}
\newcommand{\Frac}{\displaystyle\frac}
\newcommand{\eps}{\varepsilon}
\newcommand{\re}{\mathbb{R}}
\newcommand{\ce}{\mathbb{C}}
\newcommand{\ovl}{\overline}
\newcommand{\ds}{\displaystyle}

\newcommand{\mfrac}[2]{\mbox{$\frac{#1}{#2}$}}

\newcommand{\ment}[1]{\lceil #1 \rceil}
\newcommand{\gment}[1]{\left\lceil #1 \right\rceil}

\DeclareMathOperator{\mmod}{mod}

\DeclareMathOperator{\PSL}{PSL}
\DeclareMathOperator{\Res}{Res}

\newcommand{\Og}{\mathcal{O}}

\newcommand{\lqqd}{\hfill{$\blacksquare$}}

\def\alfabetico{\renewcommand{\theenumi}{\alph{enumi}}
\renewcommand{\labelenumi}{\theenumi)}}

\def\romano{\renewcommand{\theenumi}{\roman{enumi}}
\renewcommand{\labelenumi}{(\theenumi)}}

\newcommand{\bsh}{\backslash}

\newcommand{\itemref}[1]{\eqref{#1}}
\newenvironment{enumeratei}{\begin{enumerate}%
[\upshape (i)]}{\end{enumerate}}

\newenvironment{enumerateu}{\begin{enumerate}%
[\upshape 1)]}{\end{enumerate}}

\newcommand{\dds}[1]{\frac{\partial #1}{\partial s}}

\newcommand{\tq}{\mid}

\title{Majoration du nombre de z\'eros d'une fonction m\'eromorphe en dehors
d'une droite verticale et applications}
\author{Oswaldo Vel\'asquez Casta\~n\'on
\thanks{Institut de Math\'ematiques de Bordeaux, UMR 5251, Universit\'e Bordeaux 1, 351 Cours de la Lib\'eration,
33405 Talence Cedex, France. Courriel: \texttt{Oswaldo.Velasquez@math.u-bordeaux1.fr}}
}

\date{D\'ecembre 2007}

\begin{document}

\maketitle

\begin{abstract}
On \'etudie la r\'epartition des z\'eros des fonctions de la forme
$f(s) = h(s) \pm h(2a-s)$, o\`u $h(s)$ est une fonction m\'eromorphe r\'eelle sur l'axe r\'eel, $a$ un nombre r\'eel.
Un de nos r\'esultats \'etablit des conditions
suffisantes pour que tous les z\'eros de $f(s)$, sauf un nombre fini,
se trouvent sur la droite $\Re s = a$, appel\'ee {\it droite critique} pour la fonction $f(s)$, et qu'ils soient simples,
pourvu que tous les z\'eros de $h(s)$, sauf un nombre fini, soient dans le demi-plan $\Re s < a$.
Ce r\'esultat peut \^etre vu comme une g\'en\'eralisation de la condition n\'ecessaire de stabilit\'e de la fonction $h(s)$,
dans le th\'eor\`eme d'Hermite-Biehler.
On applique ces r\'esultats \`a l'\'etude de translat\'ees de la fonction
z\^eta de Riemann et de fonctions $L$, et des int\'egrales de s\'eries d'Eisenstein, entre autres.
\end{abstract}

\renewcommand{\abstractname}{Abstract}
\begin{abstract}
We study the distribution of the zeros of functions of the form $f(s) = h(s) \pm h(2a-s)$, where $h(s)$ is a
meromorphic function, real on the real line, $a$ a real number.
One of our results establishes sufficient conditions under which
all but finitely many of the zeros of $f(s)$
lie on the line $\Re s = a$, called the {\it critical line} for the function $f(s)$,
and be simple, given that all but finitely many of the zeros of $h(s)$
lie on the half-plane $\Re s < a$.
This results can be regarded as a generalization of the necessary condition of stability for the
function $h(s)$, in the Hermite-Biehler theorem.
We apply this results to the study of translations of the Riemann Zeta Function and $L$ functions, and integrals of Eisenstein Series,
among others.
\end{abstract}


\section{Motivation et notations} \label{notaciones}

Le probl\`eme de montrer que les z\'eros de certaines fonctions sont align\'es
est li\'e \`a l'hypoth\`ese de Riemann, l'un des plus grands probl\`emes ouverts des
math\'ematiques.
On rappelle quelques faits sur ce probl\`eme, cf. \cite[\S\S 1.1, 2.1, 3.3]{ztitchmarsh}.
La fonction z\^eta de Riemann, est la fonction complexe d\'efinie par
$$ \zeta(s) = \summ_{n=1}^\infty \Frac{1}{n^s}, \quad \Re s > 1 $$
et \'etendue \`a tout le plan complexe par prolongement analytique,
sauf au point $s = 1$, o\`u $\zeta(s)$ a un p\^ole simple de
r\'esidu $1$. Elle satisfait \`a l'\'equation fonctionnelle
$$ \zeta^*(s) = \zeta^*(1-s) $$
pour $s \in \ce \setminus \{0, 1\}$, o\`u $\zeta^*(s)$ est la fonction z\^eta compl\'et\'ee, d\'efinie par
$\zeta^*(s) = \pi^{-s/2} \Gamma(s/2) \zeta(s)$, o\`u $\Gamma(s)$ est la fonction gamma. Cette nouvelle fonction a deux p\^oles simples, en $s=0$, $s=1$,
de r\'esidus $1$. On peut r\'e\'ecrire l'\'equation fonctionnelle sous la forme
$\xi(s) = \xi(1-s)$, o\`u $\xi(s) = \frac{1}{2} s (s-1) \zeta^*(s)$ est une fonction enti\`ere de genre $1$.

On voit facilement \`a partir du produit eul\'erien de $\zeta(s)$ et de l'\'equation fonctionnelle que les fonctions
$\zeta^*(s)$ et $\xi(s)$ ont tous leurs z\'eros dans la bande $0 \leq \Re s \leq 1$. Pour la fonction z\^eta $\zeta(s)$
les p\^oles de $s \Gamma(s/2)$ sont aussi des z\'eros; ce sont les entiers n\'egatifs pairs (appel\'es z\'eros triviaux de la
fonction z\^eta). La non annulation des fonctions pr\'ec\'edentes sur la droite $\Re s = 1$ (donc aussi sur $\Re s=0$)
a \'et\'e d\'emontr\'ee par Hadamard et de la Vall\'ee-Poussin
(th\'eor\`eme des nombres premiers). On peut \'etendre cette r\'egion de non-annulation au del\`a de $\sigma \geq 1$;
un des ces r\'esultats est que pour une constante $A>0$, $\sigma \geq 1- A/\log|\tau|$, $|\tau|$ assez grand
\cite[\S 3.11]{ztitchmarsh}
\begin{equation} \label{nonannul}
\Frac{1}{\zeta(s)} = O \bigl(\log |\tau|\bigr).
\end{equation}
On peut trouver un meilleur r\'esultat en \cite[\S 6.19]{ztitchmarsh}.

Les z\'eros non triviaux de la fonction z\^eta de Riemann sont donc dans la bande $0 < \Re s < 1$,
et r\'epartis sym\'etriquement par rapport \`a la droite $\Re s = 1/2$, appel\'ee droite critique. L'hypoth\`ese de Riemann est
l'assertion que tous ces z\'eros sont sur la droite $\Re s = 1/2$.
On conjecture de plus que ces z\'eros sont simples.

Soit $a \in \re$, $h(s)$ une fonction m\'eromorphe dans le plan complexe, r\'eelle sur l'axe r\'eel.
On d\'efinit la fonction
$$ f(s) = h(s) \pm h(2a-s). $$
Cette fonction satisfait \`a l'\'equation $f(2a-s) = \pm f(s)$.
En particulier, les z\'eros de $f(s)$ sont sym\'etriques par rapport \`a la droite $\Re s = a$,
appel\'ee {\it droite critique}. Un lien entre la r\'epartition des z\'eros de $h(s)$ et $f(s)$
est \'etabli par le th\'eor\`eme d'Hermite-Biehler (cf. \cite[Part III, Lecture 27]{levin}, \cite{cebotarev}):
sous des bonnes conditions de croissance sur $h(s)$, l'alignement et la simplicit\'e
des z\'eros des fonctions $h(s)\pm h(2a-s)$ sont \'equivalents
\`a l'absence de z\'eros de $h(s)$ dans le demi-plan $\Re s \geq a$.
Cette derni\`ere condition, connue par le nom de {\it stabilit\'e} de $h(s)$,
appara\^it dans l'\'etude du comportement asymptotique des
solutions de certaines \'equations diff\'erentielles (cf. \cite{BC1963}).
Ceci est le cas o\`u $a=0$; les fonctions li\'ees \`a l'hypoth\`ese de Riemann
correspondent au cas $a=1/2$ (l'\'etude des probl\`emes correspondants \`a deux valeurs
de $a$  sont \'equivalents).

Dans ce travail, on \'etude les cons\'equences d'une relaxation de la condition de stabilit\'e,
par l'introduction de z\'eros de $h(s)$ dans le demi-plan $\Re s\geq a$.
Sous certaines conditions, on obtient une borne effective
pour le nombre de z\'eros de $f(s)$ en dehors de la droite critique,
ainsi que des informations sur la simplicit\'e de ces z\'eros.
L'id\'ee consiste \`a compter les z\'eros de $f(s)$ dans une bande
$0 < \Im(s) < T$, \`a l'aide de la m\'ethode standard (cf. \cite[\S\S 9.3, 9.4]{ztitchmarsh}),
puis comparer l'estimation trouv\'ee avec celle du nombre de z\'eros de $f(s)$ sur la droite critique.
Cette id\'ee est inspir\'ee de l'article de P. R. Taylor \cite{taylor}, \'el\`eve de Titchmarsh. \\

Fixons quelques notations. Dor\'enavant, pour $s \in \ce$, on note $\sigma = \Re s$  et $\tau = \Im s$.
On suppose qu'il existe $\sigma_0 > a$ tel que $f(s) \neq 0$ pour $\sigma \geq \sigma_0$. Par l'\'equation fonctionnelle,
$f(s) \neq 0$ lorsque $2a-\sigma \geq \sigma_0$, c'est-\`a-dire, lorsque $\sigma \leq 2a-\sigma_0$. Donc, tous les z\'eros de $f(s)$
sont dans la bande $|\sigma-a| < \sigma_0-a$.
\'Etant donn\'e $T > 0$, on note
$$ N(T) = \# \{s \in \ce \tq f(s)= 0, 0 < \tau < T \}, $$
le nombre des z\'eros de $f(s)$ tel que $0 < \tau < T$, et
$$ N_0(T) = \# \{s \in \ce \tq f(s)= 0, s= a+i\tau, 0 < \tau < T\}$$
le nombre de ces z\'eros sur la droite critique $\sigma = a$ (on tient compte de la multiplicit\'e des z\'eros
dans $N(T)$ et $N_0(T)$). On note $N_0'(T)$ le
nombre des z\'eros critiques, sans tenir compte de leur multiplicit\'e. Il est \'evident que
$0 \leq N_0'(T) \leq N_0(T) \leq N(T)$.

Si l'on pose
$$L(T) = \# \{s \in \ce \tq f(s) = 0, \sigma > a, 0 < \tau < T\}, $$
l'application $s \mapsto \overline{2a-s}$ \'etablit une
correspondance entre les z\'eros de $f(s)$ tels que $\sigma > a$, $0 < \tau
< T$ et ceux tels que $\sigma < a$, $0 < \tau < T$, donc
$$ N(T)-N_0(T) = \# \{s \in \ce \tq f(s) = 0, \sigma \neq a, 0 <\tau<T \} = 2 L(T). $$
Cela montre que $N(T)-N_0(T)$ est un nombre pair positif ou nul. On a de plus
$$ 0 \leq N(T)-N_0(T) \leq N(T)-N_0'(T). $$
C'est cette derni\`ere quantit\'e celle que l'on cherche \`a estimer \`a partir
d'informations sur $h(s)$.

Cet article est organis\'e de la fa\c{c}on suivante. Dans un premier temps, on consid\`ere le cas o\`u
il y a une quantit\'e finie de z\'eros de $h(s)$ dans le demi-plan $\sigma \geq a$.
Au \S \ref{seccerosrecta}, on obtient une premi\`ere minoration pour le nombre de z\'eros sur la droite critique $\sigma = a$
de la fonction $f(s)$. Au \S \ref{princcc}, et sous certaines conditions de croissance sur $h(s)$,
on obtient le r\'esultat central de l'article, le th\'eor\`eme \ref{teoalin}, qui fournit
une borne pour le nombre de z\'eros de $f(s)$ en dehors de la droite critique $\sigma=a$, en fonction
du nombre de z\'eros de la fonction $h(s)$ dans le demi-plan $\sigma \geq a$.
En cons\'equence du th\'eor\`eme, on montre au \S \ref{repglobale} que la borne inf\'erieure donn\'ee au \S \ref{seccerosrecta}
donne une bonne estimation pour le nombre de z\'eros de $f(s)$, et on obtient la repartition
globale des z\'eros de $f(s)$. De plus, dans certains cas
on obtient une am\'elioration du th\'eor\`eme \ref{teoalin}.
Le lien entre nos r\'esultats et la th\'eorie de la stabilit\'e des polyn\^omes,
notamment le th\'eor\`eme d'Hermite-Biehler, est donn\'ee au \S \ref{policaso}.
Un exemple d'application des r\'esultats, motiv\'e par l'\'etude
d'une premi\`ere approximation \`a la fonction z\^eta de Riemann, est donn\'ee au \S \ref{zeta2}.
La r\'epartition des z\'eros des sommes et des diff\'erences de translat\'ees
de la fonction z\^eta compl\'et\'ee commence au \S \ref{aquiestimation}. La n\'ecessit\'e de
travailler avec une quantit\'e infinie de z\'eros \`a droite de la droite critique
nous am\`ene \`a la premi\`ere extension du th\'eor\`eme \ref{teoalin}.
On g\'en\'eralise les r\'esultats du \S \ref{seccerosrecta} au \S \ref{inficritica}, puis on \'etablit
l'extension souhait\'ee au \S \ref{infidensidad}.
On \'etudie des fonctions provenantes de certaines int\'egrales
de s\'eries d'Eisenstein sur la surface modulaire aux \S \ref{lweng}.
Une nouvelle extension du th\'eor\`eme \ref{teoalin} motiv\'ee au \S \ref{lweng},
nous permettra d'obtenir des r\'esultats sous des conditions de croissance affaiblies pour $h(s)$, au \S \ref{faibledd}.
On utilise aussi nos r\'esultats pour simplifier des r\'esultats r\'ecents sur
la r\'epartition des z\'eros d'approximations des fonctions z\^eta d'Epstein, au \S \ref{nonefectepstein}.
Finalement, on g\'en\'eralise au \S \ref{sanssymetrie} nos r\'esultats au cas o\`u $h(s)$
n'a pas la sym\'etrie r\'eelle, et on applique ces g\'en\'eralisations \`a l'\'etude
de combinaisons de translat\'ees de s\'eries $L$ de Dirichlet. Une annexe
\`a la fin du travail lie nos r\'esultats \`a d'autres obtenus sous un point de vue
leg\`erement diff\'erent.

Les calculs num\'eriques et les graphes des fonctions consid\'er\'ees au \S \ref{aquiestimation} ont \'et\'e r\'ealis\'es
\`a l'aide du syst\`eme PARI/GP \cite{PARI2}.

\section{Z\'eros sur la droite critique} \label{seccerosrecta}

Soit $h(s)$ holomorphe et sans z\'ero sur $\sigma=a$. \`A partir d'une valeur quelconque de $\arg h(a)$,
on d\'efinit la fonction
$$ \varphi(\tau) = \arg h(a+i \tau)$$
par variation continue de l'argument de $h(s)$ le long du segment qui lie $a$ et $a+i\tau$. Une telle fonction
est appel\'ee {\it fonction de phase} de $h(s)$.

Commen\c{c}ons par un lemme simple, dont on fournit la preuve en vue d'une modification ult\'erieure.

\begin{lema} \label{compcreciente}
Soit $g: [0,\infty) \to \re$ une fonction continue telle que $g(0)
\leq 0$. Posons
$$ \tilde{g}(T) = \# \{ x \tq 0 \leq x < T, \, g(x) \equiv 0\mmod 1\}.$$
Alors, pour $T > 0$,
$$ \tilde{g}(T) \geq g(T). $$
\end{lema}

{\bf D\'emonstration.} Si $\tilde{g}(T)=\infty$, le r\'esultat est \'evident. Sinon, soient $x_1 <  \dots < x_n  $ les
nombres $x_i \in [0,T[$ tels que $g(x_i) \equiv 0 \mmod 1$, par
d\'efinition $n=\tilde{g}(T)$. Par continuit\'e, $g(x_1) \leq 0$,
$\bigl|g(x_{i})-g(x_{i-1})\bigr| \leq 1$ pour tout $i=2,\dots,n$ et
$\bigl|g(T)-g(x_n)\bigr|\leq 1$. On obtient
alors
$$
g(T)  \leq  g(T)-g(x_1) = g(T)-g(x_n) + \summ_{i=2}^n \bigl( g(x_i)-g(x_{i-1}) \bigr) \leq n
= \tilde{g}(T).
$$
Cela conclut la preuve. \hfill{$\blacksquare$\\

En fait, on a montr\'e ci-dessus que $\ment{g(T)} \leq \tilde{g}(T)$.

Il faut remarquer que si l'on cherche \`a en d\'eduire une bonne borne inf\'erieure de $\tilde{g}(T)$,
la fonction $g(T)$ doit \^etre presque croissante. Mais ce lemme nous fournit au moins une premi\`ere
estimation pour le nombre de z\'eros de $f(s)$ sur la droite critique, en utilisant une fonction de phase de $h(s)$.

\begin{lema} \label{rectacritica} Soit $a \in \re$ fix\'e, $h(s)$ une fonction m\'eromorphe dans le plan complexe,
r\'eelle sur l'axe r\'eel, sans z\'eros ni p\^oles sur la droite critique $\sigma = a$, avec $h(a)> 0$.
On consid\`ere $\arg h(a+i\tau)$, la d\'etermination continue de l'argument de $h(s)$ sur la droite critique $\sigma=a$ telle que $\arg h(a) = 0$.
Soit $f(s) = h(s) \pm h(2a-s)$. Le nombre de z\'eros de $f(s)$ tels que $0< \tau < T$, sans compter
les multiplicit\'es, sur la droite $\sigma = a$ est minor\'e par
\begin{equation} \label{dsuperior}
N_0'(T) \geq \frac{1}{\pi}\arg h(a+iT) -u_{\pm},
\end{equation}
pour $T > 0$, o\`u $ u_+ = \frac{1}{2}$, $u_- = 1$. La multiplicit\'e de $a$ comme z\'ero de $f(s)$ est
$n_{f,a} = 0$ dans le cas $+$, $n_{f,a} \geq 1$ et impaire dans le cas $-$.
\end{lema}

{\bf D\'emonstration.}
Si $s= a+i\tau$,
$$
\begin{array}{rcl}
f(a+i\tau) &=& h(a+i\tau) \pm h(a-i\tau) = h(a+i\tau) \pm \ovl{h(a+i\tau)} \\
    &=& \bigl|h(a+i\tau)\bigr| e^{i \arg h(a+i\tau)} \pm \bigl|h(a+i\tau)\bigr| e^{-i \arg h(a+i\tau)}.
\end{array}
$$
On consid\`ere deux cas:
\begin{enumerate}
\romano
\item Lorsque $f(s) = h(s)-h(2a-s)$,
$$ f(a+i\tau) = 2 i \bigl|h(a+i\tau)\bigr| \sin \bigl(\arg h(a+i\tau)\bigr) $$
et $h(a+i\tau) \neq 0$, donc $f(a+i\tau) = 0$ si et seulement si
$ \arg h(a+i\tau) \equiv 0 \mmod \pi $.
Donc
$$ N_0'(T) = \# \{\tau \tq  0 < \tau < T,  \mfrac{1}{\pi}\arg h(a+i\tau) \equiv 0 \mmod 1 \}. $$
Par le lemme \ref{compcreciente}, et comme $\arg h(a) = 0$,
$$ N_0'(T) =
\# \{ \tau \tq 0 \leq \tau < T, \mfrac{1}{\pi}\arg h(a+i\tau) \equiv 0 \mmod 1 \}-1
        \geq \mfrac{1}{\pi}\arg h(a+iT) - 1. $$
Dans ce cas $f(a) = 0$, et comme $f(\sigma) = -f(2a-\sigma)$, la
fonction $f(s)$ change de signe sur la droite r\'eelle au point
$s=a$, donc $n_{f,a} \geq 1$ et $n_{f,a}$ est impaire.

\item Lorsque $f(s) = h(s)+h(2a-s)$,
$$ f(a+i\tau) = 2 \bigl|h(a+i\tau)\bigr| \cos \bigl(\arg h(a+i\tau)\bigr), $$
donc $f(a+i\tau) = 0$ si et seulement si $ \arg h(a+i\tau) \equiv \mfrac{\pi}{2} \mmod \pi$.
D'apr\`es le lemme \ref{compcreciente}
$$ N_0'(T) = \# \{\tau \tq  0 < \tau < T,   \mfrac{1}{\pi}\arg h(a+i\tau) \equiv \mfrac{1}{2} \mmod 1 \}
     \geq  \mfrac{1}{\pi}\arg h(a+iT) - \mfrac{1}{2}.$$

On a aussi $f(a) = 2 h(a) \neq 0$, donc $n_{f,a} = 0$.
\hfill{$\blacksquare$}\\
\end{enumerate}

Si la fonction de phase $\varphi(\tau)=\arg h(a+i\tau)$ est croissante, l'estimation
dans (\ref{dsuperior}) est optimale, et
$$N_0'(T) = \gment{\mfrac{1}{\pi}\arg h(a+iT) -u_{\pm}}. $$
On peut aussi consid\'erer la fonction $ f(s) = h^2(s) - h^2(2a-s) = f^{-}(s) f^{+}(s) $.
Alors
\begin{equation} \label{intremmm}
f(a+i\tau) = 4 i \bigl|h^2(a+i\tau)\bigr| \sin \bigl(\arg h(a+i\tau)\bigr) \cos \bigl(\arg h(a+i\tau)\bigr),
\end{equation}
et les z\'eros de $f(s)$ sur la droite critique correspondent aux z\'eros de $f^+(s)$ et
$f^-(s)$ sur cette droite, de fa\c{c}on entrelac\'ee (sans tenir en compte leurs multiplicit\'es) si $\arg h(a+i\tau)$ est une fonction croissante.

Par contre, si $\arg h(a+i\tau)$ n'est pas une fonction croissante,
le lemme \ref{rectacritica} ne nous donne pas de bonne estimation pour $N_0'(T)$.
Il est donc n\'ecessaire de raffiner ce lemme.

Soient $0 \leq x_1 < x_2 < \dots < x_n < T$ dans la preuve du lemme \ref{compcreciente}.
On dit que $x_{i-1}$ est un {\it point \`a valeur enti\`ere d\'ecroissante} de $g(x)$, si $g(x_i) \leq g(x_{i-1})$.
Supposons maintenant que $g(x_i) \leq g(x_{i-1})$ pour $k$ valeurs diff\'erentes de $i$,  $1 < i \leq n$ ($g(x)$ a $k$ points \`a valeurs enti\`eres d\'ecroissantes). Soient
$$ c(T) = \# \{i \tq  1 < i \leq n, g(x_i)-g(x_{i-1}) = 0\}, $$
$$ d(T) = \# \{i \tq   1 < i \leq n, g(x_i)-g(x_{i-1}) = -1\}, $$
alors $k = c(T) + d(T)$ et
$$
\begin{array}{rcl}
g(x_n) - g(x_1) & = & \bigl(n-1-c(T)-d(T)\bigr) \cdot 1 + c(T) \cdot 0 + d(T) \cdot (-1) \\
    & = & \tilde{g}(T) - c(T) - 2 d(T)-1.
\end{array}
$$
Posons $K(T) = c(T) + 2d(T) - g(x_1)$. On a
$\tilde{g}(T) - g(x_n)  = K(T) + 1$.
En outre, on d\'efinit $\delta(T) = g(x_n) - g(T) + 1$, de fa\c{c}on que
$$ \tilde{g}(T) = g(T) + K(T) + \delta(T). $$
On a par d\'efinition $0 \leq \delta(T) \leq 2$. La fonction $\delta(T)$ est continue, sauf aux points o\`u elle (ou bien $g(T)$) prend
des valeurs enti\`eres, points o\`u elle est continue \`a gauche.

Notons certaines propri\'et\'es de $K(T)$:
\begin{enumerateu}
\item $K(T) \geq 0$, $K(T) \in \mathbb{Z}$;
\item $K(T)$ est une fonction croissante de $T$;
\item si $g(x)$ a $k$ points entiers \`a valeurs enti\`eres d\'ecroissantes dans $[0,T[$, alors $K(T) \geq k$.
En cons\'equence  $\tilde{g}(T) \geq g(T) + k$.
\end{enumerateu}

On applique ce dernier r\'esultat aux fonctions du lemme \ref{rectacritica},  $g_1(\tau) = \frac{1}{\pi} \arg h(a+i\tau)$ et $g_2(\tau) =\frac{1}{\pi} \arg h(a+i\tau)-\frac{1}{2}$,
pour obtenir l'estimation
\begin{equation} \label{sikgrande}
N_0'(T) \geq \mfrac{1}{\pi}\arg h(a+iT) -u_{\pm}+k
\end{equation}
au lieu de (\ref{dsuperior}).

\begin{lema} \label{sisolofinito} Soit $g:[0,\infty[ \to \re$ une fonction continue qui n'a qu'un nombre fini
de points \`a valeurs enti\`eres dans tout intervalle $[0,T[$, $T> 0$. Si
$g(x)$ n'a qu'un nombre fini de
points \`a valeurs enti\`eres d\'ecroissantes, alors il existe $k \in \mathbb{Z}$ tel que pour $T> 0$ assez grand,
$$ \tilde{g}(T) = \ment{g(T)} + k.$$
En particulier, $\tilde{g}(T) = g(T) + O(1)$. En particulier, si $g(T)$ est une fonction strictement croissante et $g(0) \in ]-1,0]$, alors
$k=0$ et
$$ \tilde{g}(T) = \ment{g(T)}. $$
\end{lema}

{\bf D\'emonstration.} S'il n'y a qu'un nombre fini de points \`a valeurs enti\`eres d\'ecroissantes
de $g(T)$, alors la fonction $K(T)$ est ultimement constante, $\delta(T) \neq 1$, $\delta(T) \neq 2$ pour $T$ assez grand.
S'il y a une infinit\'e de points \`a valeurs enti\`eres, on aura ultimement $\delta(T) < 1$, ce qui nous donnera $0 \leq \delta(T) < 1$.
En effet, si $\delta(T) \geq 1$ et s'il existe $T_1> T$ tel que $\delta(T_1) \in \mathbb{Z}$, alors pour la premi\`ere valeur de $T_1$ telle que cela arrive,
on a  $\delta(T_1) =1$ ou $2$ (par continuit\'e \`a gauche de $\delta(T)$), absurde.
S'il n'y a qu'une quantit\'e finie de points \`a valeurs enti\`eres, on aura soit $0 < \delta(T) < 1$ pour $T$ assez grand, soit $1 < \delta(T) < 2$ pour $T$ assez grand.
Tout cela nous permet d'\'ecrire
$$
\tilde{g}(T)+\ment{-\delta(T)} = \ment{g(T)} + K(T),
$$
o\`u $K(T)$ est constante, $\ment{-\delta(T)}$ est la constante $0$ ou $-1$, pour $T$ assez grand. \hfill{$\blacksquare$}\\

Si la fonction $\varphi(\tau)= \arg h(a+i\tau)$ satisfait \`a la condition du lemme \ref{sisolofinito},
on peut obtenir une meilleure estimation pour $N_0'(T)$ que celle donn\'ee par le lemme \ref{rectacritica} dans le cas o\`u cette
fonction n'est pas croissante. Une condition sufissante pour cela est
que $\varphi(\tau)$ soit ultimement croissante.

\section{R\'esultat principal} \label{princcc}

L'id\'ee de base de la m\'ethode de P. R. Taylor, ayant pour objectif de
d\'emontrer la proposition \ref{prtaylor}, est la suivante: d'abord, on calcule $N(T)$
\`a l'aide du principe de l'argument. Ensuite, on estime $N_0'(T)$ avec le lemme \ref{rectacritica}.
La diff\'erence $N(T) - N_0'(T)$
contient un terme d'erreur; on montre que la moyenne de la contribution de ce terme
est n\'egligeable gr\^ace \`a un th\'eor\`eme de Littlewood.
Il est \`a noter que la m\^eme id\'ee est \`a l'origine de la m\'ethode de v\'erification num\'erique
de l'hypoth\`ese de Riemann imagin\'ee par Turing \cite[\S 8.2]{edwards}.

\begin{teor} \label{teoalin} Soit $a \in \re$, $h(s)$ une fonction
m\'eromorphe sur $\ce$, r\'eelle sur la droite r\'eelle, n'ayant qu'un
nombre fini de p\^oles dans $\ce$, un nombre fini de z\'eros dans le demi-plan $\sigma > a$, holomorphe et
sans z\'ero sur la droite critique $\sigma = a$.
On d\'efinit la fonction
$$ f(s) = f^{\pm}(s) = h(s) \pm h(2a-s) $$
(en particulier $f(2a-s)=\pm f(s)$). On suppose que la fonction
$$ F(s) = \Frac{h(2a-s)}{h(s)} $$
 satisfait
\begin{enumerate}
\romano
\item \label{princcc1} pour chaque $\eta > 0$, il existe $\sigma_0 = \sigma_0(\eta) > a$ tel que $\bigl|F(s)\bigr| < \eta $
si $\sigma \geq \sigma_0$, $\tau \in \re$;
\item \label{princcc2} pour chaques $\eps > 0$ et $\sigma_0>a$, il existe une suite  $(T_n)_n$
telle que $ \ds\lim_{n \to \infty} T_n = +\infty$ et
$\bigl|F(s)\bigr| < e^{\eps|s|}$
pour $a \leq \sigma \leq \sigma_0$, $|\tau|=T_n$, $n \geq 1$.
\end{enumerate}
On note aussi, pour une fonction $m(s)$, $n_{m, \sigma > a}$ le
nombre de ses z\'eros r\'eels $\sigma$ tels que $\sigma
> a$, $N_{m, \sigma > a}$ le nombre des z\'eros (r\'eels ou complexes)
tels que $\sigma > a$, $n_{m,a}$ la multiplicit\'e de $s=a$ comme z\'ero de $m(s)$
(si c'est le cas, sinon $n_{m,a}=0$), et  $P_{m,\sigma>a}$ le nombre
des p\^oles de $m(s)$, tels que $\sigma > a$ (on consid\`ere les
multiplicit\'es dans tous les cas). Alors
\begin{equation} \label{naa}
N(T) - N_0(T) \leq N(T) - N_0'(T) \leq u_{\pm}-n_{f,\sigma > a} -
\mfrac{n_{f,a}}{2} + P_{f,\sigma > a} + N_{h,\sigma >a} - P_{h,
\sigma > a}.
\end{equation}
 pour $T > 0$, o\`u $u_+= \mfrac{1}{2}$, $u_-=1 $,
$n_{f,a}= 0$ (cas $+$) ou $n_{f,a} \geq 1$ est impaire (cas $-$). En
particulier, tous les z\'eros de $f(s)$, sauf un nombre fini, se trouvent sur la
droite $\sigma = a$ et sont simples. Le membre de gauche de
(\ref{naa}) est de plus un nombre positif pair.
\end{teor}

{\bf D\'emonstration.} Comme $h(a) \neq 0$, on peut supposer $h(a)> 0$,
puisque si $h(a)<0$, on change $h(s)$ en $-h(s)$ dans l'analyse.
On pose $g(s) = 1 \pm F(s)$ (selon le signe qui d\'efinit $f(s)$),
de telle fa\c{c}on que $f(s) = g(s) h(s)$.
Il faut remarquer que la condition {\it(\ref{princcc2})} reste valable, inchang\'ee, pour $g(s)$.

Soit $0 < \eta \leq \frac{1}{2}$. Par {\it(\ref{princcc1})},
il existe $\sigma_0 > \max\{a,0\}$ tel
que $\bigl|1-g(s)\bigr|\leq \eta$ pour $\sigma \geq \sigma_0$; puis
$\bigl|g(s)\bigr| \geq \frac{1}{2}$
et $f(s)$ n'a pas de p\^oles, $h(s)$ n'a ni de z\'eros ni de p\^oles,
pour $\sigma \geq \sigma_0$ (puisqu'ils sont en nombre fini).
Pour $T > 0$, on consid\`ere le rectangle $R$ d\'efini par
$$ -T \leq \tau \leq T, \quad 2a-\sigma_0 \leq \sigma \leq \sigma_0. $$

Suivant \cite[\S 9.9]{ztitchmarsh}, on d\'efinit le logarithme de $g(s)$ de la
fa\c{c}on suivante: si $\tau$
n'est pas l'ordonn\'ee d'un z\'ero, en partant d'une valeur pour $\log g(\sigma_0)$
(dans notre cas $\log g(\sigma_0)=0$, puisque $g(\sigma_0)>0$),
on d\'efinit $\log g(s)$ par variation continue le long des segments
qui lient $\sigma_0, \sigma_0+i\tau$ et $s$ dans l'ordre. Si $\tau$ est
l'ordonn\'ee d'un z\'ero, on d\'efinit
$$ \log g(s) = \lim_{\delta \to 0^+} \log g(s+i\delta). $$
On peut remarquer que sous les conditions donn\'ees, $\log g(s)$
ne d\'epend pas de $\sigma_0$, pour $\sigma_0$ assez grand.

L'application $\sigma \mapsto 2a-\sigma$ \'etablit une
correspondance entre les z\'eros et p\^oles $s=\sigma$ de
$f(\sigma)$ r\'eels, avec $\sigma > a$, et les z\'eros et p\^oles
avec $\sigma < a$. Donc l'exc\`es du nombre de z\'eros r\'eels
de $f(s)$ ($2n_{f,\sigma>a} + n_{f,a}$) sur le nombre de p\^oles (r\'eels ou complexes) de
$f(s)$ est
$$ 2n_{f,\sigma> a} + n_{f,a} - 2 P_{f,\sigma > a}. $$

Si $T$ n'est pas l'ordonn\'ee d'un z\'ero, par le principe de l'argument, on a
$$ 2N(T) + 2n_{f,\sigma> a} + n_{f,a} - 2 P_{f,\sigma > a} = \Frac{1}{2\pi} \Delta_R \arg f(s)+ C(T), $$
o\`u $\Delta_R$ note la variation dans le rectangle $R$, et $C(T) = 0$ pour $T$ assez grand
(lorsque $R$ contient tous les p\^oles de $f(s)$), disons $T > T_f > 0$).
Des sym\'etries $\ovl{f(s)} = f(\ovl{s})$ et $f(2a-s) = \pm f(s)$ de $f(s)$, on en d\'eduit
$$ \Delta_R \arg f(s) = 4 \Delta \arg f(s), $$
o\`u $\Delta$ d\'esigne la variation de $\sigma_0$ \`a $\sigma_0+iT$,
puis \`a $a+iT$. Puis
\begin{equation} \label{bsuperior}
\begin{array}{rcl}
N(T) &=& \mfrac{1}{\pi} \Delta \arg f(s)
-n_{f,\sigma> a} - \frac{n_{f,a}}{2} + P_{f,\sigma > a} +C(T) \\[0.2cm]
    &=& R(T)+S(T)
-n_{f,\sigma> a} - \frac{n_{f,a}}{2} + P_{f,\sigma > a}+C(T),
\end{array}
\end{equation}
o\`u
$$ R(T) = \Frac{1}{\pi} \Delta \arg h(s), \quad S(T) = \Frac{1}{\pi}
\Delta \arg g(s). $$
Si $T$ est l'ordonn\'ee d'un z\'ero, on d\'efinit $R(T) = R(T+0)$, $S(T) = S(T+0)$.
Par d\'efinition
$$S(T) = \frac{1}{\pi} \Im \log g(a+iT). $$
On prend maintenant le rectangle $R'$ d\'efini
par
$$ 0 \leq \tau \leq T, \quad a \leq \sigma \leq \sigma_0. $$
Par un th\'eor\`eme de Littlewood \cite[\S 9.9]{ztitchmarsh},
$$\int_{\partial R'} \log g(s)ds = -2 \pi i \int_{a}^{\sigma_0} \nu(\sigma)d\sigma,$$
o\`u $\nu(\sigma')$ d\'esigne l'exc\`es du nombre de z\'eros par rapport au nombre de p\^oles de $g(s)$ dans le rectangle
$ \sigma' < \sigma \leq \sigma_0$, $0 < \tau \leq T$. En particulier
$$ \Re \Bigl( \int_{\partial R'} \log g(s)ds \Bigr) = 0, $$
d'o\`u
\begin{equation} \label{intesedet}
\pi \int_0^T S(\tau) d\tau = \int_a^{\sigma_0} \log \bigl|g(\sigma+iT)\bigr|d\sigma+ \int_0^T \arg g(\sigma_0+i\tau)d\tau - I(\sigma_0),
\end{equation}
avec $I(\sigma_0)= \displaystyle \int_a^{\sigma_0} \log |g(\sigma)|d\sigma$.
Soit $\eps > 0$. La condition {\it(\ref{princcc2})} entra\^ine $\log |g(\sigma + i T_n)| < \eps (\sigma + T_n)$ pour $a \leq \sigma \leq \sigma_0$, $n \geq 1$. Puis
\begin{equation} \label{csuperior}
\pi \int_0^{T_n} S(\tau) d\tau \leq \sigma_0 \eps (\sigma_0 +  T_n) + \eta T_n - I(\sigma_0).
\end{equation}

Maintenant on estime la variation de l'argument de $h(s)$ dans le
rectangle $R''$ donn\'e par
$$ -T \leq \tau \leq T, \quad a \leq \sigma \leq \sigma_0, $$
o\`u $\sigma_0$ et $T$ sont assez grands de telle fa\c{c}on que tous les z\'eros et p\^oles de $h(s)$ restent \`a l'int\'erieur de $R''$.
Soit $\Delta_R'\arg h(s)$ la variation de l'argument de $h(s)$, obtenue
par variation continue le long des segments qui joignent $a-iT$, $\sigma_0-iT$, $\sigma_0+iT$, $a+iT$. Soit
aussi $\arg h(a+iT)$ la d\'etermination continue de l'argument de $h(s)$ sur la droite critique $\sigma = a$ telle que $\arg h(a) = 0$.
Alors
$$ \Delta_R' \arg h(s) + \arg h(a-iT) - \arg h(a+iT) = 2 \pi (N_{h,\sigma > a} - P_{h,\sigma > a}). $$
Las sym\'etrie $\ovl{h(s)} = h(\ovl{s})$ entra\^ine, comme dans le cas de $f(s)$, que
$$ \Delta_R' \arg h(s) = 2 \Delta \arg h(s)= 2 \pi R(T) \quad \textrm{ et } \quad \arg h(a-iT) = -\arg h(a+iT), $$
d'o\`u
$$ R(T) - \mfrac{1}{\pi}\arg h(a+iT) = N_{h,\sigma > a} - P_{h,\sigma > a}. $$
Donc
\begin{equation} \label{xsuperior}
R(T) - \mfrac{1}{\pi} \arg h(a+iT) =
N_{h,\sigma > a} - P_{h, \sigma > a} + D(T).
\end{equation}
pour $T > 0$, o\`u $D(T) = 0$ pour $T$ assez grand (lorsque $R$ contient tous les z\'eros et p\^oles de $h(s)$),
$T > T_h > 0$.
Maintenant, (\ref{bsuperior}) et (\ref{xsuperior}) entra\^inent
$$ N(T) = \frac{1}{\pi} \arg h(a+iT) -n_{f,\sigma> a} - \frac{n_{f,a}}{2} + P_{f,\sigma > a} + N_{h,\sigma > a} - P_{h,\sigma > a} + C(T) + D(T) + S(T), $$
ce qui avec le lemme \ref{rectacritica}, nous donne
$$ N(T) - N_0'(T) \leq B_a + C(T) +D(T) + S(T) $$
pour $T > 0$, o\`u
$$B_a = u_{\pm} -n_{f,\sigma> a} - \frac{n_{f,a}}{2} + P_{f,\sigma > a} + N_{h,\sigma > a} - P_{h,\sigma > a},$$
la borne attendue.

Il faut \'eliminer la contribution de $S(T)$ \`a la derni\`ere in\'egalit\'e.
Pour $T>0$ fix\'e, on prend $n$ assez grand pour que $T_n > \max\{T,T_f,T_h\}$. On int\`egre pour obtenir
$$
\begin{array}{rcl}
(T_n-T)\bigl( N(T)-N_0'(T) \bigr) &\leq&  \ds\int_0^{T_n} \bigl(
N(\tau)-N_0'(\tau)
\bigr) d\tau  \\
  &  \leq & \ds B_a T_n + \int_0^{T_f} C(\tau) d\tau +\int_0^{T_h} D(\tau) d\tau + \int_0^{T_n} S(\tau) d\tau\\[0.3cm]
  &  \leq & \ds B_a T_n +  \int_0^{T_f} C(\tau) d\tau +\int_0^{T_h} D(\tau) d\tau +  \eps \frac{\sigma_0^2}{\pi} +
  \eps \frac{\sigma_0 T_n}{\pi} + \frac{\eta}{\pi} T_n - \frac{I(\sigma_0)}{\pi}
\end{array}
$$
On divise par $T_n$ et on fait $n \to \infty$, d'o\`u
$$ N(T) - N_0'(T) \leq B_a + \eps \frac{\sigma_0}{\pi}  + \frac{\eta}{\pi}. $$
Puis on peut faire $\eps \to 0$ (n'ayant plus $T_n$) et, finalement, $\eta \to 0$
(n'ayant plus $\sigma_0$), pour obtenir (\ref{naa}).

En particulier, si $B_a < 2$, on a $N(T) = N_0(T)$ puisque $N(T)-N_0(T)$ est un nombre
entier pair. \hfill{$\blacksquare$}\\

{\bf Remarques.}
\begin{enumerate}
\alfabetico
\item \label{estata} La condition d'avoir $h(a) \neq 0$ est artificielle, et elle peut \^etre relax\'ee.
Soit $s=a$ un z\'ero ou un p\^ole de $h(s)$, on \'ecrit $h(s) = (s-a)^m h_1(s)$,
o\`u $m \in \mathbb{Z}$, $h_1(a) \neq 0$. Alors on peut factoriser $(s-a)^m$ et obtenir une nouvelle fonction $f_1(s)$ telle que
$f(s) = (s-a)^m f_1(s)$, en gardant le signe si $m$ est pair, et en le changeant si $m$ est impair.
Les conditions sur la nouvelle fonction $F_1(s) = h_1(2a-s)/h_1(s)$ restent les m\^emes.
On peut aussi s'occuper facilement des z\'eros et des p\^oles de $h(s)$ sur la droite critique, le cas o\`u il y en a
et ils sont en nombre fini.
Si $s_0=a+i \tau_0$ (avec $\tau_0 \neq 0$) est un z\'ero de $h(s)$, $h(a+i\tau_0) = 0$, alors $\ovl{h(a+i\tau)} = h(a-i\tau_0) = 0$.
On peut factoriser $h(s) = (s-s_0)(s-\ovl{s_0}) h_1(s) = \bigl(
(s-a)^2+\tau_0^2 \bigr) h_1(s)$, puis
$$
\begin{array}{rcl}
f(s) &=& h(s) \pm h(2a-s) = \bigl( (s-a)^2+\tau_0^2 \bigr) h_1(s) \pm \bigl( (a-s)^2+\tau_0^2 \bigr) h_1(2a-s)\\
     &=& (s-s_0)(s-\ovl{s_0}) f_1(s).
\end{array}
$$
La fonction $F(s)$ est la m\^eme, donc on revient au cas pr\'ec\'edent.
On peut \'eliminer de cette fa\c{c}on un nombre fini de z\'eros de la fonction $h(s)$. De m\^eme pour les p\^oles.

\item La fonction $F(s)$ dans le th\'eor\`eme satisfait $\bigl|F(a+i\tau)\bigr| = 1$, pour tout $\tau \in \re$.
La condition {(\it\ref{princcc1})}, est un peu forte, elle n'est pas satisfaite par les polyn\^omes, qui ont
besoin d'une d\'emarche additionnelle pour leur \'etude, que l'on verra au
\S \ref{policaso}. Mais cela \'evite aussi les cas triviaux (eg. si $p(s)$ est un polyn\^ome pair,
alors $f^-(s) = 0$). Dans le cas des polyn\^omes, on peut encore obtenir {\it(\ref{princcc1})} par perturbation du probl\`eme initial.
\item \label{phragmen}
Soient $a_1, \dots, a_m$ les z\'eros de $h(s)$ dans le demi-plan $\sigma > a$, et $b_1,\dots,b_n$ les p\^oles
de $h(s)$ dans le demi-plan $\sigma < a$, en tenant en compte les multiplicit\'es. On d\'efinit
$$ R(s) = \frac{(s-a_1)\cdots (s-a_m)(2a-s-{b_1}) \cdots (2a-s-{b_n})}{(2a-s-\ovl{a_1}) \cdots (2a-s-\ovl{a_m})(s-\ovl{b_1})\cdots (s-\ovl{b_n})}
 $$
Alors, $\bigl| R(a+i\tau) \bigr|=1$ pour $\tau \in \re$, $\ds\lim_{|s|\to \infty} \bigl| R(s) \bigr|=1$, et
la fonction
$$ F_1(s) = R(s) F(s) $$
est reguli\`ere dans le demi-plan $\sigma \geq a$. Cette fonction satisfait aussi les conditions {\it(\ref{princcc1})}, {\it(\ref{princcc2})}
du th\'eor\`eme \ref{teoalin}. En particulier, $\bigl| F_1(a+i\tau) \bigr|=1$ pour $\tau \in \re$.
Soit $\sigma_0 > a$ tel que $\bigl| F_1(\sigma+i\tau) \bigr|< 1/2$ pour $\sigma \geq \sigma_0$, $\tau \in \re$.
Dans l'intervalle $a \leq \sigma \leq \sigma_0$, la condition {\it(\ref{princcc2})} est \`equivalente \`a la condition que, pour tout
$\eps > 0$, $ \bigl| F_1(s) \bigr| = O( e^{\eps |\tau|} ) $ pour $a \leq \sigma \leq \sigma_0$, $|\tau|=T_n$, $n$
assez grand. Le principe de Phragm\'en-Lindel\"of (\cite[\S 5.65]{ftitchmarsh}, \cite[Theorem 12.8]{rudin}) s'applique
(en remarquant que la condition sur la suite de valeurs $|\tau|=T_n$ est suffisante) pour nous donner
$ \bigl| F_1(s) \bigr| < (1/2)^{\sigma_0-\sigma}<1$ pour $a < \sigma \leq \sigma_0$. Donc
\begin{enumerate}
\item[\hypertarget{phragmen1}(1)] si $h(s)$ est une fonction enti\`ere sans z\'eros dans le demi-plan $\sigma \geq a$, alors
$$ \bigl| F(s) \bigr| < 1 \mbox{ pour } \sigma > a $$
(puisque $F(s)=F_1(s)$ dans ce cas);
\item[\hypertarget{phragmen2}(2)] en g\'en\'eral, pour chaque $\eta > 0$, il existe $R> 0$ tel que
$$  \bigl| F(s)\bigr| < 1+\eta \mbox{ pour } \sigma > a, |s| > R. $$
\end{enumerate}
Une condition similaire appara\^it dans la g\'en\'eralisation du th\'eor\`eme d'Hermite-Biehler aux
fonctions enti\`eres de type exponentiel \cite[chapitre IV, \S 5]{cebotarev}, au lieu des conditions asymptotiques
du th\'eor\`eme \ref{teoalin}:
pour chaque $\eps >0$, il existe une suite de demi-cercles dans le demi-plan $\sigma > a$,
centr\'es en $s=a$ et de rayon convergeant vers l'infini, le long desquels $|F(s)| < e^{\eps |s|}$.
Sous cette condition, on peut retrouver (\hyperlink{phragmen1}{1}) et (\hyperlink{phragmen2}{2}), et (\hyperlink{phragmen1}{1}) suffit pour montrer que tous les z\'eros de la fonction $f(s)$ sont
simples et align\'es (proposition \ref{copiado}).
Certes, c'est une condition plus faible que {\it(\ref{princcc1})}, mais on pourra alors proc\'eder par perturbation
du probl\`eme initial, comme dans le cas des polyn\^omes au \S \ref{policaso}.
\item On verifiera, pour les fonctions avec lesquelles on va travailler, la condition $ F(s) = O(|\tau|^A)$
pour une constante $A > 0$, lorsque $|\tau|$ est assez grand, uniform\'ement pour $\sigma \geq a$, au lieu de {\it(\ref{princcc2})}. Cela est naturel
vu qu'on applique ce r\'esultat \`a l'\'etude de fonctions provenantes de s\'eries de Dirichlet.
Une fonction $F(s)$ avec une telle condition de croissante est appel\'ee
d'ordre fini (cf. \cite[\S 9.4]{ftitchmarsh}), ou de croissance polynomiale (cf. \cite[\S 5.2]{iwaniec})
dans le demi-plan $\sigma \geq a$.
\end{enumerate}

\begin{coro} \label{coralin} Sous les conditions du th\'eor\`eme \ref{teoalin}, si $h(s)$ est une fonction enti\`ere,
alors
\begin{equation} \label{nee}
N(T) - N_0(T) \leq N(T) - N_0'(T) \leq  u_{\pm} - \frac{n_{f,a}}{2} + N_{h,\sigma >a}-n_{f,\sigma
> a}.
\end{equation}
En particulier, si $h(s)$ n'a pas de z\'ero avec $\sigma \geq a$, alors tous les z\'eros de $f(s)$
sont sur la droite $\sigma = a$ et ils sont simples.
\end{coro}

{\bf D\'emonstration.} La formule (\ref{nee}) est evidente. Si $n_{h,\sigma > a} = 0$, comme le cot\'e droit de cette formule
est toujours positif,
$$ 0 \leq u_{\pm}- \Frac{n_{f,a}}{2} -n_{f,\sigma > a} $$
d'o\`u $n_{f,\sigma > a} = 0$, et $n_{f,a} = 0$ dans le cas $+$, $n_{f,a} = 1$ dans le cas $-$. On r\'eintroduit
ces valeurs dans la formule (\ref{nee}) pour obtenir $N(T) = N_0(T) = N_0'(T)$.
\hfill{$\blacksquare$}\\

En fait, le corollaire \ref{coralin} est \'equivalent au th\'eor\`eme \ref{teoalin}.
Un p\^ole de $h(s)$ est transform\'e en un z\'ero en multipliant par un facteur appropri\'e. La raison de garder l'enonc\'e du th\'eor\`eme pour des fonctions
m\'eromorphes est de simplifier l'application du r\'esultat.
De m\^eme on aurait pu poser $a=0$, une simple translation permettant d'en d\'eduire le cas g\'en\'eral.

\section{R\'epartition globale des z\'eros} \label{repglobale}

Sous les conditions du th\'eor\`eme \ref{teoalin}, la r\'epartition des z\'eros de $f(s)$ est essentiellement
la r\'epartition des z\'eros sur la droite critique.

\begin{teor} \label{enezero} Sous les conditions du th\'eor\`eme \ref{teoalin}, si $\arg h(a+it)$ est
une d\'etermination continue de l'argument de $h(s)$ sur la droite $\sigma=a$, avec $\arg h(a)=0$, alors
$$ N_0'(T) = \frac{1}{\pi} \arg h(a+iT) + O(1). $$
Les z\'eros des fonctions $f(s) = h(s) \pm h(2a-s)$ sont
ultimement simples et entrelac\'es. Plus pr\'ecisement,
il existe un entier $d_{f,a} \geq 0$ tel que
\begin{equation} \label{sobrelarecta}
N_0'(T) = \gment{\frac{1}{\pi} \arg h(a+iT) - u_{\pm} } + d_{f,a}
\end{equation}
pour $T > 0$ assez grand. Dans le cas particulier o\`u $\arg h(a+iT)$ est une fonction croissante, on a $d_{f,a} = 0$.
\end{teor}

{\bf D\'emonstration.} Soit $f(s) = h(s) - h(2a-s)$.
La fonction $\varphi(\tau) = \frac{1}{\pi} \arg h(a+i\tau)$ n'a qu'un nombre fini de points \`a valeurs enti\`eres d\'ecroissantes.
En effet, si ce n'est pas le cas, on aura (\ref{sikgrande}) avec $k$ grand, ce qui dans (\ref{naa}) donne
une diminution de $k$ dans la borne \`a droite. Pour $k$ assez grand, le nombre de droite de (\ref{naa}) serait
negatif, ce qui est absurde. On applique le lemme \ref{sisolofinito} pour obtenir la premi\`ere partie du r\'esultat.
La positivit\'e de $d_{f,a}$ d\'ecoule de la formule explicite obtenue.

Pour la fonction $f(s) = h(s) + h(2a-s)$, on choisit $\varphi_1(\tau) = \frac{1}{\pi} \arg h(a+i\tau) - \frac{1}{2}$,
pour obtenir le m\^eme r\'esultat qu'avec $\varphi(\tau)$.

Finalement, on consid\`ere $f(s) = h^2(s)-h^2(2a-s)= f^+(s) f^-(s)$. La croissance ultime des points \`a valeurs enti\`eres
de $\frac{1}{\pi} \arg h^2(a+i\tau) = 2 \varphi(\tau) = 2\varphi_1(\tau)+1$
correspond \`a l'entrelacement des z\'eros des fonctions $f(s) = h(s) \pm h(2a-s)$
sur la droite $\sigma=a$ (voir (\ref{intremmm})). \hfill{$\blacksquare$}\\

\noindent {\bf Remarque.} Si $h(s)$ est une fonction enti\`ere sans z\'eros avec $\sigma \geq a$, on a $d_{f,a} = 0$,
par le corollaire \ref{coralin}. En effet, il suffit d'introduire l'estimation (\ref{sobrelarecta}) dans (\ref{nee}), et voir que cela
diminue de $d_{f,a}$ la borne obtenue.

D'un autre c\^ot\'e, la remarque (\ref{phragmen}.\hyperlink{phragmen1}{1}) au th\'eor\`eme \ref{teoalin} et
la proposition \ref{copiado} entra\^inent que la fonction $\varphi(\tau) = \arg h(a+i\tau)$
est croissante, d'o\`u on aura \'egalement $d_{f,a} = 0$.
Donc il est naturel de se poser la question suivante, motiv\'ee aussi dans la discussion faite
\`a la fin du \S \ref{seccerosrecta}.\\

\noindent {\bf Question.} Sous les conditions du th\'eor\`eme \ref{teoalin}, la fonction $\varphi(\tau)=\arg h(a+i\tau)$ est-elle ultimement croissante? Cela entra\^inerait le th\'eor\`eme \ref{enezero} directement.\\

On \'etablit un lemme qui nous aidera dans la suite.

\begin{lema} \label{sobrelarectauno} Soit $h(s)$ une fonction m\'eromorphe sans z\'eros sur la droite critique $\sigma = a $.
Si $\frac{h'(a)}{h(a)} < 0$, alors $d_{f,a} \geq 1$ dans la formule (\ref{sobrelarecta}), pour la fonction $f^-(s)=h(s)-h(2a-s)$.
\end{lema}

{\bf D\'emonstration.}  Soit $\theta(\tau) = \frac{1}{\pi} \arg h(a+i\tau)$, $f(s)= f^{-}(s)$. Un calcul direct montre que
$$ \theta'(\tau) = \Re \mfrac{h'(a+i\tau)}{h(a+i\tau)} $$
et en particulier $\theta'(0) = \mfrac{h'(a)}{h(a)}$. Le premier point \`a valeur enti\`ere de $\theta(\tau)$ est $x_1=0$ avec $\theta(x_1)=0$.
Si $\theta'(0) < 0$, la fonction $\theta(\tau)$ est d\'ecroissante dans un petit intervalle
contenant $0$; par continuit\'e, le deuxi\`eme point \`a valeur enti\`ere $x_2>0$
satisfait $\theta(x_2) \leq 0$. Donc $x_1=0$ est un point \`a valeur enti\`ere d\'ecroissante
de $\theta(\tau)$. Par l'estimation (\ref{sikgrande}), $d_{f,a} \geq 1$ et
\begin{equation} \label{eeesta}
N_0'(T) \geq \mfrac{1}{\pi} \arg h(a+iT),
\end{equation}
le r\'esultat attendu. \hfill{$\blacksquare$}\\

\section{Polyn\^omes et Stabilit\'e. Le th\'eor\`eme d'Hermite-Biehler} \label{policaso}

Un polyn\^ome $p(z) = a_0 z^n + a_{1} z^{n-1} + \cdots + a_{n-1} z + a_n$ avec $a_n \neq 0$
est dit {\it stable} si tous ses z\'eros sont dans le demi-plan $\sigma < 0$. On rappelle quelques propri\'et\'es des polyn\^omes stables
\`a coefficients r\'eels (cf. \cite{postnikov}):
\begin{enumerate}
\romano
\item \label{policoeff} Tous les coefficients de $p(z)$ ont le m\^eme signe. Si $a_0 > 0$, alors $a_i > 0$ pour tout $1 \leq i \leq n$ (th\'eor\`eme de Stodola). En particulier,
$p'(0) = a_1 \neq 0$.
\item \label{poliarg} La fonction argument phase $\varphi(\tau)=\arg p(i\tau)$ est une fonction croissante de $\tau$. Si on choisit $\arg p(0) = 0$,
la contribution \`a l'argument de chaque z\'ero r\'eel (donc n\'egatif) \`a l'infini (positif) est $\pi/2$,
et la contribution d'un z\'ero non r\'eel et son conjugu\'e est $\pi$. Donc
\begin{equation} \label{contribucion}
\lim_{\tau \to +\infty} \arg p(i\tau) = n \Frac{\pi}{2}.
\end{equation}
\end{enumerate}

\begin{prop} \label{perturbacionestable} Soit $p(s)$ un polyn\^ome r\'eel stable non constant. Pour chaque $y>0$,
on d\'efinit la fonction
$$ P(s;y) = y^s p(s) \pm y^{-s} p(-s). $$
On note par $N(y;T)$ le nombre de z\'eros de $P(y;s)$ avec $0< \tau < T$, $N(y_0;T)$
le nombre de z\'eros de $P(y;s)$ avec $\sigma=a$ et $0<\tau < T$, et $N'(y_0,T)$
le nombre de ces derniers z\'eros, sans compter les multiplicit\'es. Alors
\begin{enumerate}
\item \label{nuno} Si $y \geq 1$, tous les z\'eros de la fonction $P(y;s)$
sont simples et sur la droite $\sigma = 0$. De plus
$$ 0 \leq N(y;T) - \frac{T}{\pi} \log y + u_{\pm} \leq  \frac{n}{2}. $$
\item \label{esdos} Si $0 < y < 1$, alors
$$ N(y;T)-N_0(y;T) \leq N(y;T)-N_0'(y;T) \leq n, $$
le degr\'e de $p(s)$.
\end{enumerate}
\end{prop}

{\bf D\'emonstration.} Il faut remarquer que le th\'eor\`eme \ref{teoalin} ne peut pas \^etre appliqu\'e aux fonctions $f(s) = p(s) \pm p(-s)$,
puisque
$$\ds\lim_{|s| \to \infty} \frac{p(-s)}{p(s)} = \pm1$$
et la condition {\it(\ref{princcc1})} du th\'eor\`eme n'est pas satisfaite.
Pour atteindre ce cas, on introduit la famille de fonctions
$$P(y;s) = y^s p(s) \pm y^{-s} p(-s) $$
avec le param\^etre $y > 0$. La fonction $F(y;s) = y^{-2s} {p(-s)}/{p(s)}$
satisfait $F(y;s) = O(y^{-\sigma})$ pour $|s|$ suffisament grand,
donc les conditions du th\'eor\`eme \ref{teoalin} pour $y >1 $. D'apr\`es le corollaire \ref{coralin},
\'etant $H(y;s)=y^s p(s)$ sans z\'eros dans le demi-plan $\sigma \geq 0$,
tous les z\'eros de $P(y;s)$ sont simples et align\'es sur la droite $\sigma = 0$, pour $y > 1$.

Soit $T> 0$ qui n'est pas l'ordonn\'ee d'un z\'ero de $f(s) = P(1;s)$. Comme $N(1;T)=N(T)$ est fini,
par le th\'eor\`eme de Hurwitz \cite[\S 3.45]{ftitchmarsh}, pour $y> 1$ assez proche de $1$ (d\'ependant de $T$),
$$ N(T) \leq N(y;T), \quad N(T)-N_0(T) \leq N(y;T)-N_0(y;T). $$
De la deuxi\`eme estimation, il est clair que $f(s) = p(s) \pm p(-s)$ aura tous ses z\'eros sur la droite $\sigma = 0$.
Par contre, on ne peut rien dire en g\'en\'eral \`a propos de $N_0'(T)$, donc sur la simplicit\'e des z\'eros de $f(s)$.

La formule (\ref{sobrelarecta}) nous dit que pour montrer que $N_0'(y;T) \to N_0'(T)$ il nous
faut voir la continuit\'e de l'argument de $H(y;i\tau)$ et de la constante $d_{H(y,\cdot),0}$,
en fonction du param\`etre $y$. Heureusement,
$d_{H(y,\cdot),0} = 0$ puisque l'argument de $H(y;i\tau)$ est une fonction croissante. En fait
$$ \arg H(y;i\tau)  = \tau \log y + \arg p(i\tau), \quad y \geq 1. $$
Donc, si $\frac{1}{\pi}\arg p(i T) \notin \mathbb{Z}$
$$N(T) \leq \lim_{y \to 1^+} N(y;T) = \lim_{y \to 1^+} N_0'(y;T) = N_0'(T),$$
ce qui entra\^ine la simplicit\'e des z\'eros de $f(s)$. Cela montre (\ref{nuno}).

Si $0 < y < 1$, on prend $y^{-1} > 1$ au lieu de $y$ et $p(-s)$ au lieu de $p(s)$. Par application du corollaire
\ref{coralin}, on obtient $ N(T) - N_0(T) \leq N(T) - N_0'(T) \leq n $.  \hfill{$\blacksquare$}\\

Il est \'etonnant comme on obtient un r\'esultat pour une fonction avec un nombre fini de z\'eros
\`a partir d'autres fonctions qui en ont une infinit\'e. Mais il est toujours possible d'analyser
directement le cas $y=1$, en utilisant un comptage direct des z\'eros. Le lemme \ref{rectacritica} est toujours
applicable, ce qui nous donne
$$ N_0'(T) \geq \frac{1}{\pi} \arg p(iT)  - u_{\pm}. $$
On a quatre cas \`a analyser, selon le signe $\pm$ et la parit\'e de $n$; traitons un seul cas,
les autres \'etant analogues. Soit alors $n$ pair et $f(s) = p(s) - p(-s)$. Le z\'ero $s=0$ est simple, puisque $f'(0) = 2p'(0) \neq 0$ par (\ref{policoeff}).
Donc $n_{f,a} = 1$. En fait, $f'(s)$ ne peut changer de signe que pour $s=0$, \`a cause du fait
que ses coefficients son positifs, et il est impossible d'avoir $f'(s)$ localmente constante, puisque c'est un polyn\^ome. Donc $f(s)$
est strictement croissante ou d\'ecroissante sur les intervalles $]-\infty,0[$ et $]0,+\infty[$, ce qui montre que $s=0$ est le seul z\'ero r\'eel de $f(s)$.

Pour $T > 0$ assez grand, on a $N(T) = \frac{n-1}{2}$. En utilisant (\ref{contribucion}), on obtient
$$ 0 \leq N(T) - N_0(T) \leq N(T) - N_0'(T) \leq \frac{n-1}{2} - \arg p(iT) + 1 < 1  $$
pour $T > 0$ assez grand. Donc $N(T) = N_0(T) = N_0'(T)$ et tous les z\'eros de $f(s)$ sont sur la droite $\sigma = 0$
et simples.

On reformule le r\'esultat obtenu en termes d'autres fonctions. On \'ecrit
$2q(z^2) = p(z)+p(-z)$, $2z \cdot r(z^2) = p(z)-p(-z)$. Alors
les polyn\^omes $q(z)$ et $r(z)$ ont des racines r\'eelles negatives, simples et entrelac\'ees,
la condition n\'ecessaire de stabilit\'e du polyn\^ome $p(z)$ (cf. \cite{holtz});
une condition additionnelle assure l'\'equivalence entre ces conditions dans le cas des polyn\^omes.
Ce r\'esultat d'\'equivalence est connu comme le th\'eor\`eme d'Hermite-Biehler
un th\'eor\`eme qui a \'et\'e \'etendu \`a une classe de fonctions de type exponentiel
(cf. \cite[Part III, Lecture 27]{levin}, \cite{cebotarev}), incluant les polyn\^omes.
Le corollaire \ref{coralin} (ou le th\'eor\`eme \ref{teoalin})
g\'en\'eralise la condition n\'ecessaire de stabilit\'e du th\'eor\`eme d'Hermite-Biehler pour les fonctions enti\`eres
v\'erifiant les conditions {\it(\ref{princcc1})}, {\it(\ref{princcc2})} du th\'eor\`eme \ref{teoalin}. La condition
{\it(\ref{princcc1})} laisse dehors les polyn\^omes et complique leur \'etude (les polyn\^omes sont facilement analys\'es avec
des m\'ethodes comme celles du annexe \ref{seccionestabilidad}),
mais il montre son utilit\'e dans l'\'etude d'autres fonctions,
lorsqu'on admet des z\'eros \`a droite de la droite critique,
m\^eme en quantit\'e infinie (par extension de l'argument).

\section{La fonction $\zeta_2(s)$} \label{zeta2}

Consid\'erons la suite de fonctions m\'eromorphes
$$ \zeta_n(s) = \sum_{k=1}^n \frac{1}{k^s}  - \frac{n^{1-s}}{1-s}, \quad n \geq 1.$$
Cette suite converge uniform\'ement vers la fonction z\^eta $\zeta(s)$ sur les compacts du demi-plan $\sigma > 0$ priv\'e du point $s=1$,
au vue de la formule \cite[\S 3.5]{ztitchmarsh}
$$ \zeta(s) = \zeta_n(s) + s \int_n^{+\infty} \frac{\{x\}-\frac{1}{2}}{x^{s+1}}dx - \frac{1}{2} n^{-s} $$
valable pour $\sigma > 0$, o\`u $\{ x \} = x - \lfloor x \rfloor$ est la fonction partie fractionnaire. Le
premier terme non trivial
de cette suite est celui avec $n=2$;
on veut \'etudier la r\'epartition des z\'eros de cette fonction.

On pose $h(s) = e^{\alpha s}(s-\beta)$, avec $\alpha,\beta \in \re$, et $a = 0$.
 On \'etudiera la
fonction
$$ f(s) = h(s)-h(-s) = e^{\alpha s} (s - \beta) + e^{-\alpha s} (s + \beta) = 2 s \cosh(\alpha s) - 2 \beta \sinh(\alpha s), $$
qui satisfait $f(-s) = -f(s)$. Lorsque $\alpha = \log(2)/2$, $\beta = 1$, on obtient $ f(s) = -2^{s/2}(1-s) \zeta_2(s) $.

L'annulation d'un quelconque des param\`etres $\alpha, \beta$ nous
am\`ene dans des cas triviaux. Si $\alpha = 0$, on a $f(s) = 2s$, et
si $\beta = 0$, $f(s) = 2 s \cosh(s)$, des fonctions ayant que des
z\'eros r\'eels. Le cas $\alpha < 0$ est facilement ramen\'e au cas
$\alpha > 0$ en prenant $-h(-s)$ au lieu de $h(s)$.

\begin{prop} \label{alphabeta} Soient $\alpha> 0$ , $\beta \neq 0$. Tous les z\'eros non r\'eels de la fonction
$$ f(s) = e^{\alpha s} (s - \beta) + e^{-\alpha s} (s + \beta) $$
sont simples, et sur la droite $\sigma = 0$. Quant aux z\'eros
r\'eels de $f(s)$:
\begin{enumerate}
\romano
\item Si $\beta < 0$, $s=0$ est le seul z\'ero r\'eel de $f(s)$, et il est simple;
\item Si $\beta > 0$,
\begin{enumerate}
\item pour $\alpha < \frac{1}{\beta}$, $s=0$ est le seul z\'ero r\'eel de $f(s)=0$, et il est simple;
\item pour $\alpha = \frac{1}{\beta}$, $s=0$ est le seul z\'ero r\'eel de $f(s)=0$, et il est triple;
\item pour $\alpha > \frac{1}{\beta}$, $f(s)$ a trois z\'eros r\'eels simples $s=0$, $s=\gamma>0$ et $s=-\gamma$.
\end{enumerate}
\end{enumerate}
\end{prop}

{\bf D\'emonstration.} La fonction
$$ F(s) = \frac{h(-s)}{h(s)} = e^{-2\alpha s} \Bigl( 1 + \Frac{2\beta}{s-\beta} \Bigr)$$
satisfait les conditions du th\'eor\`eme \ref{teoalin}. Donc
$$
0 \leq N(T)-N_0(T) \leq N(T) - N_0'(T) \leq 1 - n_{f,\sigma > 0} -
\Frac{n_{f,0}}{2} + N_{h,\sigma > 0} \leq \frac{3}{2}
$$
puisque $n_{f,0} \geq 1$ est impaire, et $N_{h,\sigma>0} \leq 1$.
Maintenant $N(T)-N_0(T)$ est pair, donc $N(T)=N_0(T)$ et $N_0(T) \leq N_0'(T) +1$.

\begin{enumerate}
\romano
\item Supposons $\beta < 0$. Alors $N_{h,\sigma > 0} = 0$, $n_{f,0} \geq
1$ et
$$ 0 \leq N(T)-N_0(T) \leq N(T) - N_0'(T) \leq 1 - \Frac{n_{f,0}}{2} \leq \frac{1}{2}, $$
donc $N_0'(T) = N_0(T) = N(T)$, et tous les z\'eros de $f(s)$ sont
simples et sur la droite $\sigma  = 0$ par le corollaire \ref{coralin}.
\item
On suppose maintenant $\beta > 0$. Alors $N_{h,\sigma > 0} = 1$ et
\begin{equation} \label{estimxy}
0 \leq N(T)-N_0(T) \leq N(T) - N_0'(T) \leq 2 - n_{f,\sigma > 0}- \frac{n_{f,0}}{2},
\end{equation}
On calcule $f'(0) = 2(1 - \alpha \beta)$ et $\frac{h'(0)}{h(0)} = \alpha- \frac{1}{\beta} < 0$.
\begin{enumerate}
\item Si $\alpha < \frac{1}{\beta}$, le lemme \ref{sobrelarectauno} nous permet de r\'eduire
d'une unit\'e la borne dans (\ref{estimxy}), puisque $\frac{h'(0)}{h(0)} < 0$. Donc
$$N(T)-N_0(T) \leq N(T) - N_0'(T) \leq 1 - n_{f,\sigma > 0}- \frac{n_{f,0}}{2} \leq \frac{1}{2},$$
d'o\`u $N(T) = N_0(T) = N_0'(T)$, puis $n_{f,\sigma_0}= 0$ et $n_{f,0}=1$.
\item Si $\alpha = \frac{1}{\beta}$, alors $f'(0) = 0$, d'o\`u $n_{f,0}\geq 3$. Cela dans
(\ref{estimxy}) entra\^ine $N(T) = N_0(T) = N_0'(T)$ et $n_{f,\sigma_0}= 0$.
\item Si $\alpha > \frac{1}{\beta}$, alors $f'(0) < 0$, donc il existe $x'>0$ (et proche de $0$) tel
que $f(x') < 0$. Mais
$$ \lim_{\sigma \to +\infty} f(\sigma) = +\infty, $$
d'o\`u il existe $\gamma > x'>0$ tel que $f(\gamma) = 0$, et $n_{f,\sigma > 0} \geq 1$.
Cela dans (\ref{estimxy}) nous donne encore
$N(T) = N_0(T) = N_0'(T)$ et $n_{f,0}=1$. \hfill{$\blacksquare$}
\end{enumerate}
\end{enumerate}

\begin{coro} \label{inalcanzable} Les z\'eros de la fonction
$$ \zeta_2(s) = 1 + 2^{-s} - \Frac{2^{1-s}}{1-s} $$
sont tous simples et sur la droite $\sigma = 0$.
\end{coro}

{\bf D\'emonstration.} Lorsque $\alpha = \log(2)/2$, $\beta = 1$, on obtient dans la proposition \ref{alphabeta}
$ f(s) = -2^{s/2}(1-s) \zeta_2(s) $, avec $ 0 < \alpha < \frac{1}{\beta}$. \hfill{$\blacksquare$}\\

\section{Translat\'ees de la fonction z\^eta de Riemann} \label{aquiestimation}

On consid\`ere les diff\'erences des translat\'ees de la fonction z\^eta complet\'ee,
donn\'es pour $\alpha > 0$ par
$$ f_\alpha(s)=\zeta^*(s+\alpha) - \zeta^*(s-\alpha).$$

Ces fonctions s'\'ecrivent sous la forme
$f_\alpha(s) = h_\alpha(s) - h_\alpha(1-s)$, o\`u $h_\alpha(s)= \zeta^*(s+\alpha)$. Dans
ce cas $a = \frac{1}{2}$.

Dans un premier temps, on consid\`ere $\alpha > \frac{1}{2}$.
Sous cette condition, tous les z\'eros de $h_\alpha(s)$ sont \`a gauche de la droite critique,
mais ses p\^oles vont pousser quelques z\'eros en dehors de la droite.
Ce cas a \'et\'e \'etudi\'e par P. R. Taylor en \cite[\S 1]{taylor},
dans le but d'\'etudier la r\'epartition des z\'eros de la fonction
$f_{1/2}(s)=\zeta^*\bigl(s+\frac{1}{2}\bigr)-\zeta^*\bigl(s-\frac{1}{2}\bigr)$.

\begin{prop} Soit $\alpha > \frac{1}{2}$, $\alpha^*\approx 6.81707$ tel que $(\zeta^*)'(\frac{1}{2}+\alpha^*)=0$. Tous les z\'eros non r\'eels de la fonction
$$ f_\alpha(s)= \zeta^*(s+\alpha)-\zeta^*(s-\alpha)$$
sont sur la droite $\sigma= \frac{1}{2}$. Quant aux z\'eros r\'eels:
\begin{enumerate}
 \item $s=\frac{1}{2}$ est un z\'ero simple de $f_\alpha(s)$ si $\alpha \neq \alpha^*$,
triple si $\alpha = \alpha^*$;
 \item $f_\alpha(s)$ a deux z\'eros r\'eels simples $s=\rho_\alpha $ et $s=1-\rho_\alpha$,
o\`u $\rho_\alpha \in ]1+\alpha,+\infty[$.
Ces sont les seuls z\'eros r\'eels de $f_\alpha(s)$ diff\'erents de $s=\frac{1}{2}$.
\end{enumerate}
\end{prop}

{\bf D\'emonstration.}
Suivant le th\'eor\`eme \ref{teoalin}, on pose
$$F_\alpha(s) = \Frac{\zeta^*(s-\alpha)}{\zeta^*(s+\alpha)}. $$
Par la formule complexe de Stirling \cite[\S 4.42]{ftitchmarsh}, pour $\sigma \geq \frac{1}{2}$, $|s|$
suffisament grand
$$ \log \Gamma \big( \mfrac{1}{2} s + \mfrac{1}{2} \alpha \bigr) = \bigl( \mfrac{s}{2} + \mfrac{\alpha}{2} - \mfrac{1}{2} \bigr) \log \bigl(\mfrac{s}{2}\bigr) - \mfrac{s}{2} + \mfrac{1}{2} \log 2\pi + O \bigl( \mfrac{1}{|s|} \bigr), $$
et $\zeta(s)=1+O(2^{-\sigma})$ lorsque $\sigma \geq \sigma_0$ puis
$$ \begin{array}{rcl}
\log h_\alpha(s) &=& (-\mfrac{1}{2}s-\frac{1}{2}\alpha) \log \pi + \log
\Gamma(\mfrac{1}{2}s+\mfrac{1}{2}\alpha) +\zeta(s+\alpha)\\
    &=&
    s(\mfrac{1}{2} \log (s) + O(1)),
\end{array}
$$
d'o\`u $\ds\lim_{\sigma \to \infty} h_\alpha(\sigma) =  +\infty$.

En plus, pour $|s|$ suffisament grand
\begin{equation} \label{nano}
F_\alpha(s) = (2\pi)^\alpha s^{-\alpha} \Bigl( 1 + O\bigl( |s|^{-1} \bigr) \Bigr) \Frac{\zeta(s-\alpha)}{\zeta(s+\alpha)}.
\end{equation}
Donc, pour $\sigma$ assez grand
$$ F_\alpha(s) =  O(\sigma^{-\alpha}) $$
uniform\'ement en $\tau$, la condition {\it(\ref{princcc1})} du th\'eor\`eme \ref{teoalin}.

Pour $\tau$ suffisament grand, on a $\zeta(s) = O\bigl(|\tau|^{A}\bigr)$ uniform\'ement sur un demi-plan quelconque
$\sigma \geq \sigma_\alpha$ fix\'e \cite[\S 5.1]{ztitchmarsh}
(avec $A=A(\sigma_\alpha)$), dans notre cas $\sigma_\alpha = 1/2-\alpha$.
En plus, de (\ref{nonannul}), pour $\sigma \geq 1/2$, $\alpha \geq 1/2$,
$$ \Frac{1}{\zeta(s+\alpha)} = O\bigl( \log |\tau| \bigr).$$
Donc, pour $\sigma \geq 1/2$, on a $\sigma-\alpha \geq \sigma_\alpha$ et
$\sigma+\alpha \geq 1$, d'o\`u
$$ F_\alpha(s) = O\bigl(|\tau|^{A}\log |\tau| \bigr), $$
uniform\'ement en $\sigma \geq 1/2$, la condition {\it(\ref{princcc2})} du th\'eor\`eme \ref{teoalin}.

La fonction $f_\alpha(s)$ est m\'eromorphe, avec des p\^oles simples aux
points $s=-\alpha,1-\alpha,\alpha,1+\alpha$. Le point
$s=\frac{1}{2}$ est un z\'ero simple de $f_\alpha(s)$, lorsque $\alpha \neq \alpha^*$,
et pour $\alpha = \alpha^*$ un z\'ero triple. En fait,
$f_\alpha'(\frac{1}{2}) = 2 {\zeta^*}'(\frac{1}{2}+\alpha) \neq 0$ pour $\alpha \neq \alpha^*$.

On a vu que $\ds\lim_{\sigma \to +\infty} h_\alpha(\sigma) = +\infty$, et
$\ds\lim_{\sigma \to +\infty} F_\alpha(\sigma) = 0$, donc
$\ds\lim_{\sigma \to +\infty} f_\alpha(\sigma) = +\infty$. En plus
$$
\begin{array}{rcl}
f_\alpha(1+\alpha+0)
    &=& \zeta^*(1+2\alpha) - \ds\lim_{h \to 0^+}  \zeta^*(1+h) \\
    &=& \zeta^*(1+2\alpha) - \ds\lim_{h \to 0^+}  \Frac{1}{(1+h)-1} =
-\infty. \\
\end{array}
$$
En cons\'equence, il existe un z\'ero de $f_\alpha(s)$, $\rho_\alpha \in ]1+\alpha,+\infty[$ de multiplicit\'e impaire, et $n_{f,\sigma > \frac{1}{2}} \geq 1$.

Il nous reste \`a voir que $ p_{f,\sigma > 1/2} = 2$, $p_{h,\sigma >
1/2} = 0$, et $n_{h,\sigma > 1/2} = 0$. Donc du th\'eor\`eme
\ref{teoalin}
\begin{equation} \label{amejorar}
N(T) - N_0(T) \leq N(T) - N_0'(T) \leq \mfrac{3}{2},
\end{equation}
d'o\`u $N(T) = N_0(T)$ (et $N_0(T) \leq N_0'(T) +1$)  pour tout $T >
0$. On d\'eduit de cette borne que $\rho_\alpha$ est forc\'ement simple et
$n_{f,\sigma > \frac{1}{2}} = 1$.\lqqd \\

En faisant $\alpha \to \frac{1}{2}$, le th\'eor\`eme de Hurwitz nous donne le r\'esultat
suivant, sous la forme originalement pos\'ee par P.R. Taylor.

\begin{prop}[P.R. Taylor, {\cite[\S 1]{taylor}}] \label{prtaylor} Tous les z\'eros non r\'eels de
la fonction
$$ f_{1/2}(s) = \zeta^*(s+\mfrac{1}{2}) - \zeta^*(s-\mfrac{1}{2}) $$
se trouvent sur la droite $\sigma = \frac{1}{2}$.
\end{prop}

L'\'equation (\ref{amejorar}) nous dit que, pour $\alpha > 1/2$, tous les z\'eros non r\'eels
de $f_\alpha(s)$ sont simples et sur la droite $\sigma=\frac{1}{2}$, sauf probablement deux
(un et son conjugu\'e) z\'eros doubles.
On peut utiliser le lemme \ref{sobrelarectauno} pour am\'eliorer ce r\'esultat.
On calcule
$$ r(\alpha)= \frac{h_\alpha'(\frac{1}{2})}{h_\alpha(\frac{1}{2})} = \frac{(\zeta^*)'(\frac{1}{2}+\alpha)}{\zeta^*(\frac{1}{2}+\alpha)} $$
suivant la valeur de $\alpha > \frac{1}{2}$ (voir la figure \ref{dalpha}). Pour $\alpha \to \frac{1}{2}^+$, cette valeur
approche $-\infty$.
\begin{figure}[tbh]
\unitlength=1cm
\begin{center}
\psfrag{10}[tl][tc][0.7]{{$10$}}
\psfrag{-1}[tr][tc][0.7]{$-1$}
\psfrag{1}[br][bc][0.7]{$1$}
\psfrag{-2.6205}[l][][0.7]{\begin{picture}(0,0)
    \put(6.9,-0.3){\makebox(0,0)[l]{$\alpha^*$}}
    \put(0,5){\makebox(0,0)[l]{$0$}}
\end{picture}}
\psfrag*{-2.6205}[tr][tc][0.7]{$-2.62$}
\includegraphics[width=7.25cm,height=5.05cm]{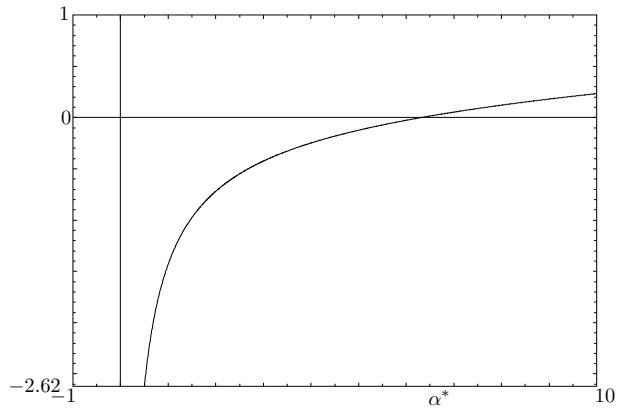}
\end{center}
\caption{La fonction $r(\alpha)$}
\label{dalpha}
\end{figure}
Lorsque $\alpha < \alpha^*$, par le lemme \ref{sobrelarectauno},
l'\'ecart $d_{f_\alpha,a}$ dans le th\'eor\`eme \ref{enezero}
satisfait $d_{f_\alpha,a} \geq 1$ et
$$ N_0'(T) \geq \mfrac{1}{\pi} \arg h(\mfrac{1}{2}+iT), $$
ce qui entra\^ine la constante $1/2$ au lieu de $3/2$ dans (\ref{amejorar}). Donc
tous les z\'eros non-r\'eels de $f_\alpha(s)$ sont simples et align\'es sur la droite
$\sigma = \frac{1}{2}$, lorsque $\alpha < \alpha^*$. La fonction $f_{\alpha^*}(s)$ a aussi
tous ses z\'eros simples et align\'es, sauf $s = \frac{1}{2}$, qui est un z\'ero triple.
Et pour $\alpha > \alpha^*$, malheureusement on ne peut pas obtenir mieux que (\ref{amejorar}).

Bien que le passage de $\alpha > \frac{1}{2}$ \`a $\frac{1}{2}$ nous fournisse la proposition \ref{prtaylor}, la
manque de stabilit\'e (par rapport au param\`etre $\alpha$) dans le cas g\'en\'eral de $N_0'(T)$ nous laisse sans information
\`a propos de la simplicit\'e des z\'eros de $f_{1/2}(s)$.
On remarque que la m\'ethode originale de P. R. Taylor ne peut pas
\^etre appliqu\'ee \`a l'\'etude de $f_{1/2}(s)$;
il est pourtant facile de la modifier pour le faire.
Dans notre cas, on peut appliquer le th\'eor\`eme \ref{teoalin} \`a la fonction
$$ \hat{f}(s) = \bigl( s - \mfrac{1}{2} \bigr) f_{1/2}(s) = \bigl( s - \mfrac{1}{2} \bigr) \zeta^* \bigl( s + \mfrac{1}{2} \bigr) +
\bigl( \mfrac{1}{2}-s  \bigr) \zeta^* \bigl(\mfrac{3}{2}+s \bigr) $$
et obtenir que presque tous les z\'eros de $f_{1/2}(s)$ sont simples, sauf probablement deux (un et son conjugu\'e, qui pourraient avoir multiplicit\'e double).

Pour la fonction
$$ \hat{f}_\alpha(s) = \zeta^*(s+\alpha) + \zeta^*(s-\alpha) $$
on a
\begin{equation} \label{imposiblemejorar}
0 \leq N(T) - N_0(T) \leq \mfrac{5}{2},
\end{equation}
(pas de z\'ero r\'eel connu pour $\hat{f}_\alpha(s)$), ce qui ne nous permet pas de montrer que tous les z\'eros de $\hat{f}_\alpha(s)$
sont align\'es. L'am\'elioration de l'estimation pr\'ec\'edente peut
\^etre obtenue en estimant la constante $d_{\hat{f},a}$.
Il faut remarquer que le lemme \ref{sobrelarectauno} ne peut pas \^etre appliqu\'e \`a cette fonction.

D'abord on analyse l'argument de $h_\alpha(s)$ sur la droite critique, lorsque $\alpha$ est proche de $\frac{1}{2}$.
Pour $\alpha = \mfrac{3}{5}$, le graphe de la fonction (figure \ref{udet})
$$ u(\tau) = \mfrac{1}{\pi} \arg h_{3/5}(\mfrac{1}{2}+i\tau) - \mfrac{1}{2}$$
\begin{figure}[tbh]
\unitlength=1cm
\begin{center}
\psfrag{21}[tl][tc][0.7]{{$21$}}
\psfrag{-10}[tr][tc][0.7]{$-10$}
\psfrag{-1.3805}[br][bc][0.7]{$-1.38$}
\psfrag{0.38053}[br][bc][0.7]{$0.38$}
\psfrag*{-1.3805}[l][][0.7]{\begin{picture}(0,0)
    \put(0,5.4){\makebox(0,0)[r]{$0$}}
    \put(0,1.4){\makebox(0,0)[r]{$-1$}}
\end{picture}}
\includegraphics[width=7.25cm,height=5.05cm]{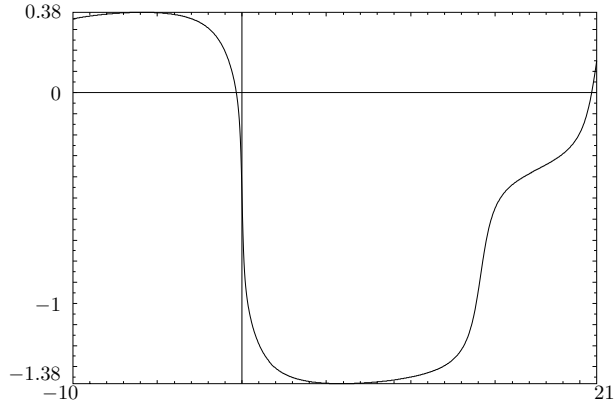}
\end{center}
\caption{La fonction $u(\tau)$}
\label{udet}
\end{figure}
nous montre que les trois premiers points $0 < x_1 < x_2 < x_3 < 21$ \`a valeurs enti\`eres de $u(\tau)$  satisfont
$u(x_1)= -1 = u(x_2)$, $u(x_3) = 0$, d'o\`u $\tilde{u}(T)$ d\'efini comme dans le lemme
\ref{compcreciente} satisfait $\tilde{u}(T) \geq u(T)+2$. Donc, l'\'ecart dans le th\'eor\`eme \ref{enezero} satisfait  $d_{f_\alpha,a} \geq 2$ et
$$ N_0'(T) = \tilde{u}(T) \geq u(T)+ 2  = \mfrac{1}{\pi} \arg h(\mfrac{1}{2}+iT) + \mfrac{3}{2}, $$
puis
$$ 0 \leq N(T)-N_0'(T) \leq \mfrac{1}{2}, $$
d'o\`u $N(T) = N_0(T) = N_0'(T)$ et tous les z\'eros de la fonction
$$ \hat{f}_{3/5}(s) = \zeta^*(s+\mfrac{3}{5}) + \zeta^*(s-\mfrac{3}{5}) $$
sont sur la droite $\sigma=\frac{1}{2}$ et ils sont simples.

Pourtant, pour des grandes valeurs de $\alpha$, l'argument de $h_\alpha(a+i\tau)$ est
une fonction croissante de $\tau>0$, et il est impossible d'am\'eliorer (\ref{imposiblemejorar})
avec nos estimations. Si on prend $\alpha$ sufisamment grand, par exemple $\alpha = 8$,
la fonction $\hat{f}_8(s) = \zeta^*(s+8)+\zeta^*(s-8)$ a un z\'ero en $ s_0 \approx 8.78369+ 1.00496 i$, \`a l'exterieur de
la droite critique, et un autre sym\'etrique $1-\ovl{s_0}$ dans le demi-plan sup\'erieur. Donc $N(T)-N_0(T) = 2$,
et l'estimation (\ref{imposiblemejorar}) est optimale.
On conclut qu'il y a exactement $4$ z\'eros hors la droite critique dans ce cas, tous simples.
En plus, $N_0(T)=N_0'(T)$, donc tous les z\'eros de $\hat{f}_8(s)$ sur la droite critique sont simples.

Il y a un z\'ero double de $f_\alpha(s)$ pour $\alpha \approx 2,61117$, en $s_0 = \frac{1}{2}+ i \tau_0 \approx \frac{1}{2}+5,4983 i$. Ceci peut \^etre d\'etect\'e facilement en remarquant que le z\'ero $\frac{1}{2}+i \tau_0$ de $f_\alpha(s)$
est un z\'ero double si et seulement si $u_\alpha'(\tau_0)=0$, o\`u
$$ u_\alpha(\tau) = \mfrac{1}{\pi} \arg h_{\alpha}(\mfrac{1}{2}+i\tau) - \mfrac{1}{2}.$$
Dans notre cas, cela correspond \`a un minimum de $u_\alpha(\tau)$.

On peut aussi faire l'analyse de $f_\alpha(s) $ sous l'hypoth\`ese de
Riemann pour $0 < \alpha < \frac{1}{2}$.

\begin{prop} Soit $0 < \alpha < \frac{1}{2}$. Sous l'hypoth\`ese de Riemann,
tous les z\'eros non r\'eels de
$$ f_{\alpha}(s) = \zeta^*(s+\alpha) - \zeta^*(s-\alpha) $$
sont sur la droite $\sigma= \frac{1}{2}$ et ils sont simples. Cette fonction a exactement
 trois z\'eros r\'eels simples $s=\frac{1}{2}$, $s=\rho_\alpha$ et $s=1-\rho_\alpha$,
o\`u $\rho_\alpha \in ]1+\alpha,+\infty[$.
\end{prop}

{\bf D\'emonstration.}
En effet, $n_{h,\sigma>
\frac{1}{2}} = 0$ et $p_{h,\sigma > \frac{1}{2}} = 1$ (on compte le
p\^ole $s=1-\alpha$ de $h_\alpha(s)$); de plus, $h_\alpha(s)$ ne s'annulle pas pour $\sigma \geq \frac{1}{2}$ sous
l'hypoth\`ese de Riemann, et
\begin{equation} \label{chigau}
\Frac{1}{\zeta(s+\alpha)}  = O\bigl(\log |\tau|\bigr)
\end{equation}
uniform\'ement pour $\sigma \geq \frac{1}{2}$, $\alpha>0$ fix\'e
\cite[Theorem 14.2]{ftitchmarsh}). Donc
$$ N(T) - N_0(T) \leq N(T)-N_0'(T) \leq \mfrac{1}{2}, $$
d'o\`u tous les z\'eros non r\'eels de $f(s)$ satisfont $\sigma =
\frac{1}{2}$ et ils sont simples. \lqqd \\

Pour obtenir des r\'esultats inconditionnels pour $0< \alpha< 1/2$, il y a deux ingr\'edients importants \`a fournir:
d'un c\^ot\'e, une borne de $\zeta(s)^{-1}$ pour $\sigma \geq 1/2$, diff\'erente de (\ref{chigau});
d'un autre c\^ot\'e,
une estimation pour le nombre de z\'eros de la fonction z\^eta dans le demi-plan $\sigma \geq 1/2+\alpha$,
supposer qu'on a un nombre fini de z\'eros dans un tel demi-plan est une condition trop forte
(cela entra\^ine l'existence d'un demi-plan $\sigma \geq \beta > 1/2$ libre de z\'eros de la fonction z\^eta).
On d\'ev\'eloppera ces outils aux \S\S \ref{inficritica}, \ref{infidensidad}.

Les fonctions $f_\alpha^{\pm}(s)=\zeta^*(s+\alpha)\pm \zeta^*(s-\alpha)$ ont \'et\'e r\'ecemment \'etudi\'ees par H. Ki dans \cite{hki}.
La m\'ethode de Titchmarsh dans \cite[\S\S 9.3, 9.4]{ztitchmarsh} donne la formule pour
le nombre de z\'eros de $f_\alpha^{\pm}(s)$ avec $0 < \tau < T$
$$ N(T) = \frac{T}{2\pi} \log \frac{T}{2\pi}  - \frac{T}{2\pi} + O(\log T), $$
o\`u la constante dans $O(\cdot)$ d\'epend de $\alpha > 0$ (la formule dans \cite[Proposition 2.2]{hki}
pour $N(T)$ explicite la d\'ependance de $\alpha$).

Comme il est naturel, H. Ki \'etudie s\'epar\'ement les cas $0 < \alpha <1/2$ et $\alpha \geq 1/2$.
On remarque une faute dans son article, dans l'\'etude du cas $\alpha \geq 1/2$, pour
confirmer
que ses r\'esultats sont coh\'erents avec les notres.
H. Ki ne s'aper\c{c}oit pas de l'existence d'un z\'ero r\'eel de $f_{\alpha}^{-}(s)$ (trouv\'e
dej\`a par P. R. Taylor dans la preuve de la proposition \ref{prtaylor}), ni des z\'eros non r\'eels et hors la droite critique de $f_\alpha^+(s)$.
Dans les derni\`eres lignes de son article, il montre que la fonction
$$ (s+\alpha)(s+\alpha-1)(s-\alpha)(s-\alpha-1) f_\alpha(s) $$
a au plus $4$ z\'eros hors la droite critique (ce qui co\"incide avec notre r\'esultat),
puis il factorise le facteur $(s+\alpha)(s+\alpha-1)(s-\alpha)(s-\alpha-1)$ en consid\'erant qu'il contient exactement $4$ z\'eros,
pour conclure que tous les z\'eros de $f_\alpha(s)$ sont sur la droite critique. Ce dernier raisonnement est erron\'e,
vu que ces facteurs correspondent aux p\^oles de $f_\alpha(s)$, et donc ne peuvent pas \^etre compt\'es commes des z\'eros de la fonction. La simplicit\'e des z\'eros n'est pas trait\'ee.

Pour $0 < \alpha < 1/2$, H. Ki localise les z\'eros de $f_{\alpha}^{\pm}(s)$ dans une r\'egion qui d\'epend
des possibles z\'eros de la fonction z\^eta dans le demi-plan $\sigma > 1/2+\alpha$ \cite[Theorem B]{hki}; il \'etudie aussi les cons\'equences de l'hypoth\`ese de Lindel\"of sur cette r\'egion \cite[Theorem C]{hki}.
Cela implique un r\'esultat de densit\'e z\'ero pour les z\'eros hors la droite $\sigma = 1/2$, dont on
en parlera au \S \ref{infidensidad}.

On peut modifier la premi\`ere fonction consid\'er\'ee et obtenir des r\'esultats similaires. On prend
$h_\alpha(s)= \xi(s+\alpha)$ avec $\alpha > 0$, la fonction $h_\alpha(s)$ est enti\`ere et on peut appliquer
le corollaire \ref{coralin}. Les estimations asymptotiques de $\xi(s+\alpha)$ sont les m\^emes
que celles de $\zeta^*(s+\alpha)$, le facteur $\frac{1}{2}(s+\alpha)(s+\alpha-1)$
ne change pas le comportement de la fonction. L'absence de p\^oles de cette fonction
simplifie les r\'esultats.

Les fonctions $\tilde{f}_\alpha(s)=\xi(s+\alpha) \pm \xi(s-\alpha)$,
ont donc tous leurs z\'eros simples, align\'es sur la droite critique $\sigma = \frac{1}{2}$,
pour $\alpha \geq \frac{1}{2}$ inconditionnellement, et sous l'hypoth\`ese de Riemann pour $0 < \alpha < \frac{1}{2}$
(une autre preuve de ce r\'esultat peut \^etre vue dans l'annexe \ref{seccionestabilidad}).
Cette fonction a \'et\'e \'etudi\'ee par Lagarias dans \cite{lagarias}, qui donne aussi
l'estimation correspondante de $N(T)$.  Lagarias
\'etudie aussi la distribution limite des espacements entre les z\'eros consecutifs des fonctions $\tilde{f}_\alpha(s)$,
et la relation avec la conjecture GUE pour la fonction z\^eta de Riemann.

\section{Infinit\'e de z\'eros sur la droite critique} \label{inficritica}

L'existence d'une quantit\'e finie de z\'eros de la fonction $h(s)$ sur la droite $\sigma = a$
n'est pas une barri\`ere pour l'application du th\'eor\`eme \ref{teoalin} (remarque \itemref{estata} au th\'eor\`eme) et
ne change pas le r\'esultat du th\'eor\`eme \ref{teoalin}. On suppose
maintenant qu'il y a une infinit\'e de z\'eros sur la droite critique $\sigma = a$.
Il est \'evident qu'on ne peut plus utiliser le lemme \ref{rectacritica}
pour essayer d'\'etendre le th\'eor\`eme \ref{teoalin}.
On obtient la g\'en\'eralisation naturelle de ce lemme avec la technique de N. Levinson dans \cite{levinson}.

Introduisons d'abord quelques notations. Pour une fonction $h(s)$ et $\sigma_0<\sigma_1$, $T>0$, on note
$$ N_h(\sigma_0,\sigma_1,T) = \# \{ s \in \ce \tq h(s)=0,  \, \sigma_0 \leq \sigma \leq \sigma_1, \,0 < \tau < T \}, $$
et $ N_h(\sigma_0,T) = N_h(\sigma_0,+\infty,T) $.

\begin{lema} \label{zerospormontones} Soit $a \in \re$, $h(s)$ une fonction m\'eromorphe, r\'eelle sur la droite r\'eelle, avec un
nombre fini de p\^oles, $h(a) \neq 0$, holomorphe sur la droite $\sigma = a$, et holomorphe
et non nulle sur la droite $\sigma  = \sigma_0$, $\sigma_0 > a$.
On d\'efinit, pour $T>0$, $R(\sigma_0,T)$ par
$$R(\sigma_0,T) = \frac{1}{\pi} \Delta\arg h(s),$$
o\`u $\Delta$ indique la variation d'argument de $h(s)$, obtenue par variation continue
le long des segments qui joignent (dans l'ordre indiqu\'e) $\sigma_0$, $\sigma_0+iT$ et $a+iT$.
On d\'enote par $n_{h,a<\sigma<\sigma_0}$ le nombre de z\'eros r\'eels de $h(s)$ dans l'intervalle
$a < \sigma < \sigma_0$ et par $P_{h,a<\sigma<\sigma_0}$ le nombre de p\^oles de $h(s)$
avec $a < \sigma < \sigma_0$.

Le nombre de z\'eros de $f(s)$ sur la droite critique $\sigma=a$ avec $0 < \tau < T$ est minor\'e par
\begin{equation}\label{cotainfncerot}
N_0'(T) \geq R(\sigma_0,T) + (-n_{h,a<\sigma<\sigma_0}+P_{h,a<\sigma<\sigma_0}) - 2 N_h(a,\sigma_0,T)- 2.
\end{equation}
\end{lema}
{\bf Remarque.} Si les z\'eros de $h(s)$ sont restreints \`a une bande finie $|\sigma-a|<\sigma_0$, avec
les notations du th\'eor\`eme \ref{teoalin}, (\ref{cotainfncerot}) devient
\begin{equation}\label{cotainfncerotmod}
N_0'(T) \geq R(T) + (-n_{h,\sigma>a}+P_{h,\sigma>a}) - 2 N_h(a,T)- 2.
\end{equation}
{\bf D\'emonstration.} Soit $T> 0$ fix\'e  tel que $h(s)$ n'a pas de z\'ero sur  la droite $\tau=T$,
et $a_1, a_2, \dots, a_M$ les diff\'erents z\'eros de $h(s)$ sur le segment qui joint $a$ et $a+iT$: $\Re(a_i) = a$ et
$$ 0 < \Im(a_1) < \Im(a_2) < \dots < \Im(a_M)<T.$$
On consid\`ere le rectangle modifi\'e $C_\eps$, form\'e par d\'eformation du rectangle
$$ C: \quad a \leq \sigma \leq \sigma_0, \quad 0 \leq \tau \leq T $$
par des demi-cercles de rayon $\eps > 0$ qui laissent \`a l'ext\'erieur les z\'eros et les p\^oles de $h(s)$
sur la fronti\`ere du rectangle non modifi\'e $C$. \\


Sur ce rectangle, en partant de la valeur en $s=\sigma_0$, l'argument de $h(s)$ est bien d\'efini.
On pose les quantit\'es
$$
\begin{array}{ll}
H_0 = \arg h \bigl( a+i(a_1-0) \bigr) - \arg h(a) \\
H_j = \arg h \bigl( a+i(a_{j+1}-0) \bigr) - \arg h \bigl( a+i(a_j+0) \bigr), & 1 \leq j \leq M-1,\\
H_M = \arg h(a+iT) - \arg h \bigl( a+i(a_M+0) \bigr).
\end{array}
$$
En appliquant le principe de l'argument sur $C_\eps$ et en prenant la limite lorsque $\eps \to 0^+$ on obtient
$$
\begin{array}{l}
\pi (-n_{h,a<\sigma<\sigma_0}+p_{h,a<\sigma<\sigma_0}) + \pi R(\sigma_0,T) - \summ_{j=0}^M H_j  \\
\hspace{2cm} - \pi \bigl( N_h(a,\sigma_0,T) - N_h(a+0,\sigma_0,T) \bigr) =
2 \pi \bigl( N_h(a+0,\sigma_0,T) - p_{h,a<\sigma<\sigma_0}'\bigr),
\end{array}
$$
o\`u $p_{h,a<\sigma<\sigma_0}$ est le nombre de p\^oles r\'eels, $p_{h,a<\sigma<\sigma_0}'$
le nombre de p\^oles non r\'eels de partie imaginaire strictement positive de $h(s)$ avec $a < \sigma <\sigma_0$,
de fa\c{c}on que $p_{h,a<\sigma<\sigma_0}+2p_{h,a<\sigma<\sigma_0}' = P_{h,a<\sigma<\sigma_0}$.

On r\'e\'ecrit la formule comme
$$
\begin{array}{rcl}
\summ_{j=0}^M H_j &=& \pi(-n_{h,a<\sigma<\sigma_0}+P_{h,a<\sigma<\sigma_0}\bigr) + \pi R(\sigma_0,T) \\
    & & \pi \bigl( N_h(a,\sigma_0,T) - N_h(a+0,\sigma_0,T) \bigr) - 2 \pi N_h(a,\sigma_0,T).
\end{array}
$$
Dans chaque intervalle (ouvert) correspondant \`a $H_j$, il y a au moins $H_j/\pi - 2$ z\'eros de $f(s)$
(par un analogue au lemme \ref{rectacritica}),
ce qui rajout\'e aux $M$ z\'eros connus de $f(s)$ (les z\'eros de $h(s)$ sur la droite critique),
nous donne
$$
\begin{array}{rcl}
N_0'(T) &\geq& M + \summ_{j=0}^M \Bigl( \Frac{H_j}{\pi} - 2 \Bigr) =
-M - 2 + \Frac{1}{\pi}  \summ_{j=0}^M H_j  \\
    & = & -M-2 + ( -n_{h,a<\sigma<\sigma_0}+P_{h,a<\sigma<\sigma_0}) + R(\sigma_0,T)  \\
    &   & + \bigl(N_h(a,\sigma_0,T)-N_h(a+0,\sigma_0,T)\bigr) - 2 N_h(a,\sigma_0,T).
\end{array}
$$
Ceci avec l'estimation
$$ N_h(a,\sigma_0,T) - N_h(a+0,\sigma_0,T) \geq M $$
nous donne le r\'esultat attendu. \lqqd

\section{Infinit\'e de z\'eros \`a droite de la droite critique et th\'eor\`emes de densit\'e} \label{infidensidad}

Consid\'erons pour l'instant les translat\'ees de la fonction z\^eta de Riemann \'etudi\'ees dans le \S \ref{aquiestimation}.
L'\'etude des fonctions telles que $f_\alpha(s)=\zeta^*(s+\alpha) \pm \zeta^*(s-\alpha)$ pour $0<\alpha < 1/2$
avec le th\'eor\`eme \ref{teoalin} est possible sous l'hypoth\`ese de Riemann ou sous l'hypoth\`ese
de finitude du nombre de z\'eros de la fonction z\^eta dans le demi-plan $ \sigma \geq \frac{1}{2}+\alpha$,
ce qui est une condition trop forte \`a demander. Le meilleur qu'on a ce sont les r\'esultats de densit\'e z\'ero,
c'est-\`a-dire, des estimations de type $ N_\zeta(\beta,T) = O(T^\theta) $, $0 \leq \theta<1$
pour le nombre de z\'eros de la fonction z\^eta avec $\sigma\geq \beta$, $0 <\tau< T$, lorsque $\beta > \frac{1}{2}$.
Mais un r\'esultat de ce type indique qu'un 0\% des z\'eros de la fonction $\zeta^*(s+\alpha)$
sont \`a droite de la droite critique $\sigma=1/2$. L'analogue du th\'eor\`eme \ref{teoalin}
devrait s'\'enoncer comme suit: \og un 0\% des z\'eros de $f_\alpha(s)$ sont
en dehors de la droite critique $\sigma=1/2$ \fg.
 En particulier, le 100\% des z\'eros de $f_\alpha(s)$ devraient \^etre sur la droite $\sigma=1/2$!

On revient au cas g\'en\'eral. Pour simplifier la notation, nous allons \'etudier seulement les fonctions $f(s)= h(s) \pm h(2a-s)$ dont la fonction $h(s)$ a les z\'eros inclus dans
une bande $|\sigma-a|<\sigma_0$.

\begin{teor} \label{teoalindensidad} Soit $h(s)$ une fonction m\'eromorphe sur $\ce$, r\'eelle sur la droite r\'eelle,
qui ne s'annulle pas pour $\sigma$ assez grand, avec un nombre fini de p\^oles.
Avec les notations du th\'eor\`eme
\ref{teoalin}, que la fonction $F(s)$ satisfait la condition {\it(\ref{princcc1})}, et
\begin{enumerate}
\item[\hypertarget{princcc2p}(ii')]  il existe une fonction croissante $\phi: \re \to \re$ telle que,
pour chaque $\sigma_0 > a$, il existe une constante $K>0$ et une suite $(T_n)_n$ telles que $\ds\lim_{n \to \infty} T_n = +\infty$,
$$ T_n \leq T_{n+1} \leq \phi(T_n) \mbox{ pour } n \geq 1, $$
et
$$ \bigl| F(s) \bigr| < e^{K|s|} \mbox{ pour } a \leq \sigma \leq \sigma_0, |\tau|=T_n, n\geq 1 . $$
\end{enumerate}
Alors, pour $T>0$
$$ N(T) - N_0'(T) \leq 4 N_h\bigl(a,\phi(2T) \bigr) + O(1). $$
\end{teor}

{\bf D\'emonstration.} On suit la d\'emarche dans la d\'emonstration du th\'eor\`eme \ref{teoalin}.
Sans perte de g\'en\'eralit\'e, on peut supposer $h(s)$ holomorphe sur la droite critique $\sigma=a$
(remarque (\ref{estata}) au th\'eor\`eme \ref{teoalin}).
Soit $\sigma_0$ assez grand pour satisfaire {\it(\ref{princcc1})} avec $\eps=\frac{1}{2}$, au m\^eme
temps que la bande $2a-\sigma_0\leq \sigma \leq \sigma_0$ contient tous les z\'eros et p\^oles de $h(s)$ et $f(s)$.
D'un autre c\^ot\'e, on prend $K>0$ et la suite $(T_n)_n$ associ\'es \`a $\sigma_0$ en {\it(\hyperlink{princcc2p}{ii'})}.

On calcule la variation d'argument de $f(s)$ le long du rectangle $R_n$ d\'efini par
$$ -T_n \leq \tau \leq T_n, \quad 2a-\sigma_0 \leq \sigma \leq \sigma_0. $$
Le th\'eor\`eme de Littlewood nous donne la formule (\ref{intesedet}), d'o\`u l'on d\'eduit que
$$ \pi \int_0^{T_n} S(\tau) d\tau \leq \sigma_0 K (\sigma_0+T_n) + \frac{1}{2} T_n - I(\sigma_0). $$
Donc, il existe $C=C(\sigma_0)>0$ tel que
\begin{equation} \label{cosuptn}
\Frac{1}{T_n} \int_0^{T_n} S(\tau) d\tau \leq C,
\end{equation}
pour $n$ assez grand, et pas forc\'ement une borne inf\'erieure, sinon on mettrait $O(1)$ au lieu de la borne $C$.

D'un autre c\^ot\'e, le lemme \ref{rectacritica} est remplac\'e par le lemme \ref{zerospormontones},
ce qui nous donne l'estimation (\ref{cotainfncerotmod})
$$ N_0'(T) \geq R(T) - 2 N_h(a,T) + O(1). $$
Puis, pour $0 < T' <T$
$$
\begin{array}{rcl}
(T-T')\bigl( N(T')-N_0'(T') \bigr) &\leq&  \ds\int_0^T \bigl(
N(\tau)-N_0'(\tau)
\bigr) d\tau  \\
  &  \leq & \ds \int_0^{T} 2 N(a,\tau) d\tau + \int_0^T S(\tau) d\tau + O(T)
\end{array}
$$
On prend $T=T_n$, $T'=T_n/2$ pour obtenir
$$ \frac{T_n}{2} \bigl( N(T_n/2)-N_0'(T_n/2) \bigr) \leq 2 N_h(a,T_n) T_n + \int_0^{T_n} S(\tau)d\tau + O(T_n).$$
En divisant par $T_n/2$ et par (\ref{cosuptn})
$$ N(T_n/2)-N_0'(T_n/2) \leq 4 N_h(a,T_n) + O(1).$$
Cette fois-ci on a fait dispara\^itre la constante $C$ vu qu'on
traite des quantit\'es positives.
Soit maintenant $T>0$ suffisamment grand. Alors, il existe $n \geq 1$ tel que $T_{n-1}/2 \leq T \leq T_n/2$, d'o\`u
$$ T_n \leq \phi(T_{n-1}) \leq \phi(2T). $$
La croissance des fonctions $N(T)-N_0'(T)$ et $N_h(a,T)$ (par rapport \`a $T$) nous permettent de conclure. \lqqd

{\bf Remarques.}
\begin{enumerate}
\item La condition {\it(\ref{princcc1})} du th\'eor\`eme \ref{teoalin} apliqu\'ee est beaucoup plus forte de ce dont on a besoin,
elle peut \^etre remplac\'e par une condition plus faible, comme dans le th\'eor\`eme \ref{teoalindom}.
\item Dans la derni\`ere d\'emonstration on \'evite de prendre la limite lorsque $T \to +\infty$;
ce qui est claire si l'on suppose que $ \ds\lim_{T \to \infty} N_h(a,T) = +\infty $ (sinon on tombe
dans le th\'eor\`eme \ref{teoalin}). Ceci explique aussi le terme d'erreur $O(1)$.
L'inter\^et du r\'esultat est le fait de pouvoir calculer l'ordre de croissance du nombre de z\'eros de $f(s)$ hors la droite critique $\sigma=a$.
\item Si $N(a,T) = O(T^A \log^B T)$ et on peut choisir $\phi(T)$ lin\'eaire en $T$, alors
$$N\bigl(a,\phi(T)\bigr) = O(T^A \log^B T).$$
\end{enumerate}

On reprend l'analyse des traslat\'ees de la fonction z\^eta de Riemann. On rappelle pour cela
des propri\'et\'es de la fonction z\^eta:
\begin{lema} \label{riemannpropi} Pour la fonction z\^eta de Riemann $\zeta(s)$, on a:
\begin{enumerate}
\item \label{densizero} 
Pour $\frac{1}{2} \leq \sigma \leq 1$ \cite[p. 275]{ivic}
$$
N_\zeta(\sigma, T) = \left\{
\begin{array}{cl}
O(T^{3(1-\sigma)/(2-\sigma)}\log^5 T), & \mfrac{1}{2} \leq \sigma \leq \mfrac{3}{4} \mbox{ (Ingham)}, \\
O(T^{3(1-\sigma)/(3\sigma-1)}\log^{44} T), & \mfrac{3}{4} \leq \sigma \leq 1 \mbox{ (Huxley)}. $$
\end{array}
\right.
$$
\item \label{valorminim} 
Il existe une constante $A>0$ telle que pour chaque $n \geq 1$, il existe $n<T_n<n+1$
tel que
$$ \bigl| \zeta(s) \bigr|  > \tau^{-A} $$
lorsque $\tau = T_n$, $-1 \leq \sigma \leq 2$ \cite[\S 9.7]{ztitchmarsh}.
\end{enumerate}
\end{lema}

\begin{prop} \label{miomio} Soit $0< \alpha < 1/2$. Pour la fonction
$f_\alpha(s)=\zeta^*(s+\alpha) \pm \zeta^*(s-\alpha)$, on
a
$$ N(T)-N_0'(T) \leq 4 N_\zeta(\mfrac{1}{2}+\alpha,2T+4) +O(1). $$
En particulier,
$$ N(T)-N_0'(T)  =
\left\{
\begin{array}{cl}
O(T^{3(1-2\alpha)/(3-2\alpha)}\log^5 T), & 0 < \alpha \leq \mfrac{1}{4}, \\
O(T^{3(1-2\alpha)/(1+6\alpha)}\log^{44} T), & \mfrac{1}{4} \leq \alpha \leq \mfrac{1}{2}. $$
\end{array}
\right.
$$
\end{prop}

{\bf D\'emonstration.} La condition {\it(\ref{princcc1})} du th\'eor\`eme \ref{teoalin} est toujours satisfaite, mais
pas la condition {\it(\ref{princcc2})}.
On combine le lemme \ref{riemannpropi} {\it(\ref{valorminim})},  et les estimations utilis\'ees dans l'\'etude du cas $\alpha \geq \frac{1}{2}$
dans l'\'equation (\ref{nano}),
pour obtenir une suite $(T_n)$, avec $n < T_n < n+1$ et une constante $B>0$ telle que
la fonction $ F_\alpha(s) = \zeta^*(s-\alpha)/\zeta^*(s+\alpha)$ satisfait
$$ F_\alpha(s) = O(\tau^B) $$
pour $ \sigma \geq \frac{1}{2}$, $\tau=T_n$, $n \geq 1$. En prenant $\phi(T)=T+2$, on a la condition {\it(\hyperlink{princcc2p}{ii'})} du th\'eor\`eme \ref{teoalindensidad}.
Cela nous donne le r\'esultat principal. Les estimations particuli\`eres
sont cons\'equence d'appliquer celles de Ingham et Huxley dans le lemme \ref{riemannpropi} {\it(\ref{densizero})}. \lqqd \\

Le r\'esultat obtenu entra\^ine une densit\'e z\'ero pour les z\'eros des fonctions $f_\alpha(s)$
hors la droite critique $\sigma=\frac{1}{2}$. On peut comparer ce r\'esultat avec celui de H. Ki \cite[Theorem B]{hki}, dont on a parl\'e \`a la fin du \S \ref{aquiestimation}, qui entra\^ine
$$N(T)-N_0(T) = O\bigl(N_\zeta(\mfrac{1}{2}+\alpha+0,T) \log T\bigr).$$
Les estimations obtenues restent les m\^emes pour les fonctions $f(s)=\xi(s+\alpha)\pm \xi(s-\alpha)$, $0 < \alpha < \frac{1}{2}$.

\section{Int\'egrales associ\'ees \`a des s\'eries d'Eisenstein} \label{lweng}

La s\'erie d'Eisenstein non-holomorphe $E(z,s)$ pour le groupe modulaire $\PSL(2,\mathbb{Z})$,
 pour $z=x+iy $ dans le demi-plan sup\'erieur $\mathbb{H} = \{z \in \ce \tq y > 0 \}$
est d\'efinie par
$$ E(z,s) =  \frac{1}{2} \summ_{(m,n)\neq (0,0)} \frac{y^s}{|mz+n|^{2s}}, $$
pour $\sigma > 1$. Cette s\'erie admet un prolongement analytique dans le plan complexe,
sauf aux points $s=0$, $s=1$ o\`u elle a des p\^oles simples, avec r\'esidu $\Res_{s=1} E(z,s) = 3/\pi$, et satisfait \`a l'\'equation fonctionnelle
$$ E^*(z,s) = E^*(z,1-s), $$
o\`u $E^*(z,s)=\pi^{-s} \Gamma(s) E(z,s)$ est la s\'erie d'Eisenstein compl\'et\'ee.
Cette s\'erie contient la fonction z\^eta de Riemann; en fait $ E(i,s) = 2 \zeta(s) L(s,\chi_{-4}) $,
o\`u $\chi_{-4}(n)=(-4/n)$; donc on peut attendre, dans certains cas, une hypoth\`ese de Riemann pour $E(z,s)$.

La s\'erie d'Eisenstein agit comme une forme modulaire de poids z\'ero, puisque
$$ E\Bigl( \frac{az+b}{cz+d},s \Bigr) = E(z,s) \mbox{ pour }
\left[ \begin{array}{cc} a  & b \\ c & d \end{array} \right] \in \PSL(2,\mathbb{Z}). $$
En particulier, $E(z+1,s) = E(z,s)$, et par la p\'eriodicit\'e par rapport \`a la variable $x$,
on a un d\'ev\'eloppement de Fourier de la s\'erie,
la formule de Chowla-Selberg \cite[\S 1]{chowlaselberg}:
$$
\begin{array}{rcl}
E(z,s) &=& \ds \zeta(2s) y^s + \sqrt{\pi} \frac{\Gamma(s-1/2)}{\Gamma(s)} \zeta(2s-1) \\
    & & + 2 \ds\sqrt{y}\frac{\pi^s}{\Gamma(s)} \summ_{n=1}^{\infty} n^{1/2-s} \summ_{d|n} d^{2s-1}
    \int_{-\infty}^{+\infty} e^{-2 \pi ny \cosh \tau} e^{(s-1/2) \tau} d\tau \cdot
    \cos (2\pi n x).
\end{array}
$$
Pour $A>0$ et $s \in \ce$, l'int\'egrale
$$ 2 K_s(2A) = \int_{-\infty}^{+\infty} e^{-2 A \cosh \tau} e^{s \tau} d\tau $$
est une fonction de Bessel. Elle est une fonction enti\`ere de $s$, et satisfait aux \'equations (cf. \cite[p.307]{polya})
\begin{equation} \label{eqfunck}
s K_s(2A) = A \bigl( K_{1+s}(2A) - K_{1-s}(2A) \bigr), \qquad K_s(2A) = K_{-s}(2A) \mbox{ pour } s \in \ce.
\end{equation}
Ces \'equations montrent que la fonction de Bessel correspond bien \`a un cas particulier des fonctions
en \'etude dans notre travail (plus pr\'ecisement, au th\'eor\`eme \ref{perturbationalpha}).  P\'olya
a montr\'e que tous les z\'eros de la fonction de Bessel sont simples, et sur la droite $\sigma = 0$ \cite[\S VI]{polya}.
La fonction diviseur
$$ \sigma_s(n)= \summ_{d|n} d^s = \prod_{p^e \| n} \frac{1-p^{(e+1)s}}{1-p^s} $$
est une fonction enti\`ere, aussi avec tous leurs z\'eros sur la droite $\sigma = 0$ et simples.
Donc, les coefficients non-nuls et non-constants dans la s\'erie de Fourier de $E(z,s)$
ont tous leurs z\'eros sur la droite $\sigma = \frac{1}{2}$.

Par contre, le terme constant du d\'ev\'eloppement de Fourier de $E^*(z,s)$
\begin{equation} \label{cambiando}
a_0(y,s) = \int_0^1 E^*(x+iy)dx =  \zeta^*(2s) y^s + \zeta^*(2-2s) y^{1-s},
\end{equation}
peut avoir deux z\'eros hors la droite critique, (corollaire \ref{mezcla}), mais qu'en nombre fini,
tous les autres restant sur cette droite et simples. La r\'epartition des z\'eros de cette fonction a \'et\'e
\'etudi\'ee par
Hejhal \cite{hejhal}, Suzuki et Lagarias \cite{lagarias2}, et finalement H. Ki \cite{hki4}.

Plus g\'en\'eralement, on peut consid\'erer une int\'egrale de $E^*(z,s)$ sur la surface modulaire $\PSL(2,\mathbb{Z} )\bsh \mathbb{H}$ par rapport
\`a une m\'esure positive, tel que la m\'esure hyperbolique $d\mu(z)=y^{-2}dx dy$
$$ \int_{\mathscr{D}} E^*(z,s) d\mu(z),$$
o\`u $\mathscr{D}$ est le domaine fondamental
$ \mathscr{D} = \bigl\{ z \in \mathbb{H} \tq |z| \geq 1, |x| \leq \frac{1}{2} \bigr\} $.
Des telles int\'egrales aparaissent dans le calcul des int\'egrales des fonctions automorphes
par la m\'ethode de Rankin-Selberg. Soit $F(z)$ une fonction continue sur $\mathbb{H}$, invariante par $\PSL(2,\mathbb{Z})$,
dont le terme constant dans le d\'ev\'eloppement de Fourier est $b_0(y) = \ds\int_0^1 F(x+iy)dx$. L'int\'egrale
de $F(z)$ sur $\mathscr{D}$ peut \^etre r\'ecup\'er\'ee par le calcul du r\'esidu (sous les bonnes conditions de convergence)
$$ \Frac{3}{\pi}  \int_{\mathscr{D}} F(z)d\mu(z) = \Res_{s=1} \int_{\mathscr{D}} F(z) E(z,s) \mu(z). $$
En fait,
$$ \int_{\mathscr{D}} F(z) E(z,s) \mu(z) = \int_0^\infty b_0(y) y^{s-2} dy = R(F,s)$$
(cf. \cite[p. 415]{zagier}), ce qui r\'eduit le calcul d'int\'egration sur le domaine fondamental \`a l'int\'egration
sur une droite. La fonction $R(F,s)$ est appel\'e
la transform\'ee de Rankin-Selberg de $F(z)$; elle est une fonction m\'eromorphe,
et satisfait \`a l'\'equation fonctionnelle $R^*(F,s) = R^*(F,1-s)$, o\`u
$R^*(F,s)= \zeta^*(2s) R(F,s)$. Pour contourner les probl\`emes de convergence qu'on rencontre dans la pratique,
Zagier introduit dans \cite{zagier} une proc\'edure de renormalisation d'int\'egrales, par l'introduction
des domaines tronqu\'es $ \mathscr{D}_T = \bigl\{ z \in \mathbb{H} \tq |z| \geq 1, |x| \leq \frac{1}{2}, y\leq T \bigr\} $.
Comme partie des r\'esultats de Zagier, on retrouve une formule de $R^*(F,s)$ incluant un terme de la forme
\begin{equation} \label{renatural}
H(y;s) = p(s) \zeta^*(2s) y^s + p(1-s) \zeta^*(2-2s) y^{1-s},
\end{equation}
o\`u $p(s)$ est une fonction rationnelle \cite[p. 419]{zagier}.

Certaines de ces int\'egrales ont attir\'e l'attention ces derniers temps.
Lin Weng a introduit r\'ecemment des fonctions z\^eta de rang $n$ non ab\'eliennes, associ\'ees \`a des corps de nombres,
en analogie \`a l'int\'egrale d'Iwasawa-Tate, dans \cite{weng}. On particularise l'\'etude au corps
des nombres rationnels. La fonction z\^eta de Weng de rang $2$ correspond \`a l'int\'egration
de la s\'erie d'Eisenstein compl\'et\'ee sur un ensemble de r\'eseaux, les r\'eseaux semi-stables de $\re^2$.
Un r\'eseau $L=L(z)= \{ mz + n \tq  m, n \in \mathbb{Z} \}$, $z \in \mathscr{D}$ est semi-stable
si, et seulement si $z \in \mathscr{D}_1$ \cite[Example 1.25]{grayson}, et la fonction z\^eta de Weng est
$$ Z_{2,\mathbb{Q}}(s) = \int_{\mathscr{D}_1} E^*(z,s)d\mu(z). $$
Suivant la renormalisation de Zagier, on consid\`ere les int\'egrales, pour $T \geq 1$, donn\'ees par
$$ Z_{2,\mathbb{Q}}^T(s) = \int_{\mathscr{D}_T} E^*(z,s)d\mu(z). $$
La m\'ethode de Rankin-Selberg nous donne alors des formules closes pour les fonctions z\^eta de Weng (cf. \cite[pp. 426--427]{zagier}, \cite[\S 1]{lagarias2}):
$$ Z_{2,\mathbb{Q}}^T(s) = \zeta^*(2s) \frac{T^{s-1}}{s-1} - \zeta^*(2-2s) \frac{T^{-s}}{s}, $$
et en particulier
$$Z_{2,\mathbb{Q}}(s) = \zeta^*(2s) \frac{1}{s-1} - \zeta^*(2-2s) \frac{1}{s}. $$

On veut maintenant \'etudier la r\'epartition des z\'eros de ces fonctions d'une fa\c{c}on unifi\'ee. Les fonctions
z\^eta de Weng et le terme constant de la z\'erie d'Eisenstein sont des cas particuliers de la formule (\ref{renatural}).
Lagarias et Suzuki proposent dans \cite{lagarias2} l'\'etude de cette fonction pour $y \geq 1$, et $p(s)$ est un polyn\^ome \`a coefficients r\'eels,
ce qui \'equivaut \`a l'\'etude avec une fonction rationnelle $p(s)$.

On fait d'abord quelques remarques sur les racines de $p(s)$. Par la remarque (\ref{estata}) du th\'eor\`eme \ref{teoalin},
on peut supposer que $s= \frac{1}{2}$ est une racine simple de $p(s)$, et qu'il n'y a pas d'autres z\'eros pour $p(s)$ sur la droite critique $\sigma = \frac{1}{2}$,
parce que ces z\'eros sont aussi des z\'eros de $H(y,s)$ de la m\^eme multiplicit\'e, et
qui peuvent \^etre facilement factoris\'es. Il faut aussi voir que $H(y;s)$ a un p\^ole en $s=1$ si
et seulement si $p(0) \neq 0$. L'application du th\'eor\`eme \ref{teoalin} et le lemme
\ref{sobrelarectauno} nous fournissent le r\'esultat suivant.

\begin{teor} \label{pabierto} Soit
$$ H^{\pm}(y;s) = p(s) \zeta^*(2s) y^s \pm p(1-s) \zeta^*(2-2s) y^{1-s}, $$
o\`u $y \geq 1$, $p(s)$ un polyn\^ome \`a coefficients r\'eels, tel que $p(s) = (2s-1) q(s)$, o\`u
$q(s)$ est un polyn\^ome sans racine sur $\sigma=\frac{1}{2}$. Alors
\begin{enumerate}
\romano
\item \label{hys1} Le nombre de z\'eros non r\'eels et de partie imaginaire positive de $H^{\pm}(y;s)$ en dehors de la droite critique $\sigma = \frac{1}{2}$
est born\'e par
$$ N(T) - N_0(T) \leq N(T) - N_0'(T) \leq u_{\pm} - n_{H,\sigma > \frac{1}{2}}  - \frac{n_{H,\frac{1}{2}}}{2} + \chi_{H,0} + N_{p,\sigma > \frac{1}{2}} - d_{H,\frac{1}{2}}, $$
o\`u $n_{H,\sigma > \frac{1}{2}}$ note le nombre de z\'eros r\'eels de $H^{\pm}(y;s)$ avec $\sigma > \frac{1}{2}$,
$n_{H,\frac{1}{2}}$ la multiplicit\'e de $s=\frac{1}{2}$ comme z\'ero de $H^{\pm}(y;s)$ ($0$ pour $H^+(y;s)$, $\geq 1$ et impaire pour $H^{-}(y;s)$),
$\chi_{H,0} = 1$ si $p(0) \neq 0$ et $\chi_{H,0}= 0$ si $p(0) = 0$, $N_{p,\sigma > \frac{1}{2}}$ est
le nombre de racines de $p(s)$ avec $\sigma > \frac{1}{2}$, et $d_{H,\frac{1}{2}} \geq 0$ est le nombre entier dans le th\'eor\`eme \ref{enezero}
(\'equation (\ref{sobrelarecta})).
\item \label{hys2} En particulier, si $s=0$ est un z\'ero de $p(s)$, et le seul z\'ero de $p(s)$ avec $\sigma \geq \frac{1}{2}$ est $s= \frac{1}{2}$, alors
tous les z\'eros de $H^{\pm}(y;s)$ sont sur la droite critique $\sigma= \frac{1}{2}$, et ils sont simples.
\item \label{hys3} Soit $y^* = 4 \pi e^{-\gamma-q'(\frac{1}{2})/q(\frac{1}{2})}$. Alors, pour la fonction $H^{-}(y;s)$,
\begin{enumerate}
\item \hypertarget{hys3a} Si $y < y^*$, on a $d_{H,\frac{1}{2}} \geq 1$;
\item \hypertarget{hys3b} Si $y = y^*$, la multiplicit\'e de $s=\frac{1}{2}$ comme racine de $H^-(y;s)$ est $\geq 3$;
\item \hypertarget{hys3c} Si $y > y^*$, il y a un z\'ero $\rho$ de $H^-(y;s)$ dans l'intervalle $]\frac{1}{2},1[$ (et $1-\rho$ dans $(0,\frac{1}{2})$) si
$q(\frac{1}{2}) q(0) > 0$ ($q(\frac{1}{2})$ et $q(0)$ ont le m\^eme signe), donc $n_{H,\sigma > \frac{1}{2}} \geq 1$.
\end{enumerate}
\end{enumerate}
\end{teor}

{\bf Remarque.} La borne obtenue dans {\it (\ref{hys1})} est $\leq \deg p  + 1$, donc uniforme en $y \geq 1$.
De m\^eme pour le nombre de z\'eros r\'eels et pour la constante $d_{H,\frac{1}{2}}$. \\

{\bf D\'emonstration.} Soit $h(y;s) = p(s) \zeta^*(2s) y^s$, alors
$$ H^{\pm}(y;s) = h(y;s) \pm h(y;1-s) $$
et on applique le th\'eor\`eme \ref{teoalin}. Les z\'eros possibles de $h(y;s)$ dans le demi-plan $\sigma \geq 1/2$ proviennent de $p(s)$,
donc ils sont en nombre fini. Les estimations qui entra\^inent les conditions {\it(\ref{princcc1})} et {\it(\ref{princcc2})} du th\'eor\`eme
\ref{teoalin}
sont les m\^emes que celles de la fonction $\zeta^*(s+\frac{1}{2})$ \'etudi\'ee au \S \ref{aquiestimation},
sauf par le facteur d\'ependant du param\`etre $y \geq 1$. Par la formule de Stirling,
la fonction $F(y;s) = {h(y;1-s)}/{h(y;s)}$ satisfait, pour $\sigma \geq 1/2$ et $|s|$ assez grand
$$
\begin{array}{rcl}
F(y;s) &=& \ds y^{1-2s} \frac{p(1-s)}{p(s)} \frac{\zeta^*(2s-1)}{\zeta^*(2s)} \\
    &=& y^{1-2s} \Frac{p(1-s)}{p(s)} \sqrt{\pi} \Frac{\Gamma \bigr(s-\frac{1}{2}\bigr)}{\Gamma(s)} \Frac{\zeta(2s-1)}{\zeta(2s)}\\
    &=& y^{1-2s} (-1)^{\deg p} \Bigl( 1 + O\bigl( |s|^{-1} \bigr) \Bigr) \sqrt{\pi} s^{-1/2} \Bigl( 1 + O\bigl( |s|^{-1} \bigr) \Bigr) \Frac{\zeta(2s-1)}{\zeta(2s)}.
\end{array}
$$
Donc
\begin{equation} \label{cocienteyyy}
{F(y;s)} = y^{1-2s} (-1)^{\deg p} \sqrt{\pi} s^{-1/2} \Bigl( 1 + O\bigl( |s|^{-1} \bigr) \Bigr) \Frac{\zeta(2s-1)}{\zeta(2s)}
\end{equation}
pour $\sigma \geq \frac{1}{2}$ et $|s|$ assez grand. Pour $\sigma \to +\infty$, cette expression est
$O(\sigma^{-1/2})$. Pour $|\tau| \to +\infty$, l'expression est $O\bigl(|\tau|^K\bigr)$
pour une constante $K > 0$. Donc, la fonction $F(y;\cdot)$ satisfait les conditions du th\'eor\`eme \ref{teoalin}, et
on obtient {\it(\ref{hys1})} en rempla\c{c}ant  $P_{H,\sigma > \frac{1}{2}} = \chi_{H,0}$, $N_{H,\sigma > \frac{1}{2}} = N_{p,\sigma > \frac{1}{2}}$
dans (\ref{naa}). La partie ({\it(\ref{hys2})} du th\'eor\`eme est une cons\'equence imm\'ediate de {\it(\ref{hys1})}.

Il nous reste \`a montrer {\it(\ref{hys3})}.
Pour obtenir {\it(\hyperlink{hys3a}{a})} on utilise le lemme \ref{sobrelarectauno}. De
$$h(y;s) = (2s-1) q(s) \zeta^*(2s) y^s = \frac{q(s)}{s} \xi(2s) y^s,$$
on a
$$
\frac{\dds{h}(y;\frac{1}{2})}{h(y;\frac{1}{2})} = \frac{q'(\frac{1}{2})}{q(\frac{1}{2})} - 2 + 2 \frac{\xi'(1)}{\xi(1)} + \log y.
$$
On utilise maintenant la formule
$$ \frac{\xi'(0)}{\xi(0)} = -\frac{\xi'(1)}{\xi(1)} = -\frac{1}{2} \gamma - 1 + \frac{1}{2} \log 4\pi $$
\cite[\S 11, p. 83]{davenport} pour obtenir
$$
\frac{\dds{h}(y;\frac{1}{2})}{h(y;\frac{1}{2})} = \frac{q'(\frac{1}{2})}{q(\frac{1}{2})} + \gamma - \log 4\pi + \log y = \log y - \log y^*,
$$
donc la condition $y < y^*$ est \'equivalente \`a $ \frac{\dds{h}(y;\frac{1}{2})}{h(y;\frac{1}{2})} < 0$, d'o\`u on obtient {\it(\hyperlink{hys3a}{a})}.
On sait que $\dds{H^-}(y;\frac{1}{2}) = 2 \dds{h}(y;\frac{1}{2})$.
Si $y = y^*$, alors $\dds{H^-}(y;\frac{1}{2}) = 0$, puis $n_{H,\frac{1}{2}} \geq 2$ et puisque ce nombre est impair, on a
$n_{H,\frac{1}{2}} \geq 3$, donc {\it(\hyperlink{hys3b}{b})}. Si $y > y^*$ on ne peut rien dire dans le cas g\'en\'eral.

Sans perte de g\'en\'eralit\'e, on peut prendre $q(\frac{1}{2}) > 0$. Alors
$$h(y;\mfrac{1}{2}) = \ds\lim_{s \to \frac{1}{2}} (2s-1) \zeta^*(2s) q(s) y^s = q(1/2) y^{1/2} > 0.$$
L'analyse ant\'erieure entra\^ine que $\dds{H^-}(y;\frac{1}{2}) =  2 \dds{h}(y;\frac{1}{2}) > 0$ si $y > y^*$. Donc, il existe $x_0> \frac{1}{2}$ tel que $H^{-}(y;x_0) > 0$.
Maintenant, la condition $q(0) > 0$ nous dit que $s=0$ n'est pas racine de $p(s)$, donc $H^-(y;s)$ a un p\^ole en $s=1$,
avec r\'esidu $-q(0) < 0$, d'o\`u $ \ds\lim_{y \to 1^-} H^{-} (y;s) = -\infty $.
Donc, il existe $\rho \in ]x_0,1[$ tel que $H^{-}(y;\rho) = 0$, ce qui montre {\it(\hyperlink{hys3c}{c})}. \hfill{$\blacksquare$}\\

Maintenant on applique le th\'eor\`eme \ref{pabierto} pour obtenir comme des cas particuliers
les r\'esultats de \cite{lagarias2}.

\begin{coro}[Suzuki et Lagarias {\cite[Theorem 2]{lagarias2}}] \label{integralt} Pour chaque $T \geq 1$ fix\'e, la fonction
$$ Z_{2,\mathbb{Q}}^T(s) = \zeta^*(2s) \Frac{T^{s-1}}{s-1} - \zeta^*(2-2s) \Frac{T^{-s}}{s} $$
a tous ses z\'eros simples, et sur la droite $\sigma =\frac{1}{2}$.
\end{coro}

{\bf D\'emonstration.} Dans le th\'eor\`eme \ref{pabierto}, on pose $p(s) = (2s-1) s$. Alors
$$ H^{-}(T;s) = T(2s-1)s(s-1) I(T;s), $$
d'o\`u l'on conclut par la partie {\it(\ref{hys2})}. \hfill{$\blacksquare$}\\

En particulier, pour $T=1$ on obtient l'hypoth\`ese de Riemann pour $Z_{2,\mathbb{Q}}(s)$.

\begin{coro}[Suzuki et Lagarias {\cite[Theorem 1]{lagarias2}}] \label{riemannzetawengdos} La fonction z\^eta non-ab\'elienne de rang $2$
sur $\mathbb{Q}$
$$ Z_{2,\mathbb{Q}}(s) = \zeta^*(2s) \Frac{1}{1-s} - \zeta^*(2-2s) \Frac{1}{s} $$
a tous ses z\'eros simples, et sur la droite $\sigma =\frac{1}{2}$.
\end{coro}

R\'ecemment, Hayashi \cite{hayashi} a obtenu une formule close pour la fonction z\^eta de Weng de rang $2$ d'un corps
de nombres quelconque \cite{hayashi}, la m\^eme formule que pour le corps des nombres rationnels
o\`u on change la fonction z\^eta de Riemann compl\'et\'ee par la fonction de Dedekind compl\'et\'ee
et \`a une constante pr\`es. \`A partir la formule obtenue, il obtient l'hypoth\`ese de Riemann
pour les fonctions de rang $2$ par la m\'ethode de Lagarias de Suzuki,
la variante r\'eelle du th\'eor\`eme \ref{perturbationalpha}.

L'hypoth\`ese de Riemann pour la fonction z\^eta de rang 3
a \'et\'e \'etablie par Suzuki en \cite{suzuki}; on peut aussi appliquer
nos m\'ethodes \`a l'\'etude de cette fonction.
La fonction z\^eta de rang $1$ n'est autre que la fonction z\^eta compl\'et\'ee $\zeta^*(s)$.

On obtient maintenant une version modifi\'ee de \cite[Theorem 3]{lagarias2}:

\begin{coro} \label{mezcla} Soit $y^* = 4 \pi e^{-\gamma}$. Pour chaque $y \geq 1$, la fonction
$$ a_0(y;s) = \zeta^*(2s) y^s + \zeta^*(2-2s) y^{1-s} $$
est une fonction m\'eromorphe, qui peut \^etre prolong\'ee
analytiquement au point $s = \frac{1}{2}$. Alors
\begin{enumerate}
\alfabetico
\item Si $y < y^*$, tous les z\'eros de $a_0(y;s)$ sont sur la droite critique $\sigma=\frac{1}{2}$ et ils sont simples.
\item Si $y = y^*$, tous les z\'eros de $a_0(y;s)$ sont sur la droite critique $\sigma=\frac{1}{2}$ et ils sont simples,
sauf $s=\frac{1}{2}$ qui est un z\'ero double.
\item Si $y > y^*$, il y a exactement deux z\'eros hors la droite critique, un z\'ero r\'eel simple $\rho_y$ dans l'intervalle $]\frac{1}{2},1[$, et son
sym\'etrique $1-\rho_y$. Tous les autres z\'eros sont sur la droite critique et ils sont simples.
\end{enumerate}
\end{coro}

{\bf D\'emonstration.} On prend $p(s) = 2s-1$, d'o\`u
$$ H^{-}(y;s) = (2s-1) a_0(y;s). $$
Comme $p(0) \neq 0$, on ne peut pas appliquer la partie {\it(\ref{hys2})} du th\'eor\`eme \ref{pabierto}.
Par la partie {\it(\ref{hys2})} du m\^eme th\'eor\`eme, pour la fonction $H^{-}(y;s)$ on a
$$
\begin{array}{rcl}
N(T) - N_0(T) &\leq & N(T) - N_0'(T) \leq 2 - N_{H,\sigma > \frac{1}{2}} - \frac{n_{H,\frac{1}{2}}}{2}  - d_{H,\frac{1}{2}},\\
    & \leq & \frac{3}{2} - n_{H,\sigma > \frac{1}{2}} - d_{H,\frac{1}{2}},
\end{array}
$$
puisque $n_{H,\frac{1}{2}} \geq 1$.
On applique maintenant la partie {\it(\ref{hys3})} du th\'eor\`eme. On a $q(s) = 1$, donc $y^* = 4\pi e^{-\gamma}$ et
\begin{enumerate}
\alfabetico
\item Si $ y < y^*$, on a $d_{H,\frac{1}{2}} \geq 1$, puis forc\'ement $n_{H,\sigma > \frac{1}{2}} = 0$ et
$n_{H,\frac{1}{2}} = 1$.
\item Si $y = y^*$, on a $n_{H,\frac{1}{2}} \geq 3$, puis forc\'ement $n_{H,\frac{1}{2}} = 3$ et $n_{H,\sigma > \frac{1}{2}} = 0$.
\item Si $y > y^*$, on a $n_{H,\sigma > \frac{1}{2}} \geq 1$, donc encore $n_{H,\sigma > \frac{1}{2}} = 1$
et $n_{H,\frac{1}{2}} = 1$.
\end{enumerate}
Dans tous les cas, $N(T)=N_0(T)=N_0'(T)$. Cela finit la preuve. \hfill{$\blacksquare$}\\

{\bf Remarques.} Le r\'esultat \cite[Theorem 3]{lagarias2} \'etablit la continuit\'e de $\rho_y$
en fonction du param\`etre $y \geq 1$, en montrant que c'est une fonction croissante de $y$, et que
$\rho_y \to 1$ lorsque $y \to +\infty$. 
L'alignement des z\'eros non-r\'eels de $a_0(y;s)$ a \'et\'e demontr\'e par Hejhal en \cite[Proposition 5.3 (f)]{hejhal} en utilisant
les relations de Maa\ss-Selberg \cite[\S 2.18]{ztitchmarsh}.
La simplicit\'e des z\'eros de $a_0(y,s)$ a \'et\'e r\'ecemment montr\'ee par H. Ki en \cite{hki4}.
On remarque que ces approches traitent l'alignement et la simplicit\'e des z\'eros s\'epar\'ement.

R\'ecemment, W. M\"ueller a donn\'e en \cite{mueller} une interpr\'etation des z\'eros de $a_0(y;s)$
en termes des valeurs propres d'un ope\'rateur autoadjoint, un pseudo-laplacien agissant
sur un sous-espace de $L^2 \bigl(\PSL(2,\mathbb{Z}) \bsh \mathbb{H}\bigr)$.
Dans ce travail, il fournit aussi une extension du corollaire \ref{mezcla} au terme constant d'une s\'erie d'Eisenstein
associ\'ee \`a $\PSL(2,\Og_{\mathbb{K}})$ \cite[Theorem 0.3]{mueller}, o\`u $\Og_{\mathbb{K}}$ est l'anneau des entiers
d'un corps quadratique imaginaire ${\mathbb{K}}=\mathbb{Q}\bigl(\sqrt{-d}\bigr)$,
$d>0$ entier,
tel que le nombre de classes d'id\'eaux entiers $h(-d)=1$. Ce terme est obtenue en rempla\c{c}ant la fonction
z\^eta de Riemann compl\'et\'ee par la fonction z\^eta de Dedekind compl\'et\'ee dans (\ref{cambiando});
nos m\'ethodes s'\'etendent facilement \`a ce cadre.\\

Quant \`a la r\'epartition des z\'eros, une application directe du th\'eor\`eme \ref{enezero}
et la remarque faite au th\'eor\`eme \ref{pabierto} nous fournissent le r\'esultat suivant.

\begin{teor} Pour chaque $y \geq 1$, la fonction $H^{\pm}(y;s)$ dans le th\'eor\`eme \ref{pabierto} a tous ses z\'eros, sauf un nombre fini, simples et align\'es
sur la droite $\sigma = \frac{1}{2}$. Le nombre de ces z\'eros, en fonction du param\^etre $y \geq 1$, avec $0 < \tau < T$ est
$$ N(y;T) = \Frac{T}{\pi}  \log T - \Frac{1}{\pi} (\log \pi + 1) T + \frac{1}{\pi} (\log y) T +  O( \log T).$$
La constante dans $O(\cdot)$ ne d\'epend pas de $y \geq 1$.
Les z\'eros des fonctions $H^{\pm}(y;s)$ sont ultimement entrelac\'es.
\end{teor}

Cette estimation est la m\^eme que celle calcul\'ee par Lagarias et Suzuki \cite[Theorem 5]{lagarias2}
pour la fonction $Z_{2,\mathbb{Q}}^T(s)$. Ils montrent de plus que
la partie imaginaire de chaque z\'ero de $Z_{2,\mathbb{Q}}^T(s)$
 avec $\tau > 0$ est une fonction strictement d\'ecroissante de $T$.

Il nous reste \`a voir le cas o\`u $0 < y < 1$. De (\ref{cocienteyyy}), la fonction $ F(y;s) = \Frac{h(y;1-s)}{h(y;s)} $
satisfait, pour $\sigma$ suffisament grand, $F(y;s)=O\bigl( y^{-2\sigma} |s|^{-1/2} \bigr)$.
Puisqu'elle contient une puissance de $y$ convergeant vers l'infini lorsque $\sigma \to + \infty$,
cette fonction ne satisfait pas les conditions du th\`eor\`eme \ref{teoalin}.

La proposition suivante a \'et\'e \'enonc\'ee par D. Hehjal en \cite[p. 89]{hejhal},
et montr\'ee par H. Ki \cite[\S 2]{hki2} dans le cas o\`u $p(s)=1$.

\begin{prop} \label{hysgg} Pour $ 0 < y < 1$, la fonction $H^{\pm}(y,s)$ du th\'eor\`eme \ref{pabierto} poss\`ede des z\'eros avec parties r\'eelles arbitrairement grandes.
Pour $\sigma_0$ suffisamment grand, le nombre de z\'eros de $H(y,s)$ avec $\sigma \geq \sigma_0$ et $0 < \tau < T$
satisfait
$$ N_{H}(\sigma_0,T) = O(T). $$
\end{prop}

{\bf D\'emonstration.} Les z\'eros de $H(y;s)$ sont, essentiellement (sauf un nombre fini), ceux
de la fonction $g(s) = 1 \pm F(y;s)$. Par l'\'equation (\ref{cocienteyyy}), $g(s) = 1 \pm y^{1-2s} \sqrt{\pi} s^{-1/2} \bigl( 1 + O(\sigma^{-1}) \bigr)
$
pour $\sigma \geq \sigma_0$. Cette fonction a le m\^eme comportement
que la fonction
$$g_1(s) = 1 \pm y^{1-2s} \sqrt{\pi} s^{-1/2},
$$
comparable \`a une exponentielle-polyn\^ome sans terme principal \cite[Chapter 13]{BC1963}. Le principe de l'argument
appliqu\'e \`a $g_1(s)$ montre que son nombre de z\'eros avec $0 < \tau < T$  est $c T + O(1)$ pour $c=-\frac{1}{\pi}  \log y > 0$
\cite[Part III, Chapter 4, Problem 206.2]{polyaszego}. Le th\'eor\`eme de Rouch\'e
transmet cette information \`a $g(s)$ dans le demi-plan $\sigma \geq \sigma_0$.

Or, chaque racine de $g(s)$ v\'erifie $ |s|^{1/2} = y^{1-2 \sigma} \sqrt{\pi} \bigl( 1 + O(\sigma^{-1}) \bigr) $.
Pour une suite de z\'eros de $g(s)$ (vu qu'il y a une infinit\'e), $|s| \to + \infty$, d'o\`u de cette
derni\`ere \'egalit\'e $\sigma \to +\infty$. \lqqd\\

On parlera ensuite des autres z\'eros de $H(y,s)$, ceux dans la bande $|\sigma-a| < \sigma_0$.

\section{Dominance faible. Un r\'esultat non effectif.} \label{faibledd}

La n\'ec\'essit\'e d'\'etendre le th\`eor\`eme \ref{pabierto} au cas $0< y < 1$,
qui n'est pas de la m\^eme nature que le cas $y \geq 1$, nous motive \`a \'etablir le r\'esultat suivant.

\begin{teor} \label{teoalindom} Soit $\sigma_0 > a$, $h(s)$ une fonction m\'eromorphe, r\'eelle sur la droite r\'eelle,
avec un nombre fini de p\^oles, un nombre fini de z\'eros dans la bande  $a < \sigma < \sigma_0$;
holomorphe et non nulle sur
la droite $\sigma=\sigma_0$. On pose
$$ f(s) = f^{\pm}(s) = h(s) \pm h(2a-s). $$
On suppose que la fonction
$$ F(s) = \Frac{h(2a-s)}{h(s)} $$
satisfait
\begin{enumerate}
\item[\hypertarget{princcc1pp}(i'')] $F(s)\neq \pm 1$ sur $\sigma=\sigma_0$, et il existe $\tau_0 >0$ tel que $\bigl|F(s)\bigr|<1$ pour $\sigma=\sigma_0$, $|\tau| \geq \tau_0$;
\item[\hypertarget{princcc2pp}(ii'')] il existent $K > 0$ et une suite $(T_n)_n$ telle que $\ds\lim_{n \to \infty} T_n = +\infty$ et
$|F(s)|< e^{K|s|}$ pour $a \leq \sigma \leq \sigma_0$, $|\tau| = T_n$, $n \geq 1$.
\end{enumerate}

On note le nombre des z\'eros de $f(s)$ avec $0 < \tau < T$  dans
la bande $|\sigma-a| < \sigma_0-a$ (en comptant les multiplicit\'es) par
$$ N(T,\sigma_0) = \# \{ s \in \ce \tq |\sigma-a|< \sigma_0-a, \; 0 < \tau < T \}. $$
Alors, il existe une constante $B=B(\sigma_0)$ telle que
$$ N(T,\sigma_0) - N_0'(T) \leq B $$
pour $T > 0$.
Donc tous les z\'eros de $f(s)$ dans la bande $|\sigma-a| < \sigma_0-a$, sauf un nombre fini, sont
sur la droite $\sigma=a$ et ils sont simples.
\end{teor}

{\bf Remarque.} En comparant avec les conditions du th\'eor\`eme \ref{teoalin}, la condition
{\it(\hyperlink{princcc1pp}{i''})} est l'analogue de {\it(\ref{princcc1})} avec $\eta=1$, alors que
{\it(\hyperlink{princcc2pp}{ii''})} est l'analogue de {\it(\ref{princcc2})} avec un $K>0$ fix\'e au lieu d'un $\eps>0$ arbitraire.

{\bf D\'emonstration.} On suit la d\'emarche de la d\'emonstration du th\'eor\`eme \ref{teoalin}.
D'apr\`es la remarque (\ref{estata}) \`a ce th\'eor\`eme, on peut supposer que $F(s)$
est reguli\`ere et non nulle sur la droite $\sigma=a$.

On consid\`ere le rectangle $R$ d\'efini par
$$ -T \leq \tau \leq T, \qquad 2a-\sigma_0 \leq \sigma \leq \sigma_0, $$
tel que $f(s)$ n'a pas ni de z\'ero ni de p\^ole sur sa fronti\`ere. On calcule $N(T,\sigma_0)$ comme
on le fait pour $N(T)$ dans le th\'eor\`eme \ref{teoalin}, en calculant la variation d'argument des fonctions $g(s)$ et $h(s)$ (avec $f(s)=h(s)g(s)$) le long du
chemin donn\'e par les deux segments qui lient dans l'ordre $\sigma_0$, $\sigma+i T$, $s$
$$ R(T) = \frac{1}{\pi} \Delta \arg h(s), \qquad S(T) = \frac{1}{\pi} \Delta \arg g(s). $$
On remarque que ces deux quantit\'es d\'ependent de $\sigma_0$ en g\'en\'eral, dans ce cas fixe. L'application du lemme de Littlewood
nous fournit la formule
$$ \pi \int_0^{T_n} S(\tau) d\tau = \int_a^{\sigma_0} \log \bigl|g(\sigma+iT_n)\bigr|d\sigma+ \int_0^{T_n} \arg g(\sigma_0+i\tau)d\tau - I(\sigma_0). $$
De {\it(\hyperlink{princcc1pp}{i''})}, pour $\tau \geq \tau_0$, on aura
$ \Re g(s) \geq 1 - |F(s)| > 0$, donc il existe $\Delta = \Delta(\sigma_0) \in \mathbb{Z}$ tel que
$$ \bigl| \arg g(\sigma_0+i \tau) - 2 \pi \Delta \bigr| < \pi $$
pour $\tau \geq \tau_0$.
Puis
$$ \pi \int_0^{T_n} S(\tau) d\tau \leq \sigma_0 K (\sigma_0+T_n) +
\int_0^{\tau_0} \arg g(\sigma_0+i \tau)d\tau + (T_n-\tau_0)(2\Delta + 1)\pi -{I(\sigma_0)} $$
au lieu de (\ref{csuperior}). La contribution correspondante de la moyenne de $S(T)$, obtenue en divisant
cette derni\`ere in\'egalit\'e par $\pi T_n$ et en faisant $n \to +\infty$, est $ \leq \sigma_0 K + 2 \Delta + 1 $,
donc finie.
La quantit\'e analogue de $B_a$ est obtenue en calculant le bilan entre le nombre s z\'eros et les p\^oles de $f(s)$ et $h(s)$
avec $\sigma < \sigma_0$, donc qu'on note par $B_a(\sigma_0)$.
L'argument de d\'emonstration du th\'eor\`eme \ref{teoalin} nous donne finalement
\begin{equation} \label{prch}
N(\sigma_0,T)-N_0'(T) \leq B_a(\sigma_0) + \sigma_0 K + 2 \Delta + 1,
\end{equation}
ce qui montre le r\'esultat. \lqqd\\

Le point faible du th\'eor\`eme \ref{teoalindom} est l'innefectivit\'e de la borne \ref{prch},
qui d\'epend de la quantit\'e $\Delta$, difficile \`a calculer; c'est la raison
pour laquelle on ne garde
pas la borne obtenue dans l'\'enonc\'e du th\'eor\`eme.\\

On peut maintenant compl\'eter l'information du th\'eor\`eme \ref{pabierto} et la proposition \ref{hysgg} dans le cas
$0 < y < 1$.

\begin{prop} Soit $0 < y < 1$, $H(y;s)$ comme dans le th\'eor\`eme \ref{pabierto}. Alors,
pour tout $\sigma_0> 1/2$, tous les z\'eros de $H(y;s)$ dans la bande
$|\sigma-1/2| < \sigma_0-1/2$, sauf un nombre fini, sont sur la droite $\sigma=1/2$ et ils sont simples.
\end{prop}

{\bf D\'emonstration.} Il suffit de v\'erifier la condition {\it(\hyperlink{princcc1pp}{i''})} du th\`eor\`eme \ref{teoalindom}
pour $\sigma=\sigma_0$ suffisament grand. On pose $ F(y;s) = \frac{h(y;1-s)}{h(y;s)} $.
Par (\ref{cocienteyyy}), pour $\sigma$ suffisament grand, $F(y;s) = O\bigl( y^{-2\sigma} |s|^{-1/2} \bigr)$.
Donc, pour $\sigma_0$ assez grand, $\bigl| F(s) \bigr| = O_{\sigma_0} \bigl( |\tau|^{-1/2} \bigr)$ pour $ \sigma = \sigma_0$, $|\tau| \geq \tau_0$.

D'un autre c\^ot\'e, pour $\sigma_0$ fix\'e, la puissance $y^{-2\sigma}$ est born\'ee pour $a \leq \sigma \leq \sigma_0$,
ce qui suffit pour v\'erifier {\it(\hyperlink{princcc2pp}{ii''})} avec (\ref{cocienteyyy}) et les arguments habituels. Le r\'esultat se suit directement. \lqqd\\

Le dernier r\'esultat n'est pas nouveau, il a \'et\'e obtenu par H. Ki,
il correspond \`a un cas particulier de la proposition \ref{tambieneste}. Notre
m\'ethode de d\'emonstration
diff\`ere de celle de H. Ki dans le sens qu'on n'a pas d'analyse de deriv\'ee \`a faire, pour un r\'esultat sur la simplicit\'e des z\'eros etudi\'es.
La g\'en\'eralit\'e du th\'eor\`eme \ref{teoalindom} nous permettra de l'appliquer
pour \'etudier la r\'epartition
des z\'eros des approximations de la fonction z\^eta d'Epstein, au \S \ref{nonefectepstein}.

\section{Approximations de la fonction z\^eta d'Epstein}\label{nonefectepstein}

Soit $Q(u,v) = a u^2 + b uv + cv^2$ avec $a>0$, $b, c \in \re$ une forme quadratique d\'efinie positive, c'est-\`a-dire avec $\Delta = 4ac-b^2 > 0$.
La fonction z\^eta d'Epstein associ\'ee \`a $Q(u,v)$ est donn\'ee par
$$ Z_Q(s) = \summ_{(m,n) \neq (0,0)} Q(m,n)^{-s}, $$
pour $\sigma > 1$. Cette fonction peut \^etre
prolong\'ee analytiquement dans le plan complexe, sauf en $s=1$,
o\`u elle a un p\^ole simple, et
elle satisfait \`a l'\'equation fonctionnelle
$$
\biggl( \frac{\sqrt{\Delta}}{2\pi} \biggr)^s \Gamma(s) Z_Q(s) = \biggl( \frac{\sqrt{\Delta}}{2\pi} \biggr)^{1-s} \Gamma(1-s) Z_Q(1-s).
$$
Cette fonction est essentiellement une s\'erie d'Eisenstein; on change notre point de vue
pour des raisons historiques. Lorsque
les coefficients $a, b, c$ sont des entiers, qu'on suppose avec $(a,b,c)=1$, l'\'etude de cette fonction est li\'ee \`a l'\'etude
du corps quadratique imaginaire $\mathbb{Q}\bigl(\sqrt{-d}\bigr)$,
o\`u $d=\Delta/4$ si $\Delta/4$ est entier, $d=\Delta$ sinon.
Soit $h(-d)$ le nombre de classes d'id\'eaux entiers de $\mathbb{Q}\bigl(\sqrt{-d}\bigr)$
(ou des classes d'\'equivalence de formes quadratiques de discriminant $-\Delta$).
Si $h(-d)=1$, la fonction z\^eta d'Epstein est, \`a une constante pr\`es, la fonction z\^eta de Dedekind
du corps quadratique $\mathbb{Q}\bigl(\sqrt{-d}\bigr)$, elle a un produit eul\'erien,
et contient la fonction z\^eta de Riemann comme facteur
(on reprend l'exemple de la s\'erie d'Eisenstein $E(i,s)$, le cas o\`u $Q(u,v)=u^2+v^2$).
Dans ce cas, on attend que tous les z\'eros de $Z_Q(s)$ soient sur la droite critique $\sigma=1/2$,
une hypoth\`ese de Riemann pour $Z_Q(s)$. Par contre, celui-ci est le cas moins fr\'equent,
puisque les seuls entiers positifs $d$ pour lesquels $h(-d)=1 $ sont les nombres de Heegner $d=3,4,7,8,11,19,43,67,163$.
Si $h(-d)>1$, la fonction z\^eta d'Epstein n'a pas de produit eul\'erien, et elle
a une infinit\'e de z\'eros en dehors de la droite $\sigma=1/2$. Plus pr\'ecisement,
il a \'et\'e montr\'e par Davenport et Heilbronn dans \cite{davenheil1} et \cite{davenheil2},
que, lorsque $h(-d)>1$, la fonction $Z_Q(s)$ a une infinit\'e de z\'eros dans le demi-plan $\sigma >1$,
arbitrairement proches de la droite $\sigma=1$, et par Voronin dans \cite{voronin}
que dans une bande $\sigma_1 \leq \sigma \leq \sigma_2$ avec $1/2 < \sigma_1 < \sigma_2 < 1$,
cette fonction a au moins $c(\sigma_1,\sigma_2)T$ z\'eros avec $0 <\tau < T$ pour $T$ suffisament grand, o\`u $c(\sigma_1,\sigma_2)>0$.
Bien s\^ur, cette fonction a une infinit\'e de z\'eros sur la droite $\sigma=1/2$,
ce qui a \'et\'e montr\'e par Potter et Titchmarsh \cite{potter}, et
le nombre des z\'eros avec $ 0 < \tau < T$ est asymptotiquement  $ c T \log T$
(une formule explicite bien s\^ur, analogue \`a celle pour la fonction z\^eta de Riemann).

On reprend la formule de Chowla-Selberg
$$
\begin{array}{rcl}
Z_Q(s) &=& 2 \zeta(2s) a^{-s} + \Frac{2^{2s} a^{s-1} \sqrt{\pi}}{\Gamma(s) \Delta^{s-1/2}} \zeta(2s-1) \Gamma\Bigl( s-\frac{1}{2} \Bigr) \\
& & + \Frac{4 \pi^s 2^{s-1/2}}{\sqrt{a} \Gamma(s) \Delta^{s/2-1/4}}
\summ_{n=1}^\infty n^{s-1/2} \summ_{d|n} d^{1-2s} \cdot \cos
\Frac{n\pi b}{a} \int_{-\infty}^{+\infty} e^{-(\pi n
\sqrt{\Delta}/a) \cosh \tau} e^{(s-1/2)\tau} d\tau,
\end{array}
$$
et on obtient des approximations de la fonction z\^eta d'Epstein en prenant les sommes partielles
des termes dans cette formule. Pour $n \geq 1$ entier, on d\'efinit $Z_{Q,n}(s)$ par
$$
\begin{array}{rcl}
Z_{Q,n}(s) &=& 2 \zeta(2s) a^{-s} + \Frac{2^{2s} a^{s-1} \sqrt{\pi}}{\Gamma(s) \Delta^{s-1/2}} \zeta(2s-1) \Gamma\Bigl( s-\frac{1}{2} \Bigr) \\
& & + \Frac{4 \pi^s 2^{s-1/2}}{\sqrt{a} \Gamma(s) \Delta^{s/2-1/4}}
\summ_{k=1}^n k^{s-1/2} \summ_{d|k} d^{1-2s} \cdot \cos
\Frac{k\pi b}{a} \int_{-\infty}^{+\infty} e^{-(\pi k \sqrt{\Delta}/a) \cosh \tau} e^{(s-1/2)\tau} d\tau.
\end{array}
$$
Ces fonctions ont \'et\'e pr\'ec\'edemment \'etudi\'ees par D. Hehjal dans \cite[\S 5]{hejhal},
et reprises par H. Ki dans \cite{hki2}, qui am\'eliore ses r\'esultats. On va red\'emontrer les r\'esultats de H. Ki avec nos techniques.
On d\'efinit pour cela une classe de fonctions qui g\'en\'eralise $Z_{Q,n}(s)$. Soient
$\lambda > 0$, $A_j>0$, $\lambda_j>0$ et $b_j \in \re$ pour $j=1,\dots,n$.
On note par $G(s)$ la fonction
\begin{equation} \label{clasegdes}
\begin{array}{rcl}
G(s) &=& \zeta(2s) \lambda^{s} + \sqrt{\pi} \Frac{\Gamma(s-\frac{1}{2})}{\Gamma(s)} \zeta(2s-1) \lambda^{1-s} \\
&& + \Frac{\pi^s}{\Gamma(s)}
\summ_{k=1}^n b_k (\lambda_k^{s-1/2}+\lambda_k^{1/2-s}) \int_{-\infty}^{+\infty} e^{-2 A_k \cosh \tau} e^{(s-1/2)\tau} d\tau.
\end{array}
\end{equation}
On voit facilement que $\frac{1}{2} \bigl( \sqrt{\Delta}/2 \bigr)^s Z_{Q,n}(s)$ est un cas particulier de $G(s)$, o\`u $\lambda = \sqrt{\Delta}/(2a)$ et
sous les conditions aditionnelles $0 < \lambda_j \leq 1$ et $\lambda_j A_j \geq \pi \lambda$ pour $j=1,\dots,n$.

La parit\'e de la fonction de Bessel et l'\'equation fonctionnelle de la fonction z\^eta de Riemann entra\^inent
l'\'equation fonctionnelle pour $G(s)$
$$ \pi^{-s} \Gamma(s) G(s) = \pi^{-1+s} \Gamma(1-s) G(1-s), $$
qui \`a son tour, entra\^ine celle de la fonction $Z_{Q,n}(s)$
$$ \biggl( \frac{\sqrt{\Delta}}{2\pi} \biggr)^s \Gamma(s) Z_{Q,n}(s) =
\biggl( \frac{\sqrt{\Delta}}{2\pi} \biggr)^{1-s} \Gamma(1-s)
Z_{Q,n}(1-s). $$

En fait, l'\'equation fonctionnelle de $G(s)$ est une cons\'equence directe de l'egalit\'e
$$ \pi^{-s} \Gamma(s) G(s) = w(s) + w(1-s), $$
o\`u
\begin{equation} \label{nuevaw}
w(s) = \lambda^s \pi^{-s} \Gamma(s) \zeta(2s) + \sum_{k=1}^n b_k \lambda_k^{s-\frac{1}{2}} 2K_{s-\frac{1}{2}}(2A_k).
\end{equation}

L'estimation suivante est analogue \`a celle dans \cite[Proposition 3.5]{hki2}.

\begin{prop} \label{propkbessel} Il existe $R> 0$ tel que
$$\frac{K_{s-1/2}(2A)}{\Gamma(s)} = O \Bigl( \frac{A^{-\sigma}}{|s|^{1/2}} \Bigr)$$
pour $|s|> R$, $\sigma \geq \frac{1}{2}$.
\end{prop}

{\bf D\'emonstration.} La formule dans \cite[p. 309]{polya}
\begin{equation} \label{vvpolya}
\begin{array}{rcl}
K_s(2A) &=& A^{-s} \Gamma(s) \Bigl( 1 + \summ_{n=1}^\infty \Frac{A^{2n}}{n!(1-s)(2-s)\cdots(n-s)} \Bigr) \\
& & + \Gamma(-s) \Bigl( 1 + \summ_{n=1}^\infty \Frac{A^{2n}}{n!(1+s)(2+s)\cdots(n+s)} \Bigr)
\end{array}
\end{equation}
permet de lier les comportements asymptotiques de la fonction de Bessel et de la fonction gamma.
Il a \'et\'e montr\'e par P\'olya dans \cite[pp. 309--310]{polya} que la fonction enti\`ere $ \psi_A(s) = A^s {K_s(2A)}/{\Gamma(s)} - 1 $
satisfait $\ds\lim_{|s| \to \infty} \psi_a(s) = 0$ dans tout demi-plan $\sigma \geq \eps$, $\eps >0$ fix\'e.
Le m\^eme calcul montre que $\psi_a(s)$ est born\'ee dans le demi-plan $\sigma \geq 0$. Or
$$ \Frac{K_{s-1/2}(2A)}{\Gamma(s)} = \bigl(1+\psi_a(s-1/2)\bigr) A^{1/2-s} \Frac{\Gamma(s-1/2)}{\Gamma(s)}. $$
On applique la formule de Stirling pour obtenir le r\'esultat. \lqqd\\

Ensuite on \'enonce et d\'emontre une variante du r\'esultat de H. Ki \cite[Theorem 1]{hki2}.

\begin{prop} \label{hkicopia} Soient $\lambda > 0$, $A_j>0$, $\lambda_j>0$ et $b_j \in \re$ pour $j=1,\dots,n$, et $G(s)$ la fonction
d\'efinie par (\ref{clasegdes}). Alors, pour tout $\sigma_0 > 1/2$, tous les z\'eros de $G(s)$ dans la bande
$|\sigma-a| < \sigma_0-a$, sauf un nombre fini, sont simples et sur la droite
$\sigma = \frac{1}{2}$.

Sous les conditions additionnelles $\lambda \geq 1$, $0 < \lambda_j \leq 1$ et $\lambda_j A_j \geq \pi \lambda$ pour $j=1,\dots,n$,
tous les z\'eros non triviaux de $G(s)$, sauf un nombre fini, sont simples et sur la droite $\sigma = \frac{1}{2}$.
\end{prop}

\begin{coro}[H. Ki {\cite[Corollary 1]{hki2}}] \label{hkidos} Pour les fonctions $Z_{Q,n}(s)$, on a le r\'esultat de la proposition
\ref{hkicopia}. Si $\sqrt{\Delta}/(2a) \geq 1$, tous les z\'eros non triviaux de $Z_{Q,n}(s)$, sauf un nombre fini, sont simples
et sur la droite $\sigma = \frac{1}{2}$.
\end{coro}

{\bf D\'emonstration (de la proposition).} Les z\'eros triviaux des fonctions $G(s)$ et $Z_{Q,n}(s)$ proviennent de la fonction gamma, comme pour la
fonction z\^eta de Riemann. On va appliquer les th\'eor\`emes \ref{teoalin} et \ref{teoalindom}.
Il est clair que les seuls p\^oles de $w(s)$ sont $s=0$, $s=1/2$. On pose
$$
\mu= \max_{1 \leq k \leq n}\lambda_k^{-1} A_k^{-1} \pi \lambda^{-1}, \quad
\nu= \max_{1 \leq k \leq n} \lambda_k A_k^{-1} \pi \lambda^{-1}.
$$On factorise $w(s)$ par
$$ w(s) = \Bigl( \frac{\lambda}{\pi} \Bigr)^s \Gamma(s) \Biggl( \zeta(2s) + \Bigl( \frac{\lambda}{\pi} \Bigr)^{-s} \sum_{k=1}^n 2b_k \lambda_k^{s-\frac{1}{2}} \frac{K_{s-\frac{1}{2}}(2A_k)}{\Gamma(s)}\Biggr). $$
D'apr\`es la proposition \ref{propkbessel}, pour $|s|$ suffisament grand, $\sigma \geq 1/2$, $1 \leq k \leq n$
$$ \Bigl( \frac{\lambda}{\pi} \Bigr)^{-s}
\lambda_k^{s-\frac{1}{2}} \frac{K_{s-\frac{1}{2}}(2A_k)}{\Gamma(s)} =
O \Bigl( \lambda^{-\sigma} \pi^{\sigma} \lambda_k^{\sigma} \frac{A_k^{-\sigma}}{|s|^{1/2}} \Bigr) =
O \Bigl( \frac{\nu^{\sigma}}{|s|^{1/2}} \Bigr),
$$
d'o\`u
$$ w(s) = \Bigl( \frac{\lambda}{\pi} \Bigr)^s \Gamma(s) \Bigl( \zeta(2s) + O\bigl( \nu^{\sigma} |s|^{-1/2} \bigr) \Bigr)
$$
et $\zeta(2s)^{-1} = O\bigl( |s|^{1/4} \bigr)$ pour $\sigma \geq 1/2$ , puis
\begin{equation} \label{wsolo}
w(s)= \Bigl( \frac{\lambda}{\pi} \Bigr)^s \Gamma(s) \zeta(2s) \Bigl( 1 + O\bigl( \nu^{\sigma} |s|^{-1/4} \bigr) \Bigr),
\end{equation}
pour $\sigma \geq 1/2$ et $|s|$ suffisament grand.

Si l'on \'ecrit $F(s) = w(1-s)/w(s)$,
$$
\begin{array}{rcl}
F(s) &=&
\displaystyle
\frac{\lambda^{1-s}\pi^{\frac{1}{2}-s} \Gamma(s-\frac{1}{2})\zeta(2s-1) + \displaystyle\sum_{k=1}^n b_k \lambda_k^{\frac{1}{2}-s} 2{K_{s-\frac{1}{2}}(2A_k)}} {\Bigl( \frac{\lambda}{\pi} \Bigr)^s \Gamma(s) \zeta(2s)\Bigl( 1 + O\bigl( \nu^\sigma \sigma^{-1/4} \bigr) \Bigr)}\\
    &=&
\displaystyle
\frac{\ds \lambda^{1-2s} \sqrt{\pi} \frac{\Gamma(s-\frac{1}{2})}{\Gamma(s)} \zeta(2s-1) + \displaystyle \sum_{k=1}^n b_k
\lambda_k^{\frac{1}{2}-s} \lambda^{-s} \pi^s 2\frac{K_{s-\frac{1}{2}}(2A_k)}{\Gamma(s)}} {\zeta(2s) \Bigl(1 + O\bigl( \nu^\sigma |s|^{-1/4} \bigr)\Bigr)}.
\end{array}
$$
Or, encore par la proposition \ref{propkbessel}, on a pour $|s|$ assez grand, $\sigma \geq 1/2$, $1 \leq k \leq n$
$$ \lambda_k^{\frac{1}{2}-s} \lambda^{-s} \pi^s \frac{K_{s-\frac{1}{2}}(2A_k)}{\Gamma(s)} =
O \Bigl( \lambda^{-\sigma} \pi^{\sigma} \lambda_k^{-\sigma} \frac{A_k^{-\sigma}}{|s|^{1/2}} \Bigr) =
O \Bigl( \frac{\mu^{\sigma}}{|s|^{1/2}} \Bigr).
$$
Puis, pour $\sigma\geq 1/2$, $|s|$ suffisament grand
\begin{equation} \label{woverw}
\begin{array}{rcl}
F(s) &=&
\displaystyle
\frac{O\bigl(\lambda^{-2\sigma} |s|^{-1/2}\bigr) \zeta(2s-1) + O\bigl(\mu^\sigma |s|^{-1/2}\bigr)}{\zeta(2s) \Bigl(1 + O\bigl( \nu^\sigma |s|^{-1/4} \bigr)\Bigr)}.
\end{array}
\end{equation}

\begin{enumerate}
\alfabetico
\item \label{bostezo} Soit $\sigma_0> 1/2$ tel que $w(s)\neq 0$ et $F(s) \neq \pm 1$ sur $\sigma = \sigma_0$,
et $\zeta(s) \asymp 1$ pour $\sigma \geq \sigma_0$ ($\sigma_0> 2$). Dans
(\ref{woverw}), on fait $|\tau| \to +\infty$ pour $s=\sigma_0+i\tau$, et on obtient $F(s) = O_{\sigma_0}\bigl( |\tau|^{-1/2}\bigr)$, d'o\`u
la condition {\it(\hyperlink{princcc1pp}{i''})} du th\'eor\`eme \ref{teoalindom}. Pour $1/2 \leq \sigma \leq \sigma_0$,
les puissances $C^\sigma$ avec $C$ constante, agissent comme des constantes. Puis, les estimations
$\zeta(2s-1) = O\bigl(|\tau|^A\bigr)$ et $\zeta(2s)^{-1} = O\bigr( \log |\tau| \bigr)$ pour $\sigma \geq 1/2$
et $|\tau|$ assez grand, entra\^inent que $F(s)= O_{\sigma_0}\bigl( |\tau|^A \log |\tau| \bigr)$
pour $1/2 \leq \sigma \leq \sigma_0$, $|\tau| \geq \tau_0$, donc la condition {\it(\hyperlink{princcc2pp}{ii''})}
du th\'eor\`eme \ref{teoalindom}. La premi\`ere partie de la proposition d\'ecoule de ce r\'esultat.
\item On suppose maintenant les conditions $\lambda \geq 1$, $0 < \lambda_j \leq 1$ et $\lambda_j A_j \geq \pi \lambda$ pour $j=1,\dots,n$. Alors,
 $\mu \leq 1$ et $\nu \leq 1$. La fonction $w(s)$ n'a qu'un nombre fini de z\'eros avec $\sigma \geq 1/2$,
 puisque dans (\ref{wsolo}) $w(s) = \lambda^s \zeta^*(2s) \Bigl( 1 + O \bigl( |s|^{-1/4} \bigr) \Bigr)$,
 puis $w(s) \neq 0 $ pour $|s|$ assez grand, $\sigma \geq 1/2$. L'\'equation (\ref{woverw}) devient
 $$
F(s) = \displaystyle
\frac{O\bigl(|s|^{-1/2}\bigr) \zeta(2s-1) + O\bigl(|s|^{-1/2}\bigr)}{\zeta(2s) \Bigl(1 + O\bigl(|s|^{-1/4} \bigr)\Bigr)}.
$$
Lorsque $\sigma \to +\infty$, on a $F(s) = O(\sigma^{-1/2})$ uniform\'ement en $\tau$, donc la condition {\it(\ref{princcc1})}
du th\'eor\`eme \ref{teoalin}.
Lorsque $\tau \to +\infty$, on a $F(s) = O\bigl( |\tau|^A \log |\tau| \bigr)$ pour $\sigma \geq 1/2$, $|\tau|$
assez grand, donc la condition {\it(\ref{princcc2})} du th\'eor\`eme \ref{teoalin} (on fait on a dej\`a montr\'e la condition souhait\'ee,
mais on a mieux qu'en \itemref{bostezo}). Le th\'eor\`eme \ref{teoalin} nous donne
la seconde partie de la proposition, ainsi que le corollaire \ref{hkidos}. \hfill{$\blacksquare$}
\end{enumerate}

Dans le cas $Q(u,v) = u^2+v^2$, on a $\lambda = \sqrt{\Delta}/(2a) =1$, d'o\`u l'inter\^et de ce r\'esultat.
L'avantage du th\'eor\`eme \ref{teoalin} sur le th\'eor\`eme \ref{teoalindom} est
la majoration effective du nombre de z\'eros de $G(s)$ en dehors de la droite critique,
le cas o\`u on pourrait estimer le nombre de z\'eros de la fonction $w(s)$ d\'efinie par
(\ref{nuevaw}) dans le demi-plan $\sigma \geq 1/2$.

Dans le cas
o\`u $0 < \lambda < 1$, H. Ki calcule aussi le nombre de z\'eros des fonctions $G(s)$
hors une bande $|\sigma-a| < \sigma_0-a$ fix\'ee et $0 < \tau < T$, pour montrer qu'il est $O(T)$,
en utilisant une approximation des fonctions de Bessel \cite[Proposition 3.5]{hki2}, pour approcher $G(s)$ avec une fonction semblable \`a une
exponentielle-polyn\^ome dans le demi-plan $\sigma \geq \sigma_0$, comme dans la
proposition \ref{hysgg}.\\

H. Ki a aussi \'etudi\'e des fonctions z\^eta d'Epstein modifi\'ees, dans \cite{hki3}.
On cite le texte original.
Soient $L_1(s), \dots, L_n(s)$ des s\'eries de Dirichlet qui peuvent \^etre continu\'ees analytiquement
au plan complexe, sauf en un nombre fini de points. On suppose que
\begin{enumerate}
\romano
\item \label{ele1} $y > 0$;
\item \label{ele2} $\alpha$ et $\beta$ sont des polyn\^omes tels que $\deg (\alpha) \geq \deg(\beta)+1$ et $\beta(s) = \beta(1-s)$;
\item \label{ele3} $a_k$ est r\'eel pour $k=1,\dots,n$;
\item \label{ele4} $\pi^{-s} \Gamma(s) L_k(s) = \pi^{-1+s} \Gamma(1-s) L_k(1-s)$, pour $k=1,2,\dots,n$, $s \in \ce \setminus \{0,1\}$;
\item \label{ele5} il existe $\eps > 0$ tel que $L_k(s) = O(|s|^{1-\eps})$, pour $k=1,2,\dots,n$, $\sigma \geq \frac{1}{2}$.
\end{enumerate}

On d\'efinit $Z(s;y,\alpha,\beta,L_1,\dots,L_n)$ par
$$ Z(s;y,\alpha,\beta,L_1,\dots,L_n) = \alpha(s) \zeta(2s) + \alpha(1-s) \sqrt{\pi} \frac{\Gamma(s-1/2)}{\Gamma(s)} \zeta(2s-1) y^{1-2s} + y^{-s} \beta(s) \summ_{k=1}^n a_k L_k(s). $$

En fait, si
$$ w(s) = \alpha(s) y^s \zeta^*(2s) + \frac{\beta(s)}{2} \summ_{k=1}^n a_k \pi^{-s} \Gamma(s) L_k(s), $$
on a
$$ \pi^{-s} \Gamma(s) y^s Z(s;y,\alpha,\beta,L_1,\dots,L_n) = w(s) + w(1-s). $$
Puis, on a l'\'equation fonctionnelle
$$ \Bigl( \Frac{y}{\pi} \Bigr)^s \Gamma(s) Z(s;y,\alpha,\beta,L_1,\dots,L_n) =
\Bigl( \Frac{y}{\pi} \Bigr)^{1-s} \Gamma(1-s) Z(1-s;y,\alpha,\beta,L_1,\dots,L_n)$$

\begin{prop}[H. Ki {\cite[Theorem 1.1]{hki3}}] \label{tambieneste} Pour tout $\sigma_0 > 1/2$, tous les z\'eros non-triviaux
de la fonction $Z(s;y,\alpha,\beta,L_1,\dots,L_n)$ dans la bande $|\sigma-1/2| < \sigma_0-1/2$, sauf un nombre fini,
sont simples et sur la droite $\sigma=1/2$.
\end{prop}

\begin{coro}[H. Ki {\cite[Corollary 1.2]{hki3}}]
Si $y \geq 1$, tous les z\'eros de $Z(s;y,\alpha,\beta,L_1,\dots,L_n)$, sauf
un nombre fini, sont simples et sur la droite $\sigma = \frac{1}{2}$.
\end{coro}

{\bf D\'emonstration.}
Par (\ref{ele2}) et (\ref{ele5}), il existe $\eps > 0$ tel que
$$
\begin{array}{rcl}
w(s) &=& \alpha(s) y^s \pi^{-s} \Gamma(s) \Bigl( \zeta(2s) + y^{-s} \ds \frac{\beta(s)}{2\alpha(s)} \summ_{k=1}^n a_k L_k(s) \Bigr)\\
    &=& \alpha(s) y^s \pi^{-s} \Gamma(s) \Bigl( \zeta(2s) + O\bigl(y^{-\sigma}|s|^{-\eps}\bigr) \Bigr)
\end{array}
$$
et $\zeta(2s)=O\bigl( |s|^{\eps/2} \bigr)$ pour $\sigma \geq 1/2$, puis

\begin{equation} \label{rererew}
w(s) = \alpha(s) y^s \pi^{-s} \Gamma(s) \zeta(2s) \Bigl( 1 + O\bigl(y^{-\sigma}|s|^{-\eps/2}\bigr) \Bigr)
\end{equation}
pour $\sigma \geq \frac{1}{2}$, $|s|$ assez grand.

D'un autre c\^ot\'e, $F(s) = w(1-s)/w(s)$ satisfait
$$
\begin{array}{rcl}
F(s) &=&
\Frac{\alpha(1-s) y^{1-s} \pi^{\frac{1}{2}-s} \Gamma(s-\frac{1}{2}) \zeta(2s-1) + \frac{\beta(s)}{2} \summ_{k=1}^n a_k \pi^{-s} \Gamma(s) L_k(s)}
{\alpha(s) y^s \pi^{-s} \Gamma(s) \Bigl( \zeta(2s) + O \bigl( y^{-\sigma} |s|^{-\eps} \bigr) \Bigr)} \\
    &=&
\ds\frac{\alpha(1-s)}{\alpha(s)} y^{1-2s} \pi^{1/2}
\frac{ \zeta(2s-1) \Frac{\Gamma(s-\frac{1}{2})}{\Gamma(s)} + \summ_{k=1}^n \frac{a_k}{2\pi^{1/2}} \frac{\beta(s)}{\alpha(s)} L_k(s)}{\zeta(2s)\Bigl( 1+O\bigl( y^{-2\sigma}|s|^{-\eps/2} \bigr)\Bigr) }\\
    &=&
\ds\frac{\alpha(1-s)}{\alpha(s)} y^{1-2s} \pi^{1/2} \Frac{ \zeta(2s-1) O\bigl(|s|^{-\frac{1}{2}} \bigr) + O \bigl( |s|^{-\eps} \bigr)}{\zeta(2s)\Bigl( 1+O\bigl( y^{-\sigma}|s|^{-\eps/2} \bigr)\Bigr)}.
\end{array}
$$
Donc, pour $\sigma \geq 1/2$, $|s|$ assez grand
\begin{equation}  \label{rerererew}
F(s) = (-1)^{\deg \alpha} \Bigl( 1 + O\bigl(|s|^{-1}\bigr)\Bigr) y^{1-2s} \pi^{1/2} \Frac{ \zeta(2s-1) O\bigl(|s|^{-\frac{1}{2}} \bigr) + O \bigl( |s|^{-\eps} \bigr)}{\zeta(2s)\Bigl( 1+O\bigl( y^{-\sigma}|s|^{-\eps/2} \bigr)\Bigr)}.
\end{equation}

\begin{enumerate}
\alfabetico
\item \label{pftl} Soit $\sigma_0 > 1/2$ suffisament grand, tel que $w(s) \neq 0$ et $F(s) \neq \pm 1$ sur $\sigma=\sigma_0$, et $\zeta(s) \asymp 1$ pour $\sigma \geq \sigma_0$.
Lorsque $|\tau| \to +\infty$, $s=\sigma_0+i \tau$, on obtient de (\ref{rerererew}), $F(s) = O_{\sigma_0}\bigl( |\tau|^{-\eps} \bigr)$,
d'o\`u la condition {\it(\hyperlink{princcc1pp}{i''})} du th\'eor\`eme \ref{teoalindom}. Pour $1/2 \leq \sigma \leq \sigma_0$, la puissance $y^{-2\sigma}$ est born\'ee;
ceci avec $\zeta(2s-1) = O\bigl( |\tau|^A \bigr)$ et $\zeta(2s)^{-1} = O\bigl( \log |\tau| \bigr)$ nous donne
$ F(s) = O_{\sigma_0} \bigl( |\tau|^A \log |\tau| \bigr)$ pour $1/2 \leq \sigma \leq \sigma_0$,
donc la condition {\it(\hyperlink{princcc2pp}{ii''})} du th\'eor\`eme \ref{teoalindom}. Cela montre la proposition.
\item On suppose maintenant que $y \geq 1$. L'\'equation (\ref{rererew}) se simplifie et on obtient
$$ w(s) = \alpha(s) y^s \zeta^*(2s) \Bigl( 1 + O\bigl(|s|^{-\eps/2}\bigr) \Bigr) $$
pour $\sigma \geq 1/2$, $|s|$ suffisament grand. Le nombre de z\'eros de $\alpha(s)$ \'etant fini, cela entra\^ine
que $w(s)$ n'a qu'un nombre fini de z\'eros dans le demi-plan $\sigma \geq 1/2$.
Pour $\sigma$ suffisament grand, (\ref{rerererew}) devient $F(s) = O(\sigma^{-\eps})$,
d'o\`u la condition {\it(\ref{princcc1})} du th\'eor\`eme \ref{teoalin}. Dans \itemref{pftl} on a montr\'e
une condition sufissante pour la condition {\it(\ref{princcc2})} du th\'eor\`eme \ref{teoalin},
mais on peut montrer que $ F(s) = O \bigl( |\tau|^A \log |\tau| \bigr)$ pour $\sigma \geq 1/2$ lorsque $y \geq 1$.
L'application du th\'eor\`eme \ref{teoalin} nous donne le corollaire. \lqqd
\end{enumerate}

\section{\'Etude en l'absence de sym\'etrie r\'eelle} \label{sanssymetrie}

On veut g\'en\'eraliser les r\'esultats obtenus jusqu'\`a pr\'esent pour des fonctions
qui peuvent prendre des valeurs complexes sur la droite r\'eelle.
Pour traiter ce cas, \'etant donn\'ee une fonction $h(s)$ m\'eromorphe, d\'efinie sur un domaine $D \subset \ce$, on d\'efinit
$$ \ovl{h}(s) = \ovl{h(\ovl{s})} $$
pour $s \in \ovl{D} = \{\ovl{z}, z \in D\}$.
La fonction $\ovl{h}(s)$ est aussi m\'eromorphe. Si $h(s)$ est r\'eelle sur la droite r\'eelle, alors $\ovl{h}(s)= h(s)$
pour $s \in D \cap \ovl{D}$ par le principe de r\'eflexion, pourvu que $D \cap \mathbb{R} \neq \emptyset$.

On veut maintenant \'etudier les fonctions
$$ f(s) = f^{\pm}(s) = h(s) \pm \ovl{h}(2a-s), $$
o\`u $a \in \re$ est fix\'e. Ces fonctions r\'epr\'esentent les parties r\'eelle ($f^+(s)$) et imaginaire ($f^-(s)$)
de la fonction $h(s)$, sur la droite critique $\sigma=a$. Au lieu des deux sym\'etries de $f(s)$ on ne trouve qu'une seule, \`a savoir
$$ \ovl{f}(2a-s) = \pm f(s). $$

En analogie aux notations au \S \ref{notaciones}, on note par
$$ \hat{N}(T) = \# \bigl\{s \in \ce \tq f(s)= 0, |\tau| < T \bigr\} $$
le nombre de z\'eros de $f(s)$ avec $|\tau| < T$,
$$ \hat{N}_0(T) = \# \bigl\{s \in \ce \tq f(s)= 0, s= a+i\tau, |\tau| < T \bigr\}$$
le nombre des m\^emes z\'eros, sur la droite critique $\sigma=a$, et $\hat{N}'_0(T)$
le nombre de ces derniers z\'eros, sans compter les multiplicit\'es.

D'abord, on estime le nombre des z\'eros sur la droite critique, sans compter les multiplicit\'es. La g\'en\'eralisation
du lemme \ref{rectacritica} est la suivante.

\begin{lema} \label{rectacriticac} Soit $a \in \re$ fix\'e, $h(s)$ une fonction m\'eromorphe dans le plan complexe,
r\'eelle sur l'axe r\'eel, sans z\'eros ni p\^oles sur la droite critique $\sigma = a$.
On consid\`ere $\arg h(a+i\tau)$, une variation continue de l'argument de $h(s)$ sur la droite critique $\sigma=a$.
Soit $f(s) = h(s) \pm \ovl{h}(2a-s)$. Le nombre de z\'eros de $f(s)$ avec $|\tau| < T$, sans compter
les multiplicit\'es, sur la droite $\sigma = a$ est minor\'e par
\begin{equation} \label{mejorable}
\tilde{N}_0'(T) \geq \mfrac{1}{\pi}\arg h(a+iT)  - \mfrac{1}{\pi}\arg h(a-iT) -1.
\end{equation}
\end{lema}

{\bf D\'emonstration.} On fait la preuve pour la fonction $f^-(s)$; en prenant $i h(s)$ au lieu de $h(s)$
on obtiendra $i f^+(s)$.

On peut supposer que $\arg h(a) \in ]-\pi,0]$. Si $f^-(a) = 0$,
$$
\begin{array}{rcl}
\hat{N}_0'(T) &=& \# \bigl\{ \tau \tq |\tau|< T,  \frac{1}{\pi} \arg h(a+i\tau) \equiv 0 \mmod 1 \bigr\} \\[0.2cm]
    &=& \# \bigl\{ \tau \tq 0 \leq \tau < T, \frac{1}{\pi} \arg h(a+i\tau) \equiv 0 \mmod 1 \bigr\}\\[0.2cm]
    & &      + \# \bigl\{ \tau \tq 0 \leq \tau < T, -\frac{1}{\pi} \arg h(a-i\tau) \equiv 0 \mmod 1 \bigr\} -1.
\end{array}
$$
On obtient (\ref{mejorable}) par double application du lemme \ref{compcreciente}.
De m\^eme si $f^-(a) \neq 0$. \hfill{$\blacksquare$}\\

{\bf Remarque.} Dans le cas o\`u $h(s)$ est r\'eelle sur la droite r\'eelle, on a $\tilde{N}_0'(T) = 2 N_0'(T) + 2 u_{\pm} - 1$
et $ \arg h(a-i\tau) = -h(a+i\tau)$ (o\`u $\arg h(a) = 0$), d'o\`u l'on retrouve le lemme \ref{rectacritica}.\\

Finalement, on retrouve la g\'en\'eralisation du th\'eor\`eme \ref{teoalin}, dont on ne fournit pas la preuve.

\begin{teor} \label{teoalinc} Soit $a \in \re$, $h(s)$ une fonction
m\'eromorphe sur $\ce$, avec un
nombre fini de p\^oles, un nombre fini de z\'eros dans le demi-plan $\sigma > a$, holomorphe et non nulle sur la droite critique $\sigma = a$.
On d\'efinit
$$ f(s) = f^{\pm}(s) = h(s) \pm \ovl{h}(2a-s) $$
(en particulier $\ovl{f}(2a-s)=\pm f(s)$). On suppose que la fonction
$$ F(s) = \Frac{\ovl{h}(2a-s)}{h(s)} $$
 satisfait
\begin{enumerate}
\item[\hypertarget{princcc1x}(i)] pour chaque $\eta > 0$, il existe $\sigma_0 = \sigma_0(\eta) > a$ tel que $\bigl|F(s)\bigr| < \eta $
si $\sigma \geq \sigma_0$, $\tau \in \re$;
\item[\hypertarget{princcc2x}(ii)] pour chaque $\eps>0$ et $\sigma_0 > a$, il existent deux suites $(T_n)_n$, $(T_n^*)_n$
telles que $\ds\lim_{n \to \infty} T_n = \ds\lim_{n \to \infty} T_n^* = +\infty$ et $\bigl| F(s) \bigr| < e^{\eps|s|}$
pour $a \leq \sigma \leq \sigma_0$, et $\tau = T_n$ ou $\tau = -T_n^*$, $n \geq 1$.
\end{enumerate}
Avec les notations du th\'eor\`eme \ref{teoalin},
\begin{equation} \label{naac}
\hat{N}(T) - \hat{N}_0(T) \leq \hat{N}(T) - \hat{N}_0'(T) \leq 1 + 2 P_{f,\sigma > a} + 2 N_{h,\sigma> a} - 2 P_{h,\sigma > a}.
\end{equation}
 pour $T > 0$. En particulier, presque tous les z\'eros de $f(s)$ se trouvent sur la
droite $\sigma = a$ et ils sont simples. Le membre de gauche de (\ref{naac}) est de plus un nombre positif pair.
\end{teor}

{\bf Remarque.} Dans le cas o\`u $h(s)$ est r\'eelle sur la droite r\'eelle, on a $\hat{N}(T) = 2 N(T) + 2 n_{f,\sigma >a} + n_{f,a}$,
ce qui avec la remarque faite ci-dessus, nous permet de retrouver le th\'eor\`eme \ref{teoalin}.
Une cons\'equence imm\'ediate de ce th\'eor\`eme est la r\'epartition globale des z\'eros de $f(s)$
(comme dans le th\'eor\`eme \ref{enezero}).

\begin{teor} \label{enezeroc} Sous les conditions du th\'eor\`eme \ref{teoalinc}, si $\arg h(a+it)$ est
une d\'etermination continue de l'argument de $h(s)$ sur la droite $a+i \re$, alors
$$ \hat{N}_0'(T) = \frac{1}{\pi} \arg h(a+iT) - \frac{1}{\pi} \arg h(a-iT) + O(1). $$
Les z\'eros des fonctions $f(s) = h(s) \pm \ovl{h}(2a-s)$ sont
ultimement simples et entrelac\'es.
\end{teor}

\begin{coro} \label{coralinc} Sous les conditions du th\'eor\`eme \ref{teoalinc}, si $h(s)$ est une fonction enti\`ere, alors
$$ \hat{N}(T) - \hat{N}_0(T) \leq \hat{N}(T) - \hat{N}_0'(T) \leq 1 + 2 N_{h,\sigma> a}.$$
\end{coro}

On appliquera ces r\'esultats \`a l'\'etude des translat\'ees des fonctions $L$ compl\'et\'ees sur $\mathbb{Q}$ au lieu
des fonctions $\zeta^*(s)$ ou $\xi(s)$ au \S \ref{aquiestimation}.
Malheureusement, l'estimation
de $\hat{N}(T)-\hat{N}_0'(T)$ obtenue n'est pas optimale. En fait, si $h(s)$ n'a pas de z\'ero dans le demi-plan $\sigma \geq a$,
dans le dernier corollaire $N_{h,\sigma> a} = 0$, d'o\`u $\hat{N}(T) = \hat{N}_0(T)$, ce qui nous dit que tous les z\'eros
de $f(s)$ sont align\'es sur la droite $\sigma=a$, mais par contre on a seulement que $\hat{N}_0(T) - \hat{N}_0'(T) \leq 1$,
ce qui indique que tous les z\'eros de $f(s)$ sont simples, sauf peut-\^etre un. Dans le cas r\'eel (th\'eor\`eme \ref{teoalin}), la sym\'etrie
additionnelle r\`egle ce probl\`eme. Dans le cas g\'en\'eral, il s'applique
la remarque (\ref{phragmen}.\hyperlink{phragmen1}{1}) faite au th\'eor\`eme \ref{teoalin}, et la proposition \ref{copiado} montre la simplicit\'e de tous les z\'eros de $f(s)$.

Si l'on oublie pour l'instant la n\'ecessit\'e de montrer la simplicit\'e des z\'eros de $f(s)$, on peut affaiblir
les conditions sur $h(s)$, comme dans le \S \ref{policaso}.

\begin{coro} \label{debrujingn} Soit $a \in \re$, $h(s)$ une fonction enti\`ere, avec un nombre fini de z\'eros dans le demi-plan $\sigma \geq a$,
sans z\'eros sur la droite $\sigma = a$. Avec les notations du th\`eor\`eme \ref{teoalinc}, si la fonction $F(s)$
satisfait
\begin{enumerate}
\item[\hypertarget{princcc1p}(i')] il existent $\sigma_0 \geq a$, $C >0$ tels que $\bigl|F(s)\bigr| \leq C $
si $\sigma \geq \sigma_0$, $\tau \in \re$,
\end{enumerate}
et la condition (\hyperlink{princcc2x}{ii}) du th\'eor\`eme, alors pourvu que la fonction $f^{\pm}(s)=h(s) \pm \ovl{h}(2a-s)$ soit non nulle,
$$ \hat{N}(T) - \hat{N}_0(T) \leq 2 N_{h,\sigma> a}.$$
\end{coro}

{\bf D\'emonstration.} Par la remarque (\ref{estata}) au th\'eor\`eme \ref{teoalin}, on peut \'eliminer
les z\'eros (en nombre fini) sur la droite $\sigma=a$; la proc\'edure ne change pas la quantit\'e $\hat{N}(T)-\hat{N}_0(T)$.
On introduit la famille de fonctions
$$ H(y;s) = y^s h(s) \pm y^{2a-s} \ovl{h}(2a-s) $$
pour le param\`etre $y> 1$. Il est \'evident que $f(s) = \ds\lim_{y \to 1^+} H(y;s)$ uniform\'ement
dans les compactes de $\ce$. La fonction
$$ F(y;s) = y^{2a-2s} \frac{\ovl{h}(2a-s)}{h(s)} = y^{2a-2s} F(s) $$
satisfait maintenant la condition {\it(\hyperlink{princcc1x}{i})} du th\'eor\`eme \ref{teoalinc}. En plus, les z\'eros de
$y^s h(s)$ sont les m\^emes que ceux de $h(s)$.
Donc
$$ \hat{N}(y;T) - \hat{N}_0(y;T) \leq \hat{N}(y;T) - \hat{N}_0'(y;T) \leq 1 + 2 N_{h,\sigma> a}.$$
o\`u les quantit\'es $\hat{N}(y;T)$, $\hat{N}_0(y;T)$ et $\hat{N}_0'(y;T)$ sont les analogues de
 $\hat{N}(T)$, $\hat{N}_0(T)$ et $\hat{N}_0'(T)$, respectivement, en fonction du param\`etre $y$.
 Dans cette derni\`ere inegalit\'e double,
 le membre \`a gauche est un nombre entier pair, et celui \`a droite un entier impair, donc
$$ \hat{N}(y;T) - \hat{N}_0(y;T) \leq 2 N_{h,\sigma> a}.$$
Pour $T$ fix\'e, le th\'eor\`eme de Hurwitz nous donne, pour $y >1$ assez proche de $1$,
$$ \hat{N}(T) - \hat{N}_0(T) \leq \hat{N}(y;T) - \hat{N}_0(y;T) \leq 2 N_{h,\sigma> a}.$$
Cela finit la preuve. \hfill{$\blacksquare$}\\


Le corollaire \ref{debrujingn} peut \^etre appliqu\'e \`a un polyn\^ome complexe quelconque;
on peut comparer ce r\'esultat avec ceux de de Bruijn, cfr. \cite[Lemma 2, Theorem 9A]{debruijn}.

On se demande s'il est possible d'am\'eliorer l'estimation obtenue dans le th\'eor\`eme \ref{teoalin}
pour avoir un r\'esultat aussi pr\'ecis que dans le cas r\'eel.
Pour r\'eduire d'une unit\'e l'estimation (\ref{naac}) dans un cas non trivial
(non couvert par la proposition \ref{copiado}),
le premier pas est de remarquer que
l'estimation (\ref{mejorable}) reste valable en rempla\c{c}ant le membre \`a droite par
$$
\tilde{N}_0'(T) \geq \ment{\mfrac{1}{\pi}\arg h(a+iT)} + \ment{-\mfrac{1}{\pi}\arg h(a-iT)} -1.
$$
Soit $\psi(\tau) = \mfrac{1}{\pi}\arg h(a+iT)$. Dans la d\'emonstration du th\'eor\`eme \ref{teoalin},
l'estimation de $R(T)-N_0'(T)$ fait apparaitre $\psi(\tau)-\ment{\psi(\tau)}$, ce que normalement
on borne trivialement par $0$. On peut garder cet \'ecart; sa contribution \`a l'estimation de la borne (\ref{naac})
est la moyenne ${\frac{1}{T}\int_{0}^T \bigl( \psi(\tau)-\ment{\psi(\tau)}\bigr) d\tau}$.
On r\'e\'enonce cela en termes de la fonction partie fractionnaire $\{ x \} = x - \lfloor x\rfloor$.
On sait que $\{ x \} = 1 + x - \ment{x}$ pour $x \notin \mathbb{Z}$, donc
$\{ \psi(\tau)\} =1 +\psi(\tau) - \ment{\psi(\tau)}$ presque partout (les z\'eros de $h(s)$ \'etant isol\'es), et l'estimation triviale
est $\{ \psi(\tau)\} \leq 1$.

{\bf Question.} Sous quelles conditions sur $h(s)$, la fonction $\psi(\tau)$ satisfait-elle
$$ 0 \leq \ds \limsup_{T \to +\infty} \Frac{1}{T} \int_0^T \{ \psi(\tau) \} d\tau < 1 ? $$

On fournit deux conditions suffisantes sur $\psi(\tau)$:
\begin{enumeratei}
\item \label{psiderivada} Soit $\psi(\tau)$ diff\'erentiable, $c_1, c_2>0$ tels que
$$ 0 < c_1 \leq \psi'(\tau) \leq c_2 $$
pour $\tau$ suffisament grand, alors
$$ \frac{c_1}{2c_2} \leq \liminf_{T\to +\infty} \Frac{1}{T} \int_0^T \{\psi(\tau)\} d\tau \leq
\limsup_{T\to +\infty} \Frac{1}{T} \int_0^T \{\psi(\tau)\} d\tau \leq 1 -
\frac{c_1}{2c_2}.
$$
\item \label{psiconvexa} Soit $\psi(\tau)$ convexe et croissante pour $\tau$ suffisament grand, alors
$$ \limsup_{T\to +\infty} \Frac{1}{T} \int_0^T \{\psi(\tau)\} d\tau \leq \frac{1}{2}.$$
\end{enumeratei}

La condition \itemref{psiderivada} peut \^etre appliqu\'ee pour \'etablir une extension de la proposition \ref{perturbacionestable}:
la condition est satisfaite par un polyn\^ome complexe $p(s)$ quelconque, lorsque $y > 1$.

La condition \itemref{psiconvexa} permet d'\'etendre les r\'esultats du \S \ref{lweng},
en particulier le th\'eor\`eme \ref{pabierto}, au cas o\`u $p(s)$ est un polyn\^ome complexe.
Les propri\'et\'es de croissance et convexit\'e (dans un petit intervalle) de la fonction $\psi(\tau)=\frac{1}{\pi}\arg \bigl((i\tau)\zeta^*(1+2i \tau)y^{i\tau}\bigr)$ (le cas $p(s)=1$) peuvent \^etre
vues dans \cite[Lemma 2.3]{hki4}; en particulier $\psi'(\tau) \sim \mfrac{1}{\pi} \log \tau $
et la condition \itemref{psiderivada} n'est pas satisfaite.

\subsection{Dominance faible et densit\'e}

Les th\'eor\`emes \ref{teoalindensidad} et \ref{teoalindom} se g\'en\'eralisent de fa\c{c}on \'evidente; pour les g\'en\'eraliser on fournira au m\^eme temps une combinaison des deux.

On introduit quelques notations. Comme dans le \S \ref{inficritica}, pour une fonction $h(s)$ et $\sigma_0<\sigma_1$, $T>0$, on note
$$ \tilde{N}_h(\sigma_0,\sigma_1,T) = \# \bigl\{ s \in \ce\tq  h(s)=0,  \, \sigma_0 \leq \sigma \leq \sigma_1, |\tau| < T \bigr\}, $$
et $ \tilde{N}_h(\sigma_0,T) = \tilde{N}_h(\sigma_0,+\infty,T) $.

De m\^eme, comme dans le \S \ref{faibledd}, pour la fonction $f(s) = h(s) \pm h(2a-s)$, on d\'efinit
$$\hat{N}(T,\sigma_0) = \# \bigl\{ s \in \ce\tq f(s)=0, |\sigma-a| < \sigma_0, |\tau| < T \bigr\}. $$

En utilisant une extension simple du lemme \ref{zerospormontones},
on d\'emontre le r\'esultat suivant, la version combin\'ee des th\'eor\`emes  \ref{teoalindensidad}
et \ref{teoalindom}.

\begin{teor} \label{teoalindomdensidadc} Soit $\sigma_0 > a$, $h(s)$ une fonction m\'eromorphe,
avec un nombre fini de p\^oles,
holomorphe et non-nulle sur la droite $\sigma=\sigma_0$. On pose
$$ f(s) = f^{\pm}(s) = h(s) \pm \ovl{h}(2a-s). $$
On suppose que la fonction
$$ F(s) = \Frac{\ovl{h}(2a-s)}{h(s)} $$
satisfait
\begin{enumerate}
\item[\hypertarget{princcc1xp}(i')] $F(s) \neq \pm 1$ sur la droite $\sigma=\sigma_0$, et il existe $\tau_0 >0$ tel que $|F(s)|<1$ pour $\sigma=\sigma_0$, $|\tau| \geq \tau_0$;
\item[\hypertarget{princcc2xp}(ii')] il existe une fonction croissante $\phi: \re \to \re$, une constante $K>0$ et des suites $(T_n)_n$, $(T_n^*)_n$ telles que $\ds\lim_{n \to \infty} T_n = \ds\lim_{n \to \infty} T_n^* = +\infty$,
$$ T_n \leq T_{n+1} \leq \phi(T_n), \quad T_n^* \leq T_{n+1}^* \leq \phi(T_n^*), \quad \mbox{ pour } n \geq 1, $$
     et $\bigl| F(s) \bigr| < e^{K|s|}$
pour $a \leq \sigma \leq \sigma_0$, et $\tau = T_n$ ou $\tau = -T_n^*$, $n \geq 1$.
\end{enumerate}
Alors, pour $T>0$
$$ \hat{N}(T,\sigma_0) - \hat{N}_0'(T) \leq 4 \hat{N}_h\bigl(a,\sigma_0,\phi(2T)\bigr) + O(1). $$
\end{teor}

\subsection{Translat\'ees des fonctions $L$ de Dirichlet}

Soit $L(s,\chi)$ une s\'erie de Dirichlet associ\'ee \`a un caract\`ere de Dirichlet
primitif de conducteur $N>1$,
$$ L(s,\chi)= \summ_{n=1}^\infty \Frac{\chi(n)}{n^s}$$
pour $\sigma > 1$.
On associe \`a cette s\'erie la fonction complet\'ee $\xi(s,\chi)$, d\'efinie par
$$ \xi(s,\chi) = \Bigl( \frac{N}{\pi} \Bigr)^{\frac{s}{2}} \Gamma\Bigl( \frac{s+\kappa}{2} \Bigr) L(s,\chi), $$
o\`u
$$ \kappa  = \frac{1}{2} \bigl(1-\chi(-1) \bigr) = \left\{
\begin{array}{cl}
0,  & \mbox{ si } \chi(-1) = 1,\\
1,  & \mbox{ si } \chi(-1) = -1.
\end{array}
\right.
$$
La fonction $\xi(s,\chi)$ est une fonction enti\`ere de genre $1$ qui satisfait \`a l'\'equation fonctionnelle \cite[Theorem 4.15]{iwaniec}
$$ \xi(s,\chi) = \epsilon(\chi) \xi(1-s,\ovl{\chi}), $$
o\`u $\epsilon(\chi) = i^{-\kappa} \frac{\tau(\chi)}{\sqrt{N}}$, $\tau(\chi)$ est une somme de Gauss, de
sorte que $\bigl|\epsilon(\chi)\bigr|=1$.
Les z\'eros de $\xi(s,\chi)$ sont sur la bande $0 < \sigma < 1$. L'hypoth\`ese de Riemann pour $L(s,\chi)$
est la conjecture que tous les z\'eros de $L(s,\chi)$ dans la bande $0 < \sigma < 1$ sont sur la droite $\sigma = 1/2$.

On va r\'eproduire les r\'esultats du \S \ref{aquiestimation} et du \S \ref{infidensidad} pour les fonctions $L(s,\chi)$. On a besoin de l'analogue du lemme \ref{riemannpropi} pour la fonction $L(s,\chi)$.

\begin{lema} \label{listosff} Soit $L(s,\chi)$ la s\'erie de Dirichlet associ\'ee \`a un caract\`ere primitif de conducteur $N>1$.
\begin{enumerate}
\item \label{silva} Il existe une constante $A>0$ telle que, pour chaque $\eps >0$, ${1}/{2} < \sigma \leq 1$, $T \geq 3$
$$ \hat{N}_{L(\cdot,\chi)}(\sigma,T)= O \bigl( (T^{(2+\eps)(1-\sigma)} + T^{c(\sigma)(1-\sigma)}) \log^A T \bigr), $$
o\`u  $ c(\sigma)= \min \Bigl( \frac{3}{2-\sigma}, \frac{3}{3\sigma-1} \Bigr) $.
\item Il existe une constante $A>0$ telle que pour chaque $n \geq 1$, il existe $n<T_n<n+1$ tel que
$$ \bigl| L(\chi,s) \bigr|  > \tau^{-A} $$
lorsque $\tau = T_n$, $-1 \leq \sigma \leq 2$.
\end{enumerate}
\end{lema}

{\bf D\'emonstration.}
\begin{enumerate}
\item C'est un cas particulier du th\'eor\`eme de densit\'e dans \cite[\S 10.4, p. 260]{iwaniec}.
\item On consid\`ere la formule dans \cite[p. 102]{davenport}
$$ \Frac{L'(\chi,s)}{L(\chi,s)} = \summ_{\rho, |\tau-\gamma| \leq 1} \Frac{1}{s-\rho} + O\bigl(\log|\tau|\bigr) $$
valable pour $-1 \leq \sigma \leq 2$, $|\tau| \geq 2$, la somme portant
sur les z\'eros $\rho=\beta+i \gamma$ de $L(s,\chi)$ avec $|\tau-\gamma| \leq 1$. L'int\'egration de cette formule
entre $s=\sigma+i\tau$ et $2+i\tau$ nous donne
$$ \log L(\chi,s)= \summ_{\rho, |\tau-\gamma| \leq 1} \log(s-\rho) + O\bigl(\log |\tau| \bigr), $$
pour $ -1 \leq \sigma \leq 2$, $|\tau| \geq 2$.
\`A partir d'ici, la m\'ethode de d\'emonstration du th\'eor\`eme \'equivalent pour la fonction z\^eta
de Riemann, dans \cite[\S 9.7]{ztitchmarsh}, nous donne le r\'esultat. {\lqqd}
\end{enumerate}

On peut maintenant \'enoncer le r\'esultat attendu. On remarque
que le cas {\it(\ref{silva})} avec un polyn\^ome constant $p(s)$  \`a \'et\'e \'etabli par Lagarias \cite[Lemma 5.1, Theorem 5.1]{lagarias},
qui au m\^eme temps a determin\'e la distribution limite des espacements des z\'eros ces fonctions.

\begin{teor} \label{compornada} Soit $\alpha > 0$. Soit $\xi(s,\chi)$ une fonction $L$ compl\'et\'ee associ\'ee
\`a un caract\`ere de Dirichlet $\chi$ primitif de conducteur $N >1$, $p(s)$ un polyn\^ome complexe sans
z\'ero sur la droite $\sigma = 1/2$. On consid\`ere la fonction
$$ f_{p, \chi,\alpha}^{\pm}(s) = p(s) \xi(s+\alpha,\chi) \pm \ovl{p}(1-s) \xi(s-\alpha,\chi).$$
Alors
\begin{enumerate}
\item \label{finaluno} Pour $\alpha \geq 1/2$,
$$ \hat{N}(T)-\hat{N}_0'(T) \leq  2 N_{p,\sigma>1/2}+1, $$
o\`u $N_{p,\sigma>1/2}$ est le nombre de racines de $p(s)$ dans le demi-plan $\sigma > 1/2$.
Donc tous les z\'eros de $f_{p,\chi,\alpha}^{\pm}(s)$, sauf un nombre fini, sont sur la droite $\sigma=1/2$, ils sont
simples et entrelac\'es. De m\^eme pour $0 < \alpha < 1/2$ sous l'hypoth\`ese de Riemann pour $L(s,\chi)$.

Dans le cas particulier o\`u $N_{p,\sigma>1/2} = 0$, tous les z\'eros des fonctions
$f_{p, \chi,\alpha}^{\pm}(s)$ sont sur la droite $\sigma=1/2$, ils sont simples et entrelac\'es.

\item \label{finaldos} Pour $0 < \alpha < 1/2$,
$$ \hat{N}(T)-\hat{N}_0'(T) \leq  4 \hat{N}_{L(\cdot,\chi)}\bigl(\mfrac{1}{2}+\alpha,2T+4)+O(1). $$
En particulier, il existe une constante $A>0$ telle que, pour tout $\eps > 0$,
$$ \hat{N}(T)-\hat{N}_0'(T)= O \bigl( (T^{(2+\eps)(\frac{1}{2}-\alpha)} + T^{k(\alpha)(\frac{1}{2}-\alpha)}) \log^A T \bigr), $$
pour $T > 0$, o\`u $ k(\alpha)= \min \Bigl( \frac{6}{3-2\alpha}, \frac{6}{1+6\alpha} \Bigr)$.
\end{enumerate}
\end{teor}

{\bf D\'emonstration.} Si dans l'\'equation fonctionnelle pour $L(s,\chi)$ on \'ecrit $\epsilon(\chi) = e^{2i \theta}$, alors
$$ e^{-i\theta} f_{\chi,\theta,\alpha}^{\pm}(s) = e^{-i\theta} p(s) \xi(s+\alpha,\chi) \pm e^{i\theta} \ovl{p}(1-s) \xi(1-s+\alpha,\ovl{\chi}).$$
On v\'erifiera les conditions des th\'eor\`emes \ref{teoalinc} et \ref{teoalindomdensidadc}.
On pose
$$ F(s) = e^{2i\theta} \Frac{\ovl{p}(1-s)}{p(s)} \cdot \Frac{\xi(s-\alpha,\chi)}{\xi(s+\alpha,\chi)}. $$
Par la formule de Stirling, pour $\sigma \geq 1/2$, $|s|$ suffisament grand
$$ F(s)
=  \omega \Bigl(\frac{\pi}{N} \Bigr)^{\alpha} s^{-\alpha} \Bigl( 1 + O\bigl( |s|^{-1} \bigr) \Bigr) \Frac{L(s-\alpha,\chi)}{L(s+\alpha,\chi)},
$$
o\`u $| \omega| = 1$.
Pour la s\'erie de Dirichlet $L(s,\chi) = 1 + O(2^{-\sigma})$ pour $\sigma$ suffisament grand, puis
$$ \bigl| F(s) \bigr| = O(\sigma^{-\alpha}), $$
donc la condition {\it(\hyperlink{princcc1x}{i})} du th\'eor\`eme \ref{teoalinc}.

D'un autre c\^ot\'e, la fonction $L(s,\chi)$ est polynomialment born\'ee dans un demi-plan
$\sigma \geq \sigma_\alpha = 1/2-\alpha$ \cite[Lemma 5.2]{iwaniec},
c'est-\`a-dire $L(s,\chi) = O(|\tau|^{A_0})$ avec $A_0=A_0(\sigma_\alpha)$. Le lemme \ref{listosff} nous donne
une constante $A_1> 0$ et une suite $(T_n)_n$ telles que $n < T_n < n+1$, $\bigl| L(s,\chi) \bigr|^{-1} < |\tau|^{A_1}$
pour $-1 \leq \sigma \leq 2$, $\tau = T_n$, $n \geq 1$. En appliquant le m\^eme lemme
\`a la fonction $L(s,\ovl{\chi})$, on obtient $A_2>0$ et une suite $(T_n^*)_n$ tels que $n < T_n^* < n+1$,
et $\bigl| L(s,\chi) \bigr|^{-1}= \bigl| L(\ovl{s},\ovl{\chi}) \bigr|^{-1} < |\tau|^{A_1}$
pour $-1 \leq \sigma \leq 2$, $\tau = -T_n^*$, $n \geq 1$.
Le rassemblement de ces conditions nous donne
$$ \bigl| F(s) \bigr| = O\bigl(|\tau|^B\bigr), $$
pour $\sigma \geq 1/2$, $\tau= T_n$ ou $\tau=-T_n^*$, $n \geq 1$,
d'o\`u la condition {\it(\hyperlink{princcc2x}{ii})} du th\'eor\`eme \ref{teoalinc}, ainsi que la condition
{\it(\hyperlink{princcc2xp}{ii'})} du th\'eor\`eme \ref{teoalindomdensidadc} pour $\sigma_0>2$ quelconque et $\phi(T)=T+2$,
comme dans la proposition \ref{miomio}.
Le premier th\'eor\`eme (avec la remarque (\ref{phragmen}.\hyperlink{phragmen1}{1}) au th\'eor\`eme \ref{teoalin} et
 la proposition \ref{copiado} dans le cas o\`u $N_{p,\sigma>1/2}=0$) nous donne (\ref{finaluno}), le deuxi\`eme (\ref{finaldos}), et l'estimation
dans (\ref{finaldos}) vient directement de la premi\`ere partie du lemme \ref{listosff}. \lqqd

\appendix
\section{Stabilit\'e} \label{seccionestabilidad}

Par souci de compl\'etude, on va mentionner des r\'esultats concernant le probl\`eme de la stabilit\'e d'un point
de vue plus classique. Soit $h(s)$ une fonction enti\`ere et $a \in \re$;
on dira que $h(s)$ est {\it stable par rapport \`a la droite $\sigma=a$}, si tous ses z\'eros
sont dans le demi-plan $\sigma < a$ (et tout simplement {\it stable} si $a=0$).

On r\'e\'enonce un r\'esultat originalement d\^u \`a de Branges \cite[Lemma 5]{debranges},
puis retrouv\'e par Lagarias dans {\cite[Lemma 2.2]{lagarias}.

\begin{prop} \label{copiado} Soit $h(s)$ une fonction enti\`ere stable par rapport \`a
la droite $\sigma = a$, et telle que la fonction
$$ F(s) = \frac{\ovl{h}(2a-s)}{h(s)} $$
satisfait
$$ |F(s)| < 1 \mbox{ pour } \sigma > a. $$
Alors
\begin{enumerate}
\item \label{chucuno} les fonctions $ f(s) = f^{\pm}(s) = h(s) \pm \ovl{h}(2a-s) $ ont tous leurs z\'eros sur la droite $\sigma=a$;
\item \label{chucdos} toute fonction de phase $\varphi(\tau)=\arg h(a+i\tau)$,
est une fonction strictement croissante. Donc, les z\'eros des fonctions $f^{\pm}(s)$ sont simples et entrelac\'es.
\end{enumerate}
\end{prop}

{\bf D\'emonstration.}
En fait, $\sigma > a$, $ |f(s)| \geq |h(s)| \bigl( 1 - |F(s)|) > 0 $, puis
par l'\'equation fonctionnelle, pour $\sigma < a$, on a $2a-\sigma > a$ et
$ |f(s)| = | \ovl{f}(2a-s) | = | f(2a-\ovl{s})| > 0$.
Cela montre {\it (\ref{chucuno})}.
On renvoit le lecteur \`a la d\'emonstration \cite[Lemma 2.2]{lagarias} pour la partie {\it (\ref{chucdos})}.
On avait dej\`a remarqu\'e l'entrelacement des z\'eros dans ce cas, au \S \ref{seccerosrecta}. \lqqd\\

Un argument simple montre qu'un polyn\^ome stable non constant (cfr. \ref{policaso}) satisfait les conditions de la proposition \ref{copiado}.
Soit $ p(s) = a_0 \prodd_\rho (s-\rho)$ la factorisation de $p(s)$, et $ F(s) = {\ovl{p}(-s)}/{p(s)} $. On a
$$
\bigl|F(s) \bigr|= \left| \Frac{\ovl{p}(-s)}{p(s)} \right| =
\prodd_\rho \left| \Frac{-s-\ovl{\rho}}{s-\rho} \right| =
\prodd_\rho \left| \Frac{s+\ovl{\rho}}{s-\rho} \right|.
$$
Soit $\rho=\beta+i \gamma$ un z\'ero de $p(s)$; on a $\beta < 0$. Pour $s=\sigma+i \tau$,
l'in\'egalit\'e $ |s + \ovl{\rho}| < |s-\rho| $, ou
$$ (\sigma+\beta)^2 + (\tau-\gamma)^2 < (\sigma-\beta)^2 + (\tau-\gamma)^2 $$
est \'equivalente \`a $ 4 \sigma \beta < 0 $. Puisque $\beta < 0$,  on doit avoir $\sigma > 0$.
Cela montre que l'in\'egalit\'e $\bigl|F(s)\bigr| <1$ est satisfaite pour chaque terme du produit.

Par la remarque (\ref{phragmen}.\hyperlink{phragmen1}{1}) faite au th\'eor\`eme \ref{teoalin}, valable pour une fonction sans la sym\'etrie r\'eelle,
si $h(s)$ est une fonction stable par rapport \`a la droite $\sigma=a$ qui satisfait les conditions du th\'eor\`eme
\ref{teoalinc}, alors elle satisfait les conditions de la proposition \ref{copiado}, et on obtient les m\^emes
conclusions dans les deux th\'eor\`emes pour les fonctions $f^{\pm}(s)$, voire mieux
avec la proposition \ref{copiado}.
Si la fonction $h(s)$ n'est pas stable, la remarque (\ref{phragmen}.\hyperlink{phragmen2}{2}) au th\'eor\`eme
nous donne une condition avec laquelle on ne peut rien
dire sur les z\'eros des fonctions $f^{\pm}(s)$; il est l\`a o\`u le th\'eor\`eme
\ref{teoalinc} donne une r\'eponse.

On veut \'etendre la condition de la proposition \ref{copiado} \`a des fonctions enti\`eres. La proc\'edure
men\'ee pour les polyn\^omes est aussi valable pour les fonctions de genre $0$ \cite[Part I, \S 4.2]{levin}, qui se factorisent de la m\^eme fa\c{c}on
que les polyn\^omes. Mais les fonctions qu'on rencontre dans les applications, telles que la fonction
$\xi(s)$ de Riemann ou les fonctions $L$, sont de genre $1$, et poss\`edent une factorisation moins simple.

Le r\'esultat suivant, est \'etabli par Suzuki et Lagarias dans \cite[Theorem 4, \S 2]{lagarias2}
pour des fonctions r\'eelles sur la droite r\'eelle. Il peut \^etre compar\'e au r\'esultat classique de P\'olya
\cite[Hilfssatz II]{polya}.

\begin{teor} \label{perturbationalpha}
Soit $h(s)$  une fonction enti\`ere de genre $0$ ou $1$, qui satisfait une \'equation fonctionnelle de la forme
\begin{equation} \label{simetriia}
\ovl{h}(2a-s) = e^{i \theta} h(s)
\end{equation}
pour $0 \leq \theta < 2\pi$, et telle que ses z\'eros sont inclus dans une bande de la forme
$$|\sigma-a| < b,$$
o\`u $b> 0$, et $\alpha \geq b$. Alors
$$ \left| \frac{h(s-\alpha)}{h(s+\alpha)} \right|  <  1 \mbox{ pour } \sigma > a. $$
\end{teor}

{\bf D\'emonstration.} Sans perte de g\'en\'eralit\'e, on peut traiter le cas $a=0$.
On supposera initialement que $h(0) \neq 0$.
On consid\`ere la factorisation de Weierstrass de $h(s)$ \cite[Part I, \S 4.2]{levin}(le choix $a=0$ simplifie le calcul;
en g\'en\'eral on prendrait la factorisation de $h(a+s)$ au lieu de celle de $h(s)$)
\begin{equation} \label{prodweier}
h(s) = e^{A+Bs} \prod_\rho \Bigl( 1-\frac{s}{\rho} \Bigr) e^{s/\rho},
\end{equation}
o\`u $\summ_\rho \frac{1}{|\rho|^2} < +\infty$, $|\rho| < b$ pour tout $\rho$.
L'\'equation fonctionnelle $\ovl{h}(-s) = e^{i \theta} h(s)$ entra\^ine que si $h(\rho) =0$, alors
$h(-\ovl{\rho}) = e^{-i\theta} \ovl{h(\rho)} = 0$.
Puis, on peut r\'egrouper les z\'eros $\rho=\beta+i\gamma$ par des blocs $B(\rho) = \{\rho, -\ovl{\rho}\}$ pour $\beta \neq 0$,
$B(\rho) = \{\rho\}$ pour $\beta = 0$.
Cela nous donne
$$ h(s) = e^{A+Bs} \prod_{\substack{B(\rho) \\ \beta \neq 0}}
\Bigl( 1-\frac{s}{\rho} \Bigr) \Bigl( 1+\frac{s}{\ovl{\rho}} \Bigr) e^{s(1/\rho-1/\ovl{\rho})}
\prod_{\substack{\rho \\ \beta = 0}}
\Bigl( 1-\frac{s}{\rho} \Bigr) e^{s(1/\rho)}
$$
et $1/\rho - 1/\ovl{\rho} = -2i \gamma/|\rho|^2$. Donc, on peut redistribuer les constantes et \'ecrire
$$ h(s) =
e^{A+Bs} \prod_\rho \Bigl( 1-\frac{s}{\rho} \Bigr) e^{c(\rho)s}, $$
o\`u $c(\rho)= - i \gamma/|\rho|^2$ (ce qui est important ici est le fait que $\Re\bigl( c(\rho) \bigr)=0$), et le produit (ainsi que tous les produits dor\'enavant) est conditionnellement convergent,
$$ \prod_\rho \Bigl( 1-\frac{s}{\rho} \Bigr) e^{c(\rho)s}=
\lim_{T \to \infty} \prod_{|\rho|<T} \Bigl( 1-\frac{s}{\rho} \Bigr) e^{c(\rho)s}$$

La deriv\'ee logarithmique de $h(s)$ (de (\ref{prodweier})) \'evalu\'ee en $s=0$ est
$ {h'(0)}/{h(0)} = B $, et par l'\'equation fonctionnelle
$$ B = \Frac{h'(0)}{h(0)} = -\Frac{h'(-0)}{h(-0)} = -B, $$
d'o\`u $\Re(B) = 0$.

Puis, pour $\alpha \geq b$, $s=\sigma+i \tau$
$$
\left| \frac{h(s-\alpha)}{h(s+\alpha)} \right|  =
\prod_\rho
\left| \frac{1-\frac{s-\alpha}{\rho}}{1-\frac{s+\alpha}{\rho}} \right| =
\prod_\rho
\left| \frac{\rho-s+\alpha}{\rho-s-\alpha} \right|
$$
et encore par la sym\'etrie des z\'eros
$$
\left| \frac{h(s-\alpha)}{h(s+\alpha)} \right|  =
\prod_\rho
\left| \frac{\rho-s+\alpha}{-\ovl{\rho}-s-\alpha} \right|. $$
On montrera que chacun des termes du produit
satisfait l'in\'egalit\'e attendue. L'in\'egalit\'e
\begin{equation} \label{gthneq}
|\rho- s - \alpha| < | -\ovl{\rho}-s+\alpha|
\end{equation}
et \'equivalente \`a
$$ (\beta-\sigma-\alpha)^2 + (\gamma-\tau)^2 < (-\beta-\sigma+\alpha)^2+(\gamma-\tau)^2 $$
ou $ (\beta-\alpha)\sigma < 0 $. Comme $\beta-\alpha < b-\alpha \leq 0$, alors on a (\ref{gthneq}) si et seulement si  $\sigma > 0$.
Ceci nous donne la conclusion du th\'eor\`eme dans ce cas.

Si $s=0$ est un z\'ero d'ordre $m$ de $h(s)$, on peut \'ecrire $h(s)=s^m h_1(s)$, o\`u $h_1(0) \neq 0$,
$\ovl{h_1}(-s)= e^{i \theta_1} h(s)$ ($\theta_1 = \theta$ si $m$ est pair, $\theta_1=\theta \pm \pi$ si $m$ est impaire) et on peut obtenir les in\'egalit\'es attendues pour $h_1(s)$. De plus
$$ \left| \frac{s-\alpha}{s+\alpha} \right|  <  1 \mbox{ pour } \sigma > 0, $$
ce qui avec l'in\'egalit\'e pour $h_1(s)$ nous donne le r\'esultat. \lqqd\\

En vue de l'\'equation fonctionnelle, le th\'eor\`eme \ref{perturbationalpha} \'etablit que la
fonction $h(s+\alpha)$ v\'erifie les conditions de la proposition \ref{copiado}.
Maintenant on d\'ecrit rapidement certaines applications de combiner ces r\'esultats.
On garde la notation du th\'eor\`eme \ref{perturbationalpha} pour $a=\frac{1}{2}$, $\alpha$ et $b$.
\begin{enumerate}
\item On pose $F(s) = \xi(2s - \frac{1}{2})$, ici $a=\frac{1}{2}$, $b=\alpha=\frac{1}{4}$.
La fonction $\xi(2s)=F\bigl(s+\frac{1}{4}\bigr)$ satisfait les conditions de la proposition \ref{copiado},
et on obtient le corollaire \ref{riemannzetawengdos}, l'hypoth\`ese de Riemann pour la fonction zeta de Weng de rang $2$.
De m\^eme on peut obtenir la partie (ii) du th\'eor\`eme \ref{pabierto}.
Il est pour cet objectif que le th\'eor\`eme \ref{perturbationalpha},
appliqu\'e \`a une fonction r\'eelle sur la droite r\'eelle, a \'et\'e introduit par Lagarias et Suzuki en
\cite{lagarias2}.
\item La fonction $\xi(s+\alpha)$ satisfait les conditions de la proposition \ref{copiado}
inconditionellement pour $\alpha \geq \frac{1}{2}$ ($a=\frac{1}{2}$, $b=\frac{1}{2}$), et sous l'hypoth\`ese de Riemann
pour $\alpha > 0$ (dans ce cas $b$ est quelconque). On obtient les r\'esultats sur les fonctions $\xi(s+\alpha)\pm \xi(s-\alpha)$
dont on a discut\'e \`a la fin du \S \ref{aquiestimation}.
\item Les fonctions $\xi(s+\alpha,\chi)$, associ\'ees aux charact\`eres primitifs de conducteur $N>1$, satisfont les conditions de la proposition \ref{copiado}
inconditionellement pour $\alpha \geq \frac{1}{2}$, et sous l'hypoth\`ese de Riemann pour
$L(s,\chi)$ pour $\alpha > 0$. Ces fonctions ont \'et\'e aussi consider\'ees par Lagarias dans \cite{lagarias}.
La combinaison de ces fonctions avec des polyn\^omes dont les racines sont dans le demi-plan $\sigma < 1/2$
nous donne le cas particulier dans la premi\`ere partie du th\'eor\`eme \ref{compornada}.
\item La fonction de Bessel $K_s(A)$ ($a=0$). Dans l'article \cite{polya}, P\'olya localise d'abord les z\'eros de cette fonction
\`a l'aide de la formule \ref{vvpolya}, dans la bande $|\sigma| < 1=b$, pour montrer apr\`es
qu'ils sont simples et align\'es avec les relations (\ref{eqfunck}).
\end{enumerate}

\section*{Remerciements}
Je remercie Michel Balazard, directeur de ma th\`ese de doctorat, dont ce travail fait partie, pour les conseils
opportuns,
et Arnaud Chadozeau, pour des discussions concernant le \S \ref{sanssymetrie}.

\bibliography{biblio}
\bibliographystyle{abbrv-fr}

\end{document}